\documentclass[11pt]{amsart}

\usepackage{amsmath}
\usepackage{amssymb}
\usepackage{dsfont}
\usepackage[all]{xy}

\usepackage{empheq} % for \xhookrightarrow

\usepackage{enumerate}

\usepackage{mathrsfs}

\usepackage{verbatim}

\usepackage[retainorgcmds]{IEEEtrantools} 

\usepackage{tensor}

\input xy              %this packet is used for creating diagrams.  
\xyoption{all}

\DeclareMathOperator{\KZ}{KZ}
\DeclareMathOperator{\Hom}{Hom}
\DeclareMathOperator{\End}{End}

\theoremstyle{plain}   \newtheorem{theorem00.1}{Proposition}[section]
\theoremstyle{plain}   \newtheorem{theorem00.2}[theorem00.1]{Proposition}
\theoremstyle{plain}   \newtheorem{theorem00.3}[theorem00.1]{Lemma}
\theoremstyle{plain}   \newtheorem{theorem00.4}[theorem00.1]{Proposition}
\theoremstyle{plain}   \newtheorem{theorem00.5}[theorem00.1]{Corollary}
\theoremstyle{plain}   \newtheorem{theorem00.6}[theorem00.1]{Proposition}
\theoremstyle{definition}   \newtheorem{definition00.1}{Definition}[section]
\theoremstyle{remark}   \newtheorem{remark00.1}{Remark}[section]
\theoremstyle{plain}   \newtheorem{theorem00.7}[theorem00.1]{Proposition}
\theoremstyle{plain}   \newtheorem{theorem00.8}[theorem00.1]{Proposition}
\theoremstyle{plain}   \newtheorem*{theorem00.9*}{Conjecture}
\theoremstyle{plain}   \newtheorem{theorem00.10}[theorem00.1]{Proposition}
\theoremstyle{plain}   \newtheorem{theorem00.11}[theorem00.1]{Corollary}

\theoremstyle{plain}   \newtheorem*{theorem000.1}{Theorem A}
\theoremstyle{plain}   \newtheorem*{theorem000.2}{Theorem B}
\theoremstyle{plain}   \newtheorem*{theorem000.3}{Theorem C}
\theoremstyle{plain}   \newtheorem*{theorem000.4}{Theorem D}

\theoremstyle{plain}   \newtheorem{theorem3.12}{Proposition}[section]
\theoremstyle{plain}   \newtheorem{theorem3.6}[theorem3.12] {Proposition}
\theoremstyle{plain}   \newtheorem{theorem3.14}[theorem3.12] {Proposition}
\theoremstyle{plain}   \newtheorem{theorem3.13}[theorem3.12] {Proposition}

\theoremstyle{plain}   \newtheorem{theorem5.17}{Proposition}[section]
\theoremstyle{plain}   \newtheorem{theorem5.9}[theorem5.17]{Proposition}
\theoremstyle{plain}   \newtheorem{theorem5.15}[theorem5.17]{Proposition} 
\theoremstyle{plain}   \newtheorem{theorem5.6}[theorem5.17]{Proposition} 
\theoremstyle{plain}   \newtheorem{theorem5.61}[theorem5.17]{Proposition} 
\theoremstyle{plain}   \newtheorem{theorem5.1}[theorem5.17]{Proposition} 
\theoremstyle{definition}   \newtheorem{definition5.4}{Definition}[section]
\theoremstyle{definition}   \newtheorem{definition5.3}[definition5.4]{Definition}
\theoremstyle{plain}   \newtheorem{theorem5.8}[theorem5.17]{Proposition} 
\theoremstyle{plain}   \newtheorem{theorem5.10}[theorem5.17]{Proposition}
\theoremstyle{definition}   \newtheorem{definition5.5}[definition5.4]{Definition}
\theoremstyle{plain}   \newtheorem{theorem5.11}[theorem5.17]{Proposition}
\theoremstyle{plain}   \newtheorem{theorem5.12}[theorem5.17]{Corollary}
\theoremstyle{plain}   \newtheorem{theorem5.2}[theorem5.17]{Proposition} 
\theoremstyle{plain}   \newtheorem{theorem5.18}[theorem5.17]{Proposition}
\theoremstyle{definition}   \newtheorem{definition5.7}[definition5.4]{Definition}
\theoremstyle{plain}   \newtheorem{theorem5.19}[theorem5.17]{Proposition}

\theoremstyle{plain}   \newtheorem{theorem7.1}{Proposition}[section] 
\theoremstyle{plain}   \newtheorem{theorem7.2}[theorem7.1]{Proposition} 
\theoremstyle{plain}   \newtheorem{theorem7.3}[theorem7.1]{Proposition} 
\theoremstyle{plain}   \newtheorem{theorem7.4}[theorem7.1]{Proposition} 
\theoremstyle{plain}   \newtheorem{theorem7.5}[theorem7.1]{Proposition}
\theoremstyle{plain}   \newtheorem{theorem7.6}[theorem7.1]{Corollary}

\theoremstyle{plain}   \newtheorem{theorem6.05}{Lemma}[section]
\theoremstyle{plain}   \newtheorem{theorem6.01}[theorem6.05]{Proposition}
\theoremstyle{plain}   \newtheorem{theorem6.02}[theorem6.05]{Proposition} 
\theoremstyle{plain}   \newtheorem{theorem6.03}[theorem6.05]{Proposition} 
\theoremstyle{definition}   \newtheorem{definition6.01}{Definition}[section]
\theoremstyle{plain}   \newtheorem{theorem6.04}[theorem6.05]{Proposition}

\title[An algebraic approach to the KZ-functor for rational Cherednik algebras]% 
 {An algebraic approach to the KZ-functor for rational Cherednik algebras associated with cyclic groups}

\author[S. Thelin]{Sam Thelin}
\address{%
   Institute of Algebra and Number Theory\\
   University of Stuttgart\\
   Pfaffenwaldring 57\\
   70569 Stuttgart\\
     Germany}
\email{sam.thelin@mathematik.uni-stuttgart.de}



\thanks{I would like to thank my supervisors Karin Erdmann and Rapha\"el Rouquier for their invaluable support. I also wish to thank Iain Gordon and Ivan Loseu for helpful comments. I am grateful to EPSRC for their financial support.}

\begin{document}

\begin{abstract}
In the case of rational Cherednik algebras associated with cyclic groups, we give an alternative proof that the projective object $P_{\text{KZ}}$ representing the KZ-functor is isomorphic to the $\Delta$-module associated with the coinvariant algebra for a subset of parameter values from which all parameter values can be obtained by integral translations. We also specify the exact parameter values for which this isomorphism occurs. Furthermore, we determine the action of the cyclotomic Hecke algebra on $P_{\text{KZ}}$  for these parameter values, thereby giving a complete algebraic description of the KZ-functor in this case.  
\end{abstract}

\maketitle

\section{Introduction} \label{introduction}

The rational Cherednik algebra $H_c(W, \mathfrak{h})$ associated with the complex reflection group $W$ and one of its reflection representations $\mathfrak{h}$, was introduced in \cite{EtGi} using ideas from \cite{Che}, and is a certain deformation of $S(\mathfrak{h} \oplus \mathfrak{h}^*) \rtimes W$, depending on a finite set of complex parameters $c$. 
The algebra $H_c(W, \mathfrak{h})$ has a so-called triangular decomposition (Theorem 1.3 of \cite{EtGi}), reminiscent of that of the universal enveloping algebra of a finite-dimensional semisimple complex Lie algebra $\mathfrak{g}$. As a result, it has a category  $\mathcal{O}_c(W, \mathfrak{h})$ of representations associated with it as introduced in \cite{DuOp} and \cite{BeEtGi}, inspired by and sharing many similarities with the BGG category $\mathcal{O}$ associated with $\mathfrak{g}$ defined in \cite{BGG}. In particular, it was proved in \cite{Gua} that category  $\mathcal{O}_c(W, \mathfrak{h})$ is a highest weight category in the sense of \cite{CPS}, indexed by the irreducible $W$-modules. An important problem in the representation theory of rational Cherednik algebras is to understand the structure of category  $\mathcal{O}_c(W, \mathfrak{h})$. 

A crucial tool to this end is the so called Knizhnik-Zamolodchikov functor, or KZ-functor for short, introduced in \cite{GGOR}. This functor relates the category  $\mathcal{O}_c(W, \mathfrak{h})$ to the category of finitely-generated modules over a corresponding cyclotomic Hecke algebra  $\mathcal{H}_c(W, \mathfrak{h})$ (cf. Section 4.C of \cite{BMR}), and enables results in the representation theory of Hecke algebras to be translated to the setting of rational Cherednik algebras and vice versa. The construction of the functor is geometric in nature, but the functor is exact, and is therefore represented by a projective object $P_{\KZ}$ of  $\mathcal{O}_c(W, \mathfrak{h})$ (cf. Section 5.4 of \cite{GGOR}).  The main aim of this article is to give an explicit algebraic description of $P_{\KZ}$, including the action of the Hecke algebra $\mathcal{H}_c(\mathbb{Z}/n, \mathbb{C})$ on it, in the case of the rational Cherednik algebras $H_c(\mathbb{Z}/n, \mathbb{C})$ corresponding to cyclic reflection groups, for a certain subset of parameter values. 

More specifically, we show that in category $\mathcal{O}_c(\mathbb{Z}/n, \mathbb{C})$, there is a subset $\mathcal{F}$ of parameter values for which  $P_{\KZ}$ is isomorphic to the module $\Delta_{c} (S(\mathfrak{h})_W)$ corresponding to the coinvariant algebra $S(\mathfrak{h})_W$, and that all possible parameter values can be obtained from $\mathcal{F}$ by integral translations. This has been conjectured by Rouquier to hold for rational Cherednik algebras associated with general complex reflection groups, and was proved in \cite{Los1} for the groups $G(l, 1, n) = \mathfrak{S}_n \rtimes (\mathbb{Z}/ l)^n$.  Here we obtain the result independently for cyclic reflection groups. In addition, we establish the precise parameter values for which $P_{\KZ} \cong \Delta_{c} (S(\mathfrak{h})_W)$, and thus, using the terminology introduced in \cite{Los1}, determine all the totally aspherical parameters for  $H_c(\mathbb{Z}/n, \mathbb{C})$. This work was part of the authors D.Phil thesis, completed in October 2014. Furthermore, we complete the algebraic description of the KZ-functor in this case, by giving an explicit description of the action of the cyclotomic Hecke algebra $\mathcal{H}_{\boldsymbol{c}} (\mathbb{Z}/n, \mathbb{C})$ on $\Delta_{\boldsymbol{c}}(S(\mathfrak{h})_W)$ for parameter values in $\mathcal{F}$.

The structure of this article is as follows:  In Section \ref{chapter 00}, we include some background material on rational Cherednik algebras, category $\mathcal{O}$, the KZ-functor and the coinvariant algebra, focusing in particular on the facts relevant to the rational Cherednik algebras $H_c(\mathbb{Z}/n, \mathbb{C})$, and in Section \ref{chapter 000} we list the main results. We then study the space $\mathcal{S} = \Hom_{H_c (W, \mathfrak{h})} (\Delta_{c} (S(\mathfrak{h})_W), \nabla_{c}(S(\mathfrak{h}^*)_W \otimes_{\mathbb{C}} \text{det}^{-1}_\mathfrak{h} \langle N \rangle) )$ for general complex reflection groups $W$ in Section \ref{chapter 3}, and show that its elements are in one-to-one correspondence with the elements of $\Hom_{\mathbb{C}W} (S(\mathfrak{h}^*)_W, \textnormal{Har}^*)$. From Section \ref{chapter 5} onwards, we specialise to the rational Cherednik algebras $H_c(\mathbb{Z}/n, \mathbb{C})$ associated with cyclic reflection groups. In this setting, we give a criterion for when an element of  $\Hom_{\mathbb{C}W} (S(\mathfrak{h}^*)_W, \textnormal{Har}^*)$ lifts to an isomorphism in $\mathcal{S}$. Using this criterion, we identify a distinguished element of $\Hom_{\mathbb{C}W} (S(\mathfrak{h}^*)_W, \textnormal{Har}^*)$ which lifts to an isomorphism in $\mathcal{S}$ for parameter values in a certain set $\mathcal{F}$. It thus follows that $\Delta_{c} (S(\mathfrak{h})_W)$ is a so-called tilting module (in the sense of \cite{Don}) for these parameter values. As it was proved in \cite{GGOR} that $P_{\KZ}$ is a tilting module, and in  \cite{DR}  that the collection of tilting modules contains a distinguished set of indecomposable objects which form a $\mathbb{Z}$-basis for the Grothendieck group $K_0(\mathcal{O}_c(W, \mathfrak{h}))$, it then follows that $P_{\KZ} \cong \Delta_{c} (S(\mathfrak{h})_W)$ for parameter values in $\mathcal{F}$, as $[P_{\KZ}] = [\Delta_{c} (S(\mathfrak{h})_W)]$ in $K_0(\mathcal{O}_c(W, \mathfrak{h}))$. We also show that all parameter values can be obtained from $\mathcal{F}$ by integral translations. In Section \ref{chapter 6} we then show, by examining the endomorphism rings of $P_{\KZ} $ and $\Delta_{c} (S(\mathfrak{h})_W)$, that $P_{\KZ} $ and $\Delta_{c} (S(\mathfrak{h})_W)$ are isomorphic in category $\mathcal{O}_c(\mathbb{Z}/n, \mathbb{C})$ precisely for the parameter values contained in $\mathcal{F}$. We conclude in Section \ref{chapter 7} by determining the action of the Hecke algebra  $\mathcal{H}_{\boldsymbol{c}} (\mathbb{Z}/n, \mathbb{C})$ on $\Delta_{c} (S(\mathfrak{h})_W)$ for these parameter values, thus giving a complete algebraic description of the KZ-functor in this case.

\section{Background and setup} \label{chapter 00}

In this section, we recall some background material that will be needed for what follows. For more comprehensive overviews, in addition to the articles cited below, the reader is referred to \cite{Ari}, \cite{Bro}, \cite{Gor} and \cite{Rou}. 

\subsection{Complex reflection groups}  \label{subsection00.1}

If $\mathfrak{h}$ is a finite-dimensional complex vector space, a non-identity element $s \in \text{GL}_{\mathbb{C}} (\mathfrak{h})$ is called a pseudo-reflection if it has finite multiplicative order and
fixes a hyperplane pointwise, that is if $\text{Ker} (s - \text{id}_{\mathfrak{h}})$ has codimension $1$ in $\mathfrak{h}$.  A complex reflection group $W$ is a finite subgroup of $\text{GL}_{\mathbb{C}} (\mathfrak{h})$ which is generated by pseudo-reflections. If $W$ is a complex reflection group, we denote by $\textnormal{Refl}_W$ its set of pseudo-reflections, let $N = |\textnormal{Refl}_W|$, and for $s \in \textnormal{Refl}_W$, we let $H_s = \text{Ker} (s - \text{id}_{\mathfrak{h}})$ be the corresponding reflection hyperplane. Furthermore, we let $\mathcal{A}_W = \{H_s : s \in \textnormal{Refl}_W  \}$, and for $H \in \mathcal{A}_W$, we denote by $W_H$ the cyclic subgroup consisting of the elements of $W$ that fix $H$ pointwise, and by $v_H$ an element of $\mathfrak{h}$ (unique up to multiplication by a scalar) such that $\mathfrak{h} = H \oplus \mathbb{C}v_H$. For each $s \in W_H \setminus \{ 0 \}$, we let $v_s = v_H$, and we denote by $\alpha_H$ as well as $\alpha_s$ the unique element of $\mathfrak{h}^*$ such that $\textnormal{Ker} \alpha_H = H$ and $(v_H, \alpha_H) := \alpha_H (v_H) = 1$. We also let $\lambda_s$ be the eigenvalue of $\alpha_s$ with respect to $s$, so that $\lambda_s^{-1}$ is the eigenvalue of $v_s$ with respect to $s$. The complex reflection groups that we shall be mainly occupied with are the cyclic ones, where $\mathfrak{h} = \mathbb{C}$ and $W = \mathbb{Z}/n = \langle s \rangle $ with $s.w = qw = e^{2\pi i /n}w$ for each $w \in \mathbb{C}$. We then have $\textnormal{Refl}_{\mathbb{Z}/n} =(\mathbb{Z}/n) \setminus \{ 0 \} $ and $\mathcal{A}_{\mathbb{Z}/n} = \{ H_s \} = \{ \{0 \} \}$, and letting $\{ \xi \}$ be a basis of $\mathbb{C}$ with $\{ x \}$ the corresponding dual basis of $\mathbb{C}^*$, we let $v_s = \xi$, so that $\alpha_s = x$ and $\lambda_s = q^{-1}$.

\subsection{Rational Cherednik algebras} \label{subsection00.2}
If $W \leq \text{GL}_{\mathbb{C}} (\mathfrak{h})$ is a complex reflection group, it acts on its set of pseudo-reflection $\textnormal{Refl}_W$ by conjugation, and we can therefore speak of $W$-equivariant functions $c: \textnormal{Refl}_W \rightarrow \mathbb{C}$, that is functions that are constant on conjugacy classes of pseudo-reflections. These functions can be identified with $\mathbb{C} [\textnormal{Refl}_W]^W$ in the obvious way. For $c \in \mathbb{C} [\textnormal{Refl}_W]^W$, the rational Cherednik algebra $H_c(W, \mathfrak{h})$ is the quotient of $T(\mathfrak{h} \oplus \mathfrak{h}^*) \rtimes W$ by the relations 
\begin{equation*}   \label{equation00.1}
[\xi_1, \xi_2] = 0, \quad [x_1, x_2] = 0 \quad \textnormal{and} \quad [\xi_1, x_1] =  (\xi_1, x_1) + \sum_{s \in \textnormal{Refl}_W} c(s) (\xi_1, \alpha_s)(v_s, x_1) s
\end{equation*} 
for all $\xi_1, \xi_2 \in \mathfrak{h}$ and $x_1, x_2 \in \mathfrak{h}^*$. It was proved in \cite{EtGi} that $H_c(W, \mathfrak{h})$ has a PBW-type decomposition in the sense that $S(\mathfrak{h}^*)$, $\mathbb{C}W$ and $S(\mathfrak{h})$ embed as subalgebras of $H_c(W, \mathfrak{h})$ and $H_c(W, \mathfrak{h}) = S(\mathfrak{h}^*) \otimes_{\mathbb{C}} \mathbb{C}W \otimes_{\mathbb{C}} S(\mathfrak{h})$ as a vector space, that is,  
 if $\{\xi_1, \ldots, \xi_n \}$ is a basis for $\mathfrak{h}$ and $\{x_1, \ldots, x_n \}$ is a basis for $\mathfrak{h}^*$, then
$\{ x_1^{a_1} \cdots x_n^{a_n} w \xi_1^{b_1} \cdots \xi_n^{b_n} : a_i, b_i \in \mathbb{N} \,\, \textnormal{and} \,\, w \in W   \}$
is a basis for $H_c(W, \mathfrak{h})$. Furthermore, $H_c(W, \mathfrak{h})$ has a $\mathbb{Z}$-grading, with $\mathfrak{h}^*$ sitting in degree $1$, $W$ sitting in degree $0$ and $\mathfrak{h}$ sitting in degree $-1$ (we let $S(\mathfrak{h})^k = S(\mathfrak{h})_{-k}$). As is proved after Lemma 2.5 of \cite{BeEtGi}, this grading is inner in the sense that, with 
\begin{equation*}   \label{equation00.2}
\mathbf{eu}_{\mathbb{C}} = \sum_{i=1}^n x_i \xi_i \quad \textnormal{and} \quad 
\mathbf{h} = \mathbf{eu}_{\mathbb{C}} + \sum_{s \in \textnormal{Refl}_W} \frac{c(s)}{1 - \lambda_s} s, 
\end{equation*} 
we have that $[\mathbf{h}, h] =\textnormal{deg}(h) h$ for every homogeneous $h \in H_c(W, \mathfrak{h})$. In the case when $W = \mathbb{Z}/n$ and $\mathfrak{h} = \mathbb{C}$, the action of $\mathbb{Z}/n$ on $\textnormal{Refl}_{\mathbb{Z}/n}$ is trivial, so that the set of $\mathbb{Z}/n$-equivariant functions $\textnormal{Refl}_{\mathbb{Z}/n} \rightarrow \mathbb{C}$ can be identified with $\mathbb{C}^{n-1}$, and the rational Cherednik algebra $H_c(\mathbb{Z}/n, \mathbb{C})$ is generated by $x$, $s$ and $\xi$ satisfying the relations 
\begin{equation}   \label{equation00.3}
s\xi s^{-1} = q\xi, \quad sxs^{-1} = q^{-1}x \quad \textnormal{and} \quad [\xi, x] =  1 + \sum_{i=1}^{n-1} c(s^i)  s. 
\end{equation} 
\subsection{An alternative parametrisation} \label{subsection00.3}
There is an alternative parametrisation of the rational Cherednik algebras associated with $W$ and $\mathfrak{h}$ that will sometimes be more convenient for our purposes. For $H \in \mathcal{A}_W$, let $e_H = |W_H|$, and for $i \in \mathbb{Z}$, let $\varepsilon_{H,i} = \frac{1}{e_H} \sum_{w \in W_H} \textnormal{det}(w)^i w$. Then $\{\varepsilon_{H, i} : 0 \leq i \leq e_H - 1\}$ is a complete set of primitive orthogonal idempotents for $\mathbb{C} W_H$, and therefore, for each $c \in \mathbb{C} [\textnormal{Refl}_W]^W$, there are unique constants $\{k_{H,i} : H \in \mathcal{A}_W,  0 \leq i \leq e_H  \}$ such that $k_{H, 0} = k_{H, e_H} = 0$ and 
\begin{equation}   \label{equation00.4}
\sum_{s \in W_H \setminus \{ 1 \}} c(s) s = e_H \sum_{i = 0}^{e_H - 1} (k_{H, i+1} - k_{H,i}) \varepsilon_{H,i}.  
\end{equation}
As $w \varepsilon_{H, i} w^{-1} = \varepsilon_{w(H), i} $ and $c$ is $W$-equivariant, it follows that $k_{H, i} = k_{w(H), i}$ for all $w \in W$ and $H \in \mathcal{A}_W$. Denoting by $\mathbb{C}[\mathbf{k}_{H, i}]$ the polynomial ring in indeterminates $\mathbf{k}_{H, i}$ with $H \in \mathcal{A}_W$ and $1 \leq i \leq e_H-1$ such that $\mathbf{k}_{w(H), i}= \mathbf{k}_{H, i}$ for all $w \in W$, and making the obvious identification between elements of $\mathbb{C}[\mathbf{k}_{H, i}]$ and sets $k = \{k_{H,i} :  H \in \mathcal{A}_W, 1 \leq i \leq e_H-1\,\, \textnormal{and} \,\, k_{w(H), i} = k_{H,i} \,\, \textnormal{for all} \,\, w \in W  \}$, equation \eqref{equation00.4} gives a one-to-one correspondence between $ \mathbb{C} [\textnormal{Refl}_W]^W$ and $\mathbb{C}[\mathbf{k}_{H, i}]$. We write $H_k(W, \mathfrak{h})$ with $k \in \mathbb{C}[\mathbf{k}_{H, i}]$ to indicate that we are using the alternative parametrisation. With the help of the constants $k_{H, i}$, we can also define the elements 
\begin{equation*}   \label{equation00.5}
a_H = \sum_{i=1}^{e_H -1} e_H k_{H,i} \varepsilon_{H,i} \quad \textnormal{and} \quad z = \sum_{H \in \mathcal{A}_W} a_H. 
\end{equation*}
Letting $\mathbf{eu} = \mathbf{eu}_{\mathbb{C}} - z$ (cf. Proposition 2.4 of \cite{DuOp}), it can then be seen that $ \mathbf{h} - \mathbf{eu}  = \sum_{H \in \mathcal{A}_W} \Big( \sum_{i=1}^{e_H - 1} k _{H, i} \Big) \in \mathbb{C}$, so that $[\mathbf{eu}, h] =\textnormal{deg}(h) h$ for every homogeneous $h \in H_c(W, \mathfrak{h})$. Furthermore it is proved in Lemma 2.5 of \cite{DuOp} that $z \in Z(\mathbb{C}W)$, as $w a_H w^{-1} = a_{w(H)}$, and therefore $z$ acts by multiplication by a scalar $c_E$ on every $E \in \textnormal{Irr}(W)$. These constants play an important role in the representation theory of $H_k(W, \mathfrak{h})$. Returning again to the case $W = \mathbb{Z}/n$ and $\mathfrak{h} = \mathbb{C}$, we let $E_i = \mathbb{C}.v_i$ be the element of $\textnormal{Irr}(\mathbb{Z}/n)$ such that $s.v_i = q^{-i}v_i$. Letting $k_i = k_{\{0\}, i}$ and $\varepsilon_i := \varepsilon_{\{ 0 \}, i}$, it then follows that $\varepsilon_i$ fixes $v_j$ if $i \equiv j \pmod n$ and kills it otherwise. As $z = a_{\{ 0 \}} = \sum_{i=1}^{n -1} n k_i \varepsilon_i  $, it follows that $c_i := c_{E_i} = nk_i$ for $i \in \mathbb{Z}$. The rational Cherednik algebras associated with $W = \mathbb{Z}/n$ and $\mathfrak{h} = \mathbb{C}$ can therefore also be parametrised in terms of the $\boldsymbol{c} := (c_1, \ldots, c_n) $, where $c_n = 0$ and the other $c_i$ can vary freely in $\mathbb{C}$. We write $H_{\boldsymbol{c}}(\mathbb{Z}/n, \mathbb{C})$ to indicate that we are using this parametrisation. 

\subsection{The standard and costandard modules}  \label{subsection00.4}
Two types of $H_c(W, \mathfrak{h})$-modules are of particular importance. First, for a finite-dimensional, graded $S(\mathfrak{h}) \rtimes W$-module $M$, we let
\begin{equation*}   \label{equation00.6}
\Delta_c (M) = H_c(W, \mathfrak{h}) \otimes_{S(\mathfrak{h}) \rtimes W} M. 
\end{equation*}
 By the PBW Theorem, $\Delta_c(M) \cong S(\mathfrak{h}^*) \otimes_{\mathbb{C}} M$ as a vector space, and $\Delta_c (M)$ is a graded $H_c(W, \mathfrak{h})$-module with $\Delta_c(M)_k = \bigoplus_{i + j = k} (S(\mathfrak{h}^*)_i \otimes_{\mathbb{C}} M_j)$. In particular, if $E$ is a finite-dimensional $W$-module, it can be made into a graded $S(\mathfrak{h}) \rtimes W$-module by declaring that $\mathfrak{h}$ acts by zero and letting $E$ be concentrated in degree zero. This way one can speak of the module $\Delta_c (E)$, introduced in equation (25) of \cite{DuOp} and just after Definition 2.1 of \cite{BeEtGi}. 
Next, for a finite-dimensional, graded $S(\mathfrak{h}^*) \rtimes W$-module $M$, we let 
\begin{equation*}   \label{equation00.7}
\nabla_c (M) = \textnormal{Homgr}^{\bullet}_{S(\mathfrak{h}^*) \rtimes W} (H_c(W, \mathfrak{h}), M) = \bigoplus_{i \in \mathbb{Z}} \textnormal{Homgr}^{i}_{S(\mathfrak{h}^*) \rtimes W} (H_c(W, \mathfrak{h}), M). 
\end{equation*}
The module $\nabla_c(M)$ is a graded  $H_c(W, \mathfrak{h})$-module with its $i$th degree component  $\nabla(M)_i$ equal to $\textnormal{Homgr}^{i}_{S(\mathfrak{h}^*) \rtimes W} (H_c(W, \mathfrak{h}), M)$. Furthermore, it can be seen that $\nabla_c(M)$ consists of the elements of $\Hom_{S(\mathfrak{h}^*) \rtimes W} (H_c(W, \mathfrak{h}), M)$ that are nilpotent with respect to the action of $\mathfrak{h}$ (that is elements that are killed by $\mathfrak{h}^r$ for some $r \geq 0$), which we denote by $\Hom_{S(\mathfrak{h}^*) \rtimes W} (H_c(W, \mathfrak{h}), M)^{ln}$. Similarly to above, if $E$ is a finite-dimensional $W$-module, it can be made into a graded $S(\mathfrak{h}^*) \rtimes W$-module by declaring that $\mathfrak{h}^*$ acts by zero and letting $E$ be concentrated in degree zero. Therefore one can speak of the module $\nabla_c (E)$ (cf. \cite{Gua} and \cite{GGOR}). 
The modules $\Delta_c (M)$ and $\nabla_c (M)$ belong to category $\mathcal{O}_c(W, \mathfrak{h})$ as defined in the next subsection.

\subsection{Category $\mathcal{O}$}  \label{subsection00.5}
A family of $H_c(W, \mathfrak{h})$-modules that is of particular interest is the so-called category $\mathcal{O}_c(W, \mathfrak{h})$, introduced in \cite{DuOp} and \cite{BeEtGi}. It consists of the $H_c(W, \mathfrak{h})$-modules $M$ that are finitely generated as $H_c(W, \mathfrak{h})$-modules and locally nilpotent with respect to the action of $\mathfrak{h}$, that is such that for every $m \in M$, there exists some $r(m) \geq 0$ such that $\mathfrak{h}^{r(m)}.m = 0$. Every module $M$ of category $\mathcal{O}_c(W, \mathfrak{h})$ decomposes as the direct sum of its generalised $\mathbf{eu}$-eigenspaces $\mathcal{W}_{\alpha}(M)$ (cf. \cite{Gua}), and in particular $\Delta_c(E)_i = \mathcal{W}_{i - c_E}(\Delta_c(E))$ and $\nabla_c(E)_i = \mathcal{W}_{i - c_E}(\nabla_c(E))$ for all $E \in \textnormal{Irr}(W)$. 
We let $\overline{\textnormal{Irr}}(W)$ denote a complete set of pairwise non-isomorphic $W$-modules, and endow it with a partial ordering $<_c$ by declaring that for $E,F \in \overline{\textnormal{Irr}}(W)$, we have $E <_c F$ if and only if $c_F - c_E \in \mathbb{Z}_{>0}$. It is then proved in \cite{Gua} that category $\mathcal{O}_c(W, \mathfrak{h})$ is a finite-length highest weight category in the sense of \cite{CPS}, with indexing set $\overline{\textnormal{Irr}}(W)$. The standard modules are the $\Delta_c(E)$, and the costandard modules are the $\nabla_c(E)$ with $E \in \overline{\textnormal{Irr}}(W)$, whereas the simple modules are the $L_c(E)$ that are the heads of the $\Delta_c(E)$, or equivalently the socles of the $\nabla_c(E)$. Furthermore, we denote the projective and injective indecomposable modules by $P_c(E)$ and $I_c(E)$ respectively. We let $\mathcal{F} (\Delta)$ and $\mathcal{F}(\nabla)$ denote the sets of modules that have a finite filtration with quotients isomorphic to modules $\Delta_c(E)$ and $\nabla_c(E)$ respectively, and recall that $\mathcal{F} (\Delta) \cap \mathcal{F}(\nabla)$ is the collection of tilting modules (in the sense of \cite{Don}) in category $\mathcal{O}_c(W, \mathfrak{h})$. The following important consequence of the fact that category $\mathcal{O}_c(W, \mathfrak{h})$ is a finite-length highest weight category with a finite indexing set follows from Proposition 3.1 of \cite{DR}:

\begin{theorem00.1}[{\cite[Proposition 3.1]{DR}}]  \label{theorem00.1}
There is a set $\{T(E): E \in  \overline{\textnormal{Irr}} (W)   \}$ of indecomposable tilting modules in category $\mathcal{O}_c(W, \mathfrak{h}) $  such that every module $M \in \mathcal{F}(\Delta) \cap \mathcal{F}(\nabla)$ is a direct sum of finitely many of the $T(E)$. Furthermore, $\{[T(E)]: E \in  \overline{\textnormal{Irr}} (W)   \}$ is a $\mathbb{Z}$-basis for the Grothendieck group $K_0(\mathcal{O}_c(W, \mathfrak{h}))$. 
\end{theorem00.1}

In particular, any two tilting modules that represent the same element in $K_0(\mathcal{O}_c(W, \mathfrak{h}))$ are necessarily isomorphic. We will be making use of this fact below.

In the case of cyclic reflection groups, we let $\overline{\textnormal{Irr}}(W) = \{E_i : 1 \leq i \leq n \}$ with $E_i$ defined as in Subsection \ref{subsection00.3}. We note the following result, where for two complex numbers $z_1$ and $z_2$ we write $z_1 \equiv z_2 \pmod n$ to mean that $z_1 - z_2 \in n\mathbb{Z}$:

\begin{theorem00.2}  \label{theorem00.2}
Category $\mathcal{O}_{\boldsymbol{c}}(\mathbb{Z}/n, \mathbb{C}) $ is semisimple if and only if $c_i - c_j \not\equiv j - i \pmod n$ for all $1 \leq i, j \leq n$. 
\end{theorem00.2} 

\begin{proof}
It is a consequence of Propositions 5.15 and 5.16 of \cite{GGOR}, that category $\mathcal{O}_{\boldsymbol{c}}(\mathbb{Z}/n, \mathbb{C})$ is semisimple if and only if the cyclotomic Hecke algebra $\mathcal{H}_{\boldsymbol{c}}(\mathbb{Z}/n, \mathbb{C}) $ (cf. the next subsection) is semisimple, and the claim then follows from \eqref{equation00.9} together with the fact that $c_i = nk_i$. 
\end{proof}

\subsection{Cyclotomic Hecke algebras}  \label{subsection00.6}
In this subsection we outline the definition of the cyclotomic Hecke algebra $\mathcal{H}_k(W, \mathfrak{h})$ associated with the rational Cherednik algebra $H_k(W, \mathfrak{h})$. For details, the reader is referred to \cite{BMR} Section 4.C. 

For a complex reflection group $W$ with reflection representation $\mathfrak{h}$, we let $\mathfrak{h}_{\textnormal{reg}} = \mathfrak{h} \setminus \bigcup_{H \in \mathcal{A}_W} H$. Then $W$ acts freely on $\mathfrak{h}_{\textnormal{reg}}$ and, fixing $x_0 \in \mathfrak{h}_{\textnormal{reg}}$, the braid group associated with $W$ and $\mathfrak{h}$ is the fundamental group  $\pi_1 (\mathfrak{h}_{\textnormal{reg}}/W, x_0)$. 
We call  a pseudo-reflection $s \in \textnormal{Refl}_W$ distinguished if it satisfies $\textnormal{det}_{\mathfrak{h}}(s) = \textnormal{exp} (\frac{2\pi i}{e_{H_s}})$.
As the subgroups $W_H$ are cyclic, there is precisely one distinguished pseudo-reflection associated to each hyperplane $H \in \mathcal{A}_W$, and we denote it by $s_H$. It is then proved in \cite{BMR} Theorem 2.17 that 
the braid group $\pi_1 (\mathfrak{h}_{\textnormal{reg}}/W, x_0)$ is generated by elements $\{ T_{s_H} : H \in \mathcal{A}_W \}$, where $T_{s_H}$ is a so-called $s_H$-generator of the monodromy around the image of $H$ in $\mathfrak{h}_{\textnormal{reg}}/W$.  
Furthermore, it is shown in Proposition 4.22 of \cite{BMR} that (with some mild restrictions on the reflection group $W$) the  $T_{s_H}$ satisfy certain braid relations, for the details of which we again refer the reader to \cite{BMR} Section 4.C. The cyclotomic Hecke algebra $\mathcal{H}_k(W, \mathfrak{h})$ is then defined as the quotient of the group algebra $\mathbb{C} \pi_1 (\mathfrak{h}_{\textnormal{reg}}/W, x_0)$ by the two-sided ideal generated by the set $\{ \prod_{j = 1}^{e_H } (T_{s_H} - \textnormal{det}_{\mathfrak{h}}(s_H)^{-j} \textnormal{exp} (- 2 \pi i k_{H, j})) :  H \in \mathcal{A}_W   \}$ 
(in e.g. \cite{GGOR}, the set of generators $\{ T_{s_H} : H \in \mathcal{A}_W \}$ for the cyclotomic Hecke algebra $\mathcal{H}_k(W, \mathfrak{h})$ used here are replaced by their inverses, so that $\mathcal{H}_k(W, \mathfrak{h})$ is defined in terms of the  set $\{ \prod_{j = 1}^{e_H } (T_{s_H} - \textnormal{det}_{\mathfrak{h}}(s_H)^{j} \textnormal{exp} ( 2 \pi i k_{H, j})) :  H \in \mathcal{A}_W   \}$ instead). 
In the case where $W = \mathbb{Z}/n = \langle s \rangle$ and $\mathfrak{h} = \mathbb{C}$, the generator $s$ is the unique distinguished pseudo-reflection of $ \mathbb{Z}/n$, and therefore $ \pi_1 (\mathfrak{h}_{\textnormal{reg}}/W, x_0)$ is generated by $T := T_s$. As $\mathfrak{h}_{\textnormal{reg}} / W = \mathbb{C}^{\times} / (\mathbb{Z}/n)$ is a cylinder, $ \pi_1 (\mathfrak{h}_{\textnormal{reg}}/W, x_0) \cong \mathbb{Z}$, and therefore the corresponding cyclotomic Hecke algebra satisfies 
\begin{equation}  \label{equation00.9}
\mathcal{H}_k(\mathbb{Z}/n, \mathbb{C}) \cong  \mathbb{C}[T] / \langle \prod_{j = 1}^{n } (T -q^{-j} \textnormal{exp} (- 2 \pi i k_j )) \rangle. 
\end{equation}
\subsection{The KZ-functor and $P_{\KZ}$}  \label{subsection00.7}
A critical tool for studying category $\mathcal{O}$ is the KZ-functor introduced in Section 5.3 of \cite{GGOR}. Through localising, making use of the so-called Dunkl embedding (cf. Proposition 4.5 of \cite{EtGi}) and taking monodromy, it is proved in \cite{GGOR} that there exists an exact functor $\KZ: \mathcal{O}_c(W, \mathfrak{h}) \longrightarrow  \mathcal{H}_c(W, \mathfrak{h})\textnormal{-mod}$, where $\mathcal{H}_c(W, \mathfrak{h})\textnormal{-mod}$ denotes the finitely-generated modules over the cyclotomic Hecke algebra $\mathcal{H}_c(W, \mathfrak{h})$. The definition of KZ is geometric in nature, but as the functor is exact, it is represented by a projective module $P_{\KZ}$ in $\mathcal{O}_c(W, \mathfrak{h})$. In other words,  there is an algebra homomorphism $\phi: \mathcal{H}_k(W, \mathfrak{h}) \longrightarrow \End_{H_k(W, \mathfrak{h})} (P_{\KZ})^{\textnormal{opp}}$, 
so that $P_{\KZ}$ has the structure of an $(H_k(W, \mathfrak{h}), \mathcal{H}_k(W, \mathfrak{h}))$-bimodule, giving $\Hom_{H_k(W, \mathfrak{h})} (P_{\KZ}, M)$ the structure of an $\mathcal{H}_k(W, \mathfrak{h})$-module for all modules $M$ in category $\mathcal{O}_k(W, \mathfrak{h})$, such that the functors KZ and $\Hom_{H_k(W, \mathfrak{h})} (P_{\KZ}, -)$ are naturally isomorphic. 
The aim of this article is to give an explicit algebraic description of $P_{\KZ}$ as an $(H_k(W, \mathfrak{h}), \mathcal{H}_k(W, \mathfrak{h}))$-bimodule in the special case of rational Cherednik algebras corresponding to cyclic complex reflection groups for certain parameter values. To this end, we recall a few properties of $P_{\KZ}$ that will be important below. First, as $P_{\KZ}$ is projective, it is contained in $\mathcal{F}(\Delta)$, and  it follows from the definition of $P_{\KZ}$ that it has a decomposition $P_{\KZ} = \bigoplus_{E \in \textnormal{Irr} (W)} \textnormal{dim} (\KZ(L(E))) P(E)$. Using this decomposition together with the fact that $\textnormal{dim}(\Hom_{H_k(W, \mathfrak{h})} (P(E), M)) = [M: L(E)]$ for every module $M$ in category $\mathcal{O}_k(W, \mathfrak{h})$, a calculation similar to that at the end of the proof of Proposition 5.15 of \cite{GGOR} gives the following result (where for modules $M \in \mathcal{F}(\Delta)$ and $E \in \textnormal{Irr} (W)$, we denote by $[M: \Delta(E)]$ the number of times $\Delta(E)$ appears in a $\Delta$-filtration of $M$):  

\begin{theorem00.3} \label{theorem00.3}
The $\Delta$-filtration multiplicities of $P_{\KZ}$ satisfy $[P_{\KZ} : \Delta(E) ] = \textnormal{dim}(E)$
for every $E \in \textnormal{Irr} (W)$. 
\end{theorem00.3}

Furthermore, it is proved in Proposition 5.21 of \cite{GGOR} that $P_{\KZ}$ is also injective, so that it belongs to $\mathcal{F}(\nabla)$. We therefore have the following result, we which will be crucial in what follows. 

\begin{theorem00.4}[{\cite[Proposition 5.21]{GGOR}}] \label{theorem00.4}
The projective object $P_{\KZ}$ is a tilting module. 
\end{theorem00.4}

Since $\mathcal{O}_c(W, \mathfrak{h})$ is a highest weight category,  $\{[\Delta(E)]: E \in \overline{\textnormal{Irr}} (W)  \}$ is a $\mathbb{Z}$-basis for $K_0(\mathcal{O}_c(W, \mathfrak{h}))$, and we obtain the following corollary by combining Proposition \ref{theorem00.1}, Lemma \ref{theorem00.3} and Proposition \ref{theorem00.4}: 

\begin{theorem00.5} \label{theorem00.5}
If $M$ is a tilting module in category $\mathcal{O}_c(W, \mathfrak{h})$ and $[M : \Delta(E) ] = \textnormal{dim}E$
for every $E \in \textnormal{Irr} (W)$, then $M \cong P_{\KZ}$ as $H_c(W, \mathfrak{h})$-modules. 
\end{theorem00.5}

In Sections \ref{chapter 3} and \ref{chapter 5},  we construct such a tilting module $M$ in category $\mathcal{O}_k(\mathbb{Z}/n, \mathbb{C})$ for a collection $\mathcal{F}$ of parameter values, and  Corollary \ref{theorem00.5} then implies that $M$ is isomorphic to $P_{\KZ}$ when $k \in \mathcal{F}$.   The interest in $\mathcal{F}$ comes from the fact that there are ways to relate categories $\mathcal{O}_k(W, \mathfrak{h})$ and $\mathcal{O}_{k'} (W, \mathfrak{h})$ corresponding to different parameter values $k$ and $k'$. 
 Observe that  the definition of the cyclotomic Hecke algebra $\mathcal{H}_k(W, \mathfrak{h})$ implies that it does not change if integer values are added to the $k_{H,i}$. If this is done in such a way that we still have $k'_{H, i} = k'_{w(H), i}$ for all $w \in W$, $H \in \mathcal{A}_W $ and $1 \leq i \leq e_H - 1$, then this corresponds to the cyclotomic Hecke algebra $\mathcal{H}_{k'} (W, \mathfrak{h})$ associated with the rational Cherednik algebra $H_{k'} (W, \mathfrak{h})$ for some new parameter value $k'$. One might then hope that this is an indication that the categories $\mathcal{O}_{k} (W, \mathfrak{h})$ and $\mathcal{O}_{k'} (W, \mathfrak{h})$ are equivalent. This is unfortunately not true in general, but Losev proves in Theorem 1.1 of \cite{Los} that the categories are derived equivalent: 
 
\begin{theorem00.6}[{\cite[Theorem 1.1]{Los}}] \label{theorem00.6}
If $k \in \mathbb{C}[\mathbf{k}_{H,i}]$ and $l \in \mathbb{Z}[\mathbf{k}_{H,i}]$, then there is a derived equivalence $D^b (\mathcal{O}_{k} (W, \mathfrak{h})) \stackrel{\sim}{\longrightarrow} D^b (\mathcal{O}_{k + l} (W, \mathfrak{h}))$. 
Furthermore, this equivalence respects $P_{\KZ}$. 
\end{theorem00.6}

Rouquier also proves in Theorem 5.5 of \cite{Rou} that under certain conditions on the parameter $k$, the categories $\mathcal{O}_{k} (W, \mathfrak{h})$ and $\mathcal{O}_{k +l} (W, \mathfrak{h})$ are in fact equivalent. Taken together, these results motivate the following definition. 

\begin{definition00.1}  \label{definition00.1}
\begin{itemize}
\item[(i)] For two parameter values $k, k' \in \mathbb{C}[\mathbf{k}_{H,i}]$, we write $k \sim k'$ if $k - k' \in \mathbb{Z}[\mathbf{k}_{H,i}]$ and observe that this is an equivalence relation on $\mathbb{C}[\mathbf{k}_{H,i}]$. 
\item[(ii)] We say that a set of parameter values $\mathcal{K} \subseteq \mathbb{C}[\mathbf{k}_{H,i}]$ is good if, for every parameter value $k \in \mathbb{C}[\mathbf{k}_{H,i}]$, there exists $k' \in \mathcal{K} $ such that $k \sim k'$. 
\end{itemize}
\end{definition00.1}

\begin{remark00.1}  \label{remark00.1}
If we think of $\mathcal{O}_{c} (W, \mathfrak{h})$ as parametrised in terms of $\mathbb{C}[\textnormal{Refl}_W]^W$ instead, we call a set of parameter values $\mathcal{C} \subseteq \mathbb{C}[\textnormal{Refl}_W]^W$ good if the corresponding set $\mathcal{K} \subseteq \mathbb{C}[\mathbf{k}_{H,i}]$ is good. 
\end{remark00.1}

In Section \ref{chapter 5} we give an explicit algebraic description of $P_{\KZ}$ in category $\mathcal{O}_{k} (\mathbb{Z}/n, \mathbb{C})$ for a good set of parameter values $\mathcal{F}$.

\subsection{The coinvariant algebra}  \label{subsection00.8}
In order to give an algebraic description of $P_{\KZ}$ in category $\mathcal{O}_{k} (\mathbb{Z}/n, \mathbb{C})$ for a good set of parameter values, we will make use of the coinvariant algebra. If $W$ is a complex reflection group with corresponding reflection representation $\mathfrak{h}$, the symmetric algebra $S(\mathfrak{h})$ has a natural structure of a graded $S(\mathfrak{h}) \rtimes W$-module given by $(g \otimes w). g' = g(w.g')$ for $w \in W$ and $g, g' \in S(\mathfrak{h})$, with $\mathfrak{h}$ sitting in degree $-1$.  
We let $S(\mathfrak{h})^W = \{ g \in S(\mathfrak{h}) : w.g = g \,\, \textnormal{for all} \,\, w \in W  \} $
be the subalgebra of $S(\mathfrak{h})$ consisting of the $W$-invariant elements of $S(\mathfrak{h})$, and we define the coinvariant algebra $S(\mathfrak{h})_W$ to be the quotient $S(\mathfrak{h})_W = S(\mathfrak{h}) / (S(\mathfrak{h}). S(\mathfrak{h})^W_{<0})$. This is a graded $S(\mathfrak{h}) \rtimes W$-module with the image of $\mathfrak{h}$ sitting in degree $-1$. The importance of it in our setting comes from the following result, proved in \cite{Ch}: 

\begin{theorem00.7}[{\cite[Theorem B]{Ch}}] \label{theorem00.7}
\begin{itemize} \item[(i)] There is an isomorphism $S(\mathfrak{h})_W \cong \mathbb{C} W$
as ungraded $W$-modules. 

\item[(ii)] The coinvariant algebra $S(\mathfrak{h})_W$ is concentrated in degrees between $- N$ and $0$, where we recall that $N = |\textnormal{Refl}_W|$. 
\end{itemize}

\end{theorem00.7} 

As a consequence of Proposition \ref{theorem00.7} and the fact that $H_{c} (W, \mathfrak{h})$ is free as a right $S(\mathfrak{h}) \rtimes W$-module by the PBW Theorem, we have the following: 

\begin{theorem00.8} \label{theorem00.8}
The module $\Delta_c(S(\mathfrak{h})_W)$ belongs to category $\mathcal{O}_c(W, \mathfrak{h})$ and satisfies that $[\Delta_c(S(\mathfrak{h})_W): \Delta_c(E)] = \textnormal{dim} (E)$
for every $E \in \textnormal{Irr} (W)$. 
\end{theorem00.8}

The following has been conjectured by Rouquier (unpublished): 

\begin{theorem00.9*}  \label{theorem00.9}
$\mathbf{(Rouquier)}$ For a general rational Cherednik algebra $H_c(W, \mathfrak{h})$, there exists a good set of parameter values $\mathcal{C} \subseteq \mathbb{C}[\textnormal{Refl}_W]^W$, such that $P_{\KZ} \cong \Delta_c(S(\mathfrak{h})_W)$
when $c \in \mathcal{C}$. 
\end{theorem00.9*}

We observe that it follows from Corollaries \ref{theorem00.5} and \ref{theorem00.8} that if $\Delta_c(S(\mathfrak{h})_W)$ has a $\nabla$-filtration, so that it is a tilting module, then $\Delta_c(S(\mathfrak{h})_W) \cong P_{\KZ}$. Note that this is automatically true for parameter values such that category $\mathcal{O}_c(W, \mathfrak{h}) $ is semisimple. We will prove that it holds in category $\mathcal{O}_{c} (\mathbb{Z}/n, \mathbb{C})$ for a good set of parameter values, thus confirming Rouquier's conjecture in this case. We do this by showing that $\Delta_c(S(\mathfrak{h})_W)$ is isomorphic to a module in $\mathcal{F} (\nabla)$ that we now introduce. For the reflection representation $\mathfrak{h}^*$ of $W$, we can analogously to above consider the coinvariant algebra corresponding to $\mathfrak{h}^*$, defined by $S(\mathfrak{h}^*)_W = S(\mathfrak{h}^*) / (S(\mathfrak{h}^*). S(\mathfrak{h}^*)^W_{>0})$ (observe that we take the image of $\mathfrak{h}^*$ to sit in degree $1$). Letting $\textnormal{det}_{\mathfrak{h}}$ denote the $1$-dimensional $W$-representation obtained by taking the determinant, we then have that $S(\mathfrak{h}^*)_W^{\otimes} := S(\mathfrak{h}^*)_W \otimes \textnormal{det}^{-1}_{\mathfrak{h}} \langle N \rangle$ 
is a graded $S(\mathfrak{h}^*) \rtimes W$-module concentrated in degrees between $-N$ and $0$, where $S(\mathfrak{h}^*) $ acts naturally on the first factor, leaving the second factor fixed, and $W$ acts diagonally. Similarly to Proposition \ref{theorem00.8} we then have: 

\begin{theorem00.10} \label{theorem00.10}
The module $\nabla_c(S(\mathfrak{h}^*)_W^{\otimes})$ belongs to category $\mathcal{O}_c(W, \mathfrak{h})$ and satisfies $[\nabla_c(S(\mathfrak{h}^*)_W^{\otimes}): \nabla_c(E)] = \textnormal{dim} (E)$
for every $E \in \textnormal{Irr} (W)$. 
\end{theorem00.10}

As a consequence we have: 

\begin{theorem00.11} \label{theorem00.11}
The modules $\Delta_c(S(\mathfrak{h})_W)$ and $\nabla_c(S(\mathfrak{h}^*)_W^{\otimes})$ represent the same element in the Grothendieck group $K_0(\mathcal{O}_c(W, \mathfrak{h}))$.
\end{theorem00.11}
\begin{proof}
This follows from Propositions \ref{theorem00.8} and \ref{theorem00.10}, as $[\Delta(E)] = [\nabla(E)]$ in $K_0(\mathcal{O}_c(W, \mathfrak{h}))$ for every $E \in \textnormal{Irr} (W)$ by Proposition 3.3 of \cite{GGOR}. 
\end{proof}
When $W = \mathbb{Z}/n$ and $\mathfrak{h} = \mathbb{C}$, we have that $S(\mathfrak{h})^W = \mathbb{C}[\xi]^{\mathbb{Z}/n} =  \mathbb{C}[\xi^n]$, so that $S(\mathfrak{h})_W = \mathbb{C}[\xi]/ \langle \xi^n \rangle$. Similarly, $S(\mathfrak{h}^*)_W = \mathbb{C}[x]/ \langle x^n \rangle$. 

\subsection{Harmonic polynomials}  \label{subsection00.9}

In order to find parameter values $c$ for which  $\Delta_c(S(\mathfrak{h})_W)$ is a tilting module, we will make use of the so-called harmonic polynomials (cf. Sections 3.3.3 and 4.3 of \cite{Bro}). Let $W$ be a complex reflection group with reflection representation $\mathfrak{h}$. We recall from Proposition 3.20 of \cite{Bro} that there exists a unique $\mathbb{C}$-algebra homomorphism $D:  S(\mathfrak{h}^*) \rightarrow \End_{\mathbb{C}} (S(\mathfrak{h}))$ such that $D(x)(\xi) = (\xi, x)$ for all $x \in \mathfrak{h}^*$ and $\xi \in \mathfrak{h}$, and such that $D(x)(gg') = D(x)(g)g' + gD(x)(g')$
for all $x \in \mathfrak{h}^*$ and $g, g' \in S(\mathfrak{h})$. Thus, $S(\mathfrak{h})$ has the structure of an $S(\mathfrak{h}^*)$-module, and furthermore $w.D(f)(g) = D(w.f)(w.g)$
for all  $f \in S(\mathfrak{h}^*) $, $g \in S(\mathfrak{h}) $ and $w \in W$. This gives rise to a duality $\langle -, - \rangle : S(\mathfrak{h}^*) \times S(\mathfrak{h})  \rightarrow \mathbb{C}$ given by $\langle f , g \rangle = D(f)(g)(0)$, where $D(f)(g)(0)$ denotes the coefficient of the degree zero component of $D(f)(g)$. It is then proved in Lemma 3.21 of \cite{Bro} that if $\{\xi_1, \ldots, \xi_k \}$ is a basis for $\mathfrak{h}$ and $\{x_1, \ldots, x_k  \}$ is the corresponding dual basis for $\mathfrak{h}^*$, and we let $\xi^{\boldsymbol{m}} = \xi_1^{m_1} \cdots \xi_k^{m_k}$ and $\boldsymbol{m}! = m_1 ! \cdots m_k !$ for every $\boldsymbol{m} = (m_1, \ldots, m_k) \in \mathbb{N}^k$, then $\langle x^{\boldsymbol{n}}, \xi^{\boldsymbol{m}}\rangle$ equals $\boldsymbol{m}!$ if $ \boldsymbol{n} = \boldsymbol{m}$ and $0$ otherwise. One then defines the so-called harmonic polynomials by $\textnormal{Har} := \{ g \in S(\mathfrak{h}) : \langle f , g \rangle = 0 \,\, \textnormal{for all} \,\, f \in S(\mathfrak{h}^*).S(\mathfrak{h}^*)^W_{>0}   \}$. This is a homogeneous $W$-submodule of $S(\mathfrak{h})$, and Proposition 3.24 of \cite{Bro} says that there is a decomposition $S(\mathfrak{h}) =   (S(\mathfrak{h}).S(\mathfrak{h})^W_{<0})  \oplus   \textnormal{Har}$ as graded $W$-modules. We also recall the the Jacobian associated to $W$ is the element $J := \prod_{H \in \mathcal{A}_W} v_H^{e_H - 1} \in S(\mathfrak{h})$. When $W = \mathbb{Z}/n$ and $\mathfrak{h} = \mathbb{C}$, we have that $S(\mathfrak{h}^*).S(\mathfrak{h}^*)^W_{>0} = \langle x^n \rangle$, so that $\textnormal{Har} = \textnormal{Span}_{\mathbb{C}} \{\xi^i : 0 \leq i \leq n-1  \}$, and $J = \xi^{n-1}$. 

\subsection{More general rings}  \label{subsection00.10}

This subsection will be needed for Section \ref{chapter 7} below. 
Let $k$ be a local, noetherian, commutative $\mathbb{C}$-algebra.  Rational Cherednik algebras associated with a complex reflection group $W \leq \text{GL}_{\mathbb{C}} (\mathfrak{h})$ can be considered over $k$ rather than $\mathbb{C}$. If $c: \textnormal{Refl}_W \rightarrow k$ is a $W$-equivariant function, the rational Cherednik algebra $H^k_c(W, \mathfrak{h})$ is defined, analogously to above, as the quotient of $k \otimes_{\mathbb{C}} (T(\mathfrak{h} \oplus \mathfrak{h}^*) \rtimes W)$ by the relations 
\begin{equation*}   \label{equation00.20}
[\xi_1, \xi_2] = 0, \quad [x_1, x_2] = 0 \quad \textnormal{and} \quad [\xi_1, x_1] =  (\xi_1, x_1) + \sum_{s \in \textnormal{Refl}_W} c(s) (\xi_1, \alpha_s)(v_s, x_1) s
\end{equation*} 
for all $\xi_1, \xi_2 \in \mathfrak{h}$ and $x_1, x_2 \in \mathfrak{h}^*$. Much of the above theory carries over to this setting. In particular, the results of \cite{EtGi} again apply to give a PBW-type decomposition $H^k_c(W, \mathfrak{h}) = (k \otimes_{\mathbb{C}} S(\mathfrak{h}^*)) \otimes_{k} kW \otimes_{k} (k \otimes_{\mathbb{C}} S(\mathfrak{h}))$ as $k$-modules. Furthermore, $H^k_c(W, \mathfrak{h})$ again has a $\mathbb{Z}$-grading which is inner with respect to the element $\mathbf{eu} = \mathbf{eu}_{\mathbb{C}} - z$ with $\mathbf{eu}_{\mathbb{C}}$ and $z$ defined entirely analogously to in Subsections \ref{subsection00.2} and \ref{subsection00.3} above, and there is again an alternative parametrisation as described in Subsection \ref{subsection00.3} (we sometimes write $k^k_{H, i}$ and $c^k_{E}$ to emphasise that we are working over $k$). The standard and costandard modules are defined as in Subsection \ref{subsection00.4} (observe that $H^k_c(W, \mathfrak{h})$ has subalgebras isomorphic to $k \otimes_{\mathbb{C}}(S(\mathfrak{h}) \rtimes W)$ and $k \otimes_{\mathbb{C}}(S(\mathfrak{h}^*) \rtimes W)$) and category $\mathcal{O}^k_c (W, \mathfrak{h})$ again consists of the finitely generated $H^k_c(W, \mathfrak{h})$-modules that are  locally nilpotent with respect to the action of $\mathfrak{h}$. Like category $\mathcal{O}_c (W, \mathfrak{h})$, the category $\mathcal{O}^k_c (W, \mathfrak{h})$ is a highest weight category, though now in the sense of Definition 4.11 of \cite{Rou2} (cf. Theorem 5.2 of \cite{Rou2}). The cyclotomic Hecke algebra described in Subsection \ref{subsection00.6} can be defined over $k$ as well, as long as exponentials are well-defined in $k$, in which case $\mathcal{H}^k_c(W, \mathfrak{h})$ is defined as the quotient of the group algebra $k \pi_1 (\mathfrak{h}_{\textnormal{reg}}/W, x_0)$ by the two-sided ideal generated by the set $\{ \prod_{j = 1}^{e_H } (T_{s_H} - \textnormal{det}_{\mathfrak{h}}(s_H)^{-j} \textnormal{exp} (- 2 \pi i k_{H, j})) :  H \in \mathcal{A}_W   \}$. 
We will make use of these concepts in Section \ref{chapter 7}.

\section{Results} \label{chapter 000}

In Section \ref{chapter 3}, we investigate maps between the modules $\Delta_c(S(\mathfrak{h})_W)$ and $\nabla_c(S(\mathfrak{h}^*)_W^{\otimes})$ and prove the following result: 

\begin{theorem000.1} \label{theorem000.1}
There is a one-to-one correspondence between the elements of the spaces $\textnormal{Hom}_{H_c (W, \mathfrak{h})} (\Delta_{c} (S(\mathfrak{h})_W), \nabla_{c}(S(\mathfrak{h}^*)_W^{\otimes}) )$ and $\textnormal{Hom}_{\mathbb{C}W} (S(\mathfrak{h}^*)_W, \textnormal{Har}^*)$. 
\end{theorem000.1}

Specialising to the setting of the rational Cherednik algebras $H_{c}(\mathbb{Z}/n,\mathbb{C})$ in Section \ref{chapter 5}, we use Theorem A  to obtain the following: 

\begin{theorem000.2} \label{theorem000.2}
There is a good set of parameter values $\mathcal{F}$ such that $P_{\textnormal{KZ}} \cong \Delta_{\boldsymbol{c}}(S(\mathfrak{h})_W)$ 
in category $\mathcal{O}_{\boldsymbol{c}} (\mathbb{Z}/n, \mathbb{C})$ when $\boldsymbol{c} \in \mathcal{F}$, where $\mathcal{F}$ consists of those parameters $\boldsymbol{c}$ which satisfy that $c_i - c_j = j - i$ whenever $c_i - c_j \equiv j - i \pmod n$ for $1 \leq i,j \leq n$. 
\end{theorem000.2}

Examining the endomorphism rings of $P_{\textnormal{KZ}}$ and $\Delta_{\boldsymbol{c}}(S(\mathfrak{h})_W)$ in Section \ref{chapter 6}, we obtain a precise condition on the parameter $\boldsymbol{c}$ for $P_{\textnormal{KZ}}$ and $\Delta_{\boldsymbol{c}}(S(\mathfrak{h})_W)$ to be isomorphic in category $\mathcal{O}_{\boldsymbol{c}} (\mathbb{Z}/n, \mathbb{C})$:

\begin{theorem000.3} \label{theorem000.3}
In category $\mathcal{O}_{\boldsymbol{c}} (\mathbb{Z}/n, \mathbb{C})$, the modules $P_{\textnormal{KZ}}$ and $\Delta_{\boldsymbol{c}}(S(\mathfrak{h})_W)$ are isomorphic if and only if  $\boldsymbol{c} \in \mathcal{F}$. 
\end{theorem000.3}

From Theorem B and the fact that $\mathcal{H}_{\boldsymbol{c}} (\mathbb{Z}/n, \mathbb{C})$ is commutative, it follows that for parameter values in $\mathcal{F}$, the KZ-functor induces an algebra homomorphism $\phi$ from $\mathcal{H}_{\boldsymbol{c}} (\mathbb{Z}/n, \mathbb{C})$ to $\textnormal{End}_{H_{\boldsymbol{c}} (\mathbb{Z}/n, \mathbb{C})} (\Delta_{\boldsymbol{c}}(S(\mathfrak{h})_W))$, which is in fact an isomorphism by Theorem 5.15 of \cite{GGOR}. We give an explicit description of this isomorphism in Section \ref{chapter 7}, thus specifying the $(H_k(W, \mathfrak{h}), \mathcal{H}_k(W, \mathfrak{h}))$-bimodule structure on $\Delta_{\boldsymbol{c}}(S(\mathfrak{h})_W)$ and thereby completing the algebraic description of the KZ-functor in this case. 

\begin{theorem000.4}  \label{theorem000.4}
For  $\boldsymbol{c} \in \mathcal{F}$, the isomorphism $\phi: \mathcal{H}_{\boldsymbol{c}} (\mathbb{Z}/n, \mathbb{C})  \stackrel{\sim}{\longrightarrow} \textnormal{End}_{H_{\boldsymbol{c}} (\mathbb{Z}/n, \mathbb{C})} (\Delta_{\boldsymbol{c}}(S(\mathfrak{h})_W))$ induced by the $\KZ$-functor is given by 
$\phi(\overline{T}) = \eta = s.\textnormal{exp}(\frac{2\pi i}{n} \mathbf{eu})$, where $\overline{T}$ is the natural generator of  $\mathcal{H}_{\boldsymbol{c}} (\mathbb{Z}/n, \mathbb{C}) = \mathbb{C}[T]/ \langle \prod^n_{j = 1} (T - q^{-j}.q_j^{-1}) \rangle$.  
\end{theorem000.4}

\section{The space $\text{Hom}_{H_c(W, \mathfrak{h})} (\Delta_c (S(\mathfrak{h})_W), \nabla_c (S(\mathfrak{h}^*)_W^{\otimes}) ) $} \label{chapter 3}

In order to find parameter values $c$ for which $\Delta_{c} (S(\mathfrak{h})_W)$ is a tilting module, so that $P_{\textnormal{KZ}} \cong \Delta_{c} (S(\mathfrak{h})_W)$ in category $\mathcal{O}_c (W, \mathfrak{h})$, we investigate the space 
\begin{equation*}  \label{equation3.1}
	\mathcal{S} := \text{Hom}_{H_c (W, \mathfrak{h})} (\Delta_{c} (S(\mathfrak{h})_W), \nabla_{c}(S(\mathfrak{h}^*)_W^{\otimes}) )
\end{equation*}
in this section, where we recall that $S(\mathfrak{h}^*)_W^{\otimes} = S(\mathfrak{h}^*)_W \otimes_{\mathbb{C}} \text{det}^{-1}_\mathfrak{h} \langle N \rangle$, and show that there is a one-to-one correspondence between the elements of $\mathcal{S}$ and the elements of
$\textnormal{Hom}_{\mathbb{C}W} (S(\mathfrak{h}^*)_W, \textnormal{Har}^*)$. In Section \ref{chapter 5}  we identify the elements of this space in the special setting of 
category $\mathcal{O}_c (\mathbb{Z}/n, \mathbb{C})$, and calculate their liftings to $\mathcal{S}$. It turns out that one of these elements lifts to an isomorphism between $\Delta_{c} (S(\mathfrak{h})_W)$ and $\nabla_{c}(S(\mathfrak{h}^*)_W^{\otimes})$ for a good  set $\mathcal{F}$ of parameter values  (see Definition \ref{definition00.1}),
thus demonstrating that for these parameter values, $P_{\textnormal{KZ}} \cong \Delta_{c} (S(\mathfrak{h})_W)$ in category $\mathcal{O}_c (\mathbb{Z}/n, \mathbb{C})$. 

In this section, we work with a general finite complex reflection group
$W$, with reflection representation $\mathfrak{h}$, and we specialise to category $\mathcal{O}_c (\mathbb{Z}/n, \mathbb{C})$ from Section \ref{chapter 5} onwards.  
We  demonstrate how $\mathcal{S}$ 
is related to $\textnormal{Hom}_{\mathbb{C}W} (S(\mathfrak{h}^*)_W, \textnormal{Har}^*)$ in three steps, beginning by showing that there is a one-to-one 
correspondence between elements of $\mathcal{S}$ 
and elements of  $\nabla_{c}(S(\mathfrak{h}^*)_W^{\otimes})$ that are killed by $S(\mathfrak{h})^W_{<0}$ and fixed by $W$.  

\subsection{The correspondence between $\mathcal{S}$
and $\nabla_{c}(S(\mathfrak{h}^*)_W^{\otimes})$   }  \label{subsection4.1}

We first observe the following simple consequence of induction and restriction. 

\begin{theorem3.12} \label{theorem3.12}
The elements of $\mathcal{S}$ are in one-to-one correspondence with the elements of  $\nabla_{c}(S(\mathfrak{h}^*)_W^{\otimes})$ that are killed by $S(\mathfrak{h})^W_{<0}$ and fixed 
by all elements of $W$. 
\end{theorem3.12}
\begin{proof}
As $\Delta_{c} (S(\mathfrak{h})_W) = H_c (W, \mathfrak{h}) \otimes_{S(\mathfrak{h}) \rtimes W} S(\mathfrak{h})_W$, 
it follows from induction and restriction that 
\begin{equation}  \label{equation3.3}
	\mathcal{S} \cong 
	\textnormal{Hom}_{S(\mathfrak{h}) \rtimes W} (S(\mathfrak{h})_W, \nabla_{c}(S(\mathfrak{h}^*)_W^{\otimes}) ).
\end{equation}
We recall from Subsection \ref{subsection00.8} that the action of $S(\mathfrak{h}) \rtimes W$ on $S(\mathfrak{h})$ is given by $(f \otimes w).g = f(w.g)$, and that $S(\mathfrak{h})_W = S(\mathfrak{h})/ (S(\mathfrak{h}).S(\mathfrak{h})^W_{<0})$.
We therefore have a map $\phi:  S(\mathfrak{h}) \rtimes W  \rightarrow S(\mathfrak{h})_W$ given by $\phi:  f \otimes w \mapsto f + (S(\mathfrak{h}).S(\mathfrak{h})^W_{<0})$, which is clearly onto, and whose kernel contains $(S(\mathfrak{h}).S(\mathfrak{h})^W_{<0})$ and $w-1$ for all $w \in W$. Furthermore, if $\sum^n_{i=1} f_i \otimes w_i \in \text{ker}\varphi$, then $\sum^n_{i=1} f_i  \in (S(\mathfrak{h}).S(\mathfrak{h})^W_{<0})    $ by the definition of $\varphi$, and therefore 
\begin{equation*} \label{equation3.134}
	\sum^n_{i=1} f_i \otimes w_i = \sum^n_{i=1} f_i \otimes (w_i - 1) + (\sum^n_{i=1} f_i ) \otimes 1 \in \langle S(\mathfrak{h})^W_{<0}, w-1 \rangle_{w \in W}.    
\end{equation*}
It follows that $S(\mathfrak{h})_W \cong (S(\mathfrak{h}) \rtimes W )/ \langle S(\mathfrak{h})^W_{<0}, w-1 \rangle_{w \in W}$ as $S(\mathfrak{h}) \rtimes W$-modules, and therefore \eqref{equation3.3} implies that $\mathcal{S}$ is in one-to-one correspondence with the elements of $\nabla_{c}(S(\mathfrak{h}^*)_W^{\otimes})$ that are killed by $S(\mathfrak{h})^W_{<0}$ and fixed by all elements of $W$, as claimed. 
\end{proof}

In order to understand what such elements look like, we go on to find an alternative expression for $\nabla_{c}(S(\mathfrak{h}^*)_W^{\otimes})$.  

\subsection{Another description of  $\nabla_{c}(S(\mathfrak{h}^*)_W^{\otimes})$}  \label{subsection4.2}
Recall from Subsection \ref{subsection00.4} that 
\begin{equation*} \label{equation3.20}
	\nabla_{c} (S(\mathfrak{h}^*)_W^{\otimes}) = \textnormal{Hom}_{S(\mathfrak{h}^*) \rtimes W} (H_c (W, \mathfrak{h}), S(\mathfrak{h}^*)_W^{\otimes})^{ln}. 
\end{equation*}
Denoting by $S(\mathfrak{h})^{\circledast}$ the elements of $S(\mathfrak{h})^* = \textnormal{Hom}_{\mathbb{C}} (S(\mathfrak{h}), \mathbb{C})$ that are nilpotent with respect to the natural action of $\mathfrak{h}$, we then have the following: 

\begin{theorem3.6}  \label{theorem3.6}
There is an isomorphism of vector spaces 
\begin{equation*}  \label{equation3.82} 
	\nabla_{c}(S(\mathfrak{h}^*)_W^{\otimes}) \cong S(\mathfrak{h})^{\circledast} \otimes_{\mathbb{C}} S(\mathfrak{h}^*)_W^{\otimes}.  
\end{equation*} 
The isomorphism gives $S(\mathfrak{h})^{\circledast} \otimes_{\mathbb{C}} S(\mathfrak{h}^*)_W^{\otimes}$ the structure of an 
$S(\mathfrak{h}) \rtimes W$-module, with $S(\mathfrak{h})$ acting naturally on the first factor, leaving the second factor fixed, and $W$ acting diagonally. 
\end{theorem3.6}

\begin{proof}
By the PBW-Theorem (see Subsection \ref{subsection00.2}), we have that $H_c (W, \mathfrak{h}) = (S(\mathfrak{h}^*) \rtimes W) \otimes_{\mathbb{C}} S(\mathfrak{h})$ as an $S(\mathfrak{h}^*) \rtimes W$-module, and it is therefore a consequence of induction and restriction that the maps defined by 
\begin{equation*}  \label{equation3.51}
	f \mapsto [f': y \mapsto f(1 \otimes y)]  \quad \textnormal{and} \quad g \mapsto [g': b \otimes y \mapsto b.g(y)] 
\end{equation*}
for all $y \in S(\mathfrak{h})$ and all $b \in S(\mathfrak{h}^*) \rtimes W$, provide mutually inverse isomorphisms between $\textnormal{Hom}_{S(\mathfrak{h}^*) \rtimes W} (H_c (W, \mathfrak{h}), S(\mathfrak{h}^*)_W^{\otimes})$ and $\textnormal{Hom}_{\mathbb{C}} (S(\mathfrak{h}), S(\mathfrak{h}^*)_W^{\otimes})$. It is furthermore easy to see that  these maps restrict to an isomorphism between $\textnormal{Hom}_{S(\mathfrak{h}^*) \rtimes W} (H_c (W, \mathfrak{h}), S(\mathfrak{h}^*)_W^{\otimes})^{ln}$ and $\textnormal{Hom}_{\mathbb{C}} (S(\mathfrak{h}), S(\mathfrak{h}^*)_W^{\otimes})^{ln}$, where $\textnormal{Hom}_{\mathbb{C}} (S(\mathfrak{h}), S(\mathfrak{h}^*)_W^{\otimes})^{ln}$ denotes the elements of the space $\textnormal{Hom}_{\mathbb{C}} (S(\mathfrak{h}), S(\mathfrak{h}^*)_W^{\otimes})$ that are locally nilpotent with respect to the natural action of $\mathfrak{h}$. This isomorphism gives $\textnormal{Hom}_{\mathbb{C}} (S(\mathfrak{h}), S(\mathfrak{h}^*)_W^{\otimes})^{ln}$ the structure of an $S(\mathfrak{h}) \rtimes W$-module, with $S(\mathfrak{h})$ acting naturally, and $W$ acting diagonally in the sense that $(w.f)(y) = w.f(w^{-1}.y) $ for all $w \in W$, $y \in S(\mathfrak{h})$ and $f \in \textnormal{Hom}_{\mathbb{C}} (S(\mathfrak{h}), S(\mathfrak{h}^*)_W^{\otimes})^{ln}$. 

Next, let $\mathscr{B}$ be a homogeneous basis  for $S(\mathfrak{h})$. Observe that an element $f \in S(\mathfrak{h})^*$ is contained in $S(\mathfrak{h})^{\circledast}$ if and only if there exists an integer $n$ such that $f(S(\mathfrak{h})^{>n}) = 0$. As $\mathfrak{h}$ is finite-dimensional,  it therefore follows that the set of dual vectors $\mathscr{B}^{\vee} = \{v^{\vee} : v \in \mathscr{B} \}$ corresponding to $\mathscr{B}$, is a basis for $S(\mathfrak{h})^{\circledast}$. We define a map $\alpha: \textnormal{Hom}_{\mathbb{C}} (S(\mathfrak{h}), S(\mathfrak{h}^*)_W^{\otimes})^{ln} \to S(\mathfrak{h})^{\circledast} \otimes_{\mathbb{C}} S(\mathfrak{h}^*)_W^{\otimes}$ by 
\begin{equation*} \label{equation3.64}
	\alpha: f \mapsto \sum_{v \in \mathscr{B}} v^{\vee} \otimes f(v), 
\end{equation*} 
observing that only finitely many of the terms in the sum are non-zero. We also define a map $\beta:  S(\mathfrak{h})^{\circledast} \otimes_{\mathbb{C}} S(\mathfrak{h}^*)_W^{\otimes} \to \textnormal{Hom}_{\mathbb{C}} (S(\mathfrak{h}), S(\mathfrak{h}^*)_W^{\otimes})^{ln}$ by 
\begin{equation*} \label{equation3.66}
	\beta: g \otimes m \mapsto [y \mapsto g(y)m]. 
\end{equation*} 
We show that $\alpha$ and $\beta$ are inverse to each other. 
First, let $f \in \textnormal{Hom}_{\mathbb{C}} (S(\mathfrak{h}), S(\mathfrak{h}^*)_W^{\otimes})^{ln}$. Then, for every $w \in \mathscr{B}$, 
\begin{equation*} \label{equation3.67}
	(\beta \circ \alpha)(f)(w) = \sum_{v \in \mathscr{B}} v^{\vee} (w) f(v) = f(w),  
\end{equation*} 
so that $(\beta \circ \alpha)(f) = f$ and $\beta \circ \alpha = \text{id}_{\textnormal{Hom}_{\mathbb{C}} (S(\mathfrak{h}), S(\mathfrak{h}^*)_W^{\otimes})^{ln}}$. 
Next, suppose again that $w \in \mathscr{B}$ and that 
$m \in S(\mathfrak{h}^*)_W^{\otimes}$. Then 
\begin{equation*} \label{equation3.68}
	(\alpha \circ \beta)(w^{\vee} \otimes m) = \sum_{v \in \mathscr{B}} v^{\vee} \otimes (w^{\vee} (v) m) = w^{\vee} \otimes m,   
\end{equation*} 
and since $\mathscr{B}^{\vee}$ is a basis for $S(\mathfrak{h})^{\circledast}$ it follows that $\alpha \circ \beta = \text{id}_{S(\mathfrak{h})^{\circledast} \otimes_{\mathbb{C}} S(\mathfrak{h}^*)_W^{\otimes}}$. 
Therefore $\textnormal{Hom}_{\mathbb{C}} (S(\mathfrak{h}), S(\mathfrak{h}^*)_W^{\otimes})^{ln} $ and $S(\mathfrak{h})^{\circledast} \otimes_{\mathbb{C}} S(\mathfrak{h}^*)_W^{\otimes}$ are isomorphic as vector spaces. This isomorphism gives $S(\mathfrak{h})^{\circledast} \otimes_{\mathbb{C}} S(\mathfrak{h}^*)_W^{\otimes}$ the structure of a $S(\mathfrak{h}) \rtimes W$-module, with $S(\mathfrak{h})$ acting naturally on the first factor, leaving the second factor fixed and $W$ acting diagonally. This completes the proof. 
\end{proof}

In line with Proposition \ref{theorem3.12}, we thus want to understand the elements of  $S(\mathfrak{h})^{\circledast} \otimes_{\mathbb{C}} S(\mathfrak{h}^*)_W^{\otimes}$ that are killed by $S(\mathfrak{h})^W_{<0}$ and fixed by all elements of $W$. To this end, we show in the next 
subsection that these elements of
$S(\mathfrak{h})^{\circledast} \otimes_{\mathbb{C}} S(\mathfrak{h}^*)_W^{\otimes}$ 
are in one-to-one
correspondence with the elements of $\textnormal{Hom}_{\mathbb{C}W} (S(\mathfrak{h}^*)_W, \textnormal{Har}^*)$. First, we make the following observation, which will be needed in Section \ref{chapter 5} to identify isomorphisms in $\mathcal{S}$ in the setting of category $\mathcal{O}_c (\mathbb{Z}/n, \mathbb{C})$. 

\begin{theorem3.14} \label{theorem3.14}
Under the one-to-one correspondence between $\mathcal{S}$ and the elements of $S(\mathfrak{h})^{\circledast} \otimes_{\mathbb{C}} S(\mathfrak{h}^*)_W^{\otimes}$ that are killed by $S(\mathfrak{h})^W_{<0}$ and fixed by all elements of $W$, established in Propositions \ref{theorem3.12} and \ref{theorem3.6}, the isomorphisms in $\mathcal{S}$ correspond to those elements that generate $S(\mathfrak{h})^{\circledast} \otimes_{\mathbb{C}} S(\mathfrak{h}^*)_W^{\otimes}$  as an $H_c (W, \mathfrak{h})$-module.  
\end{theorem3.14}

\begin{proof}
Observe that Corollary \ref{theorem00.11} implies that an element of $\mathcal{S}$ is an isomorphism precisely when it is an epimorphism. 
From Proposition \ref{theorem3.12}, it therefore follows that the isomorphisms in $\mathcal{S}$ correspond precisely to the elements of $\nabla_{c}(S(\mathfrak{h}^*)_W^{\otimes})$ that are killed by $S(\mathfrak{h})^W_{<0}$, fixed by all elements of $W$ and generate $\nabla_{c}(S(\mathfrak{h}^*)_W^{\otimes})$ as an $H_c (W, \mathfrak{h})$-module. The claim then follows from Proposition \ref{theorem3.6}, where $S(\mathfrak{h})^{\circledast} \otimes_{\mathbb{C}} S(\mathfrak{h}^*)_W^{\otimes}$ inherits its $H_c (W, \mathfrak{h})$-module structure from $\nabla_{c}(S(\mathfrak{h}^*)_W^{\otimes})$. 
\end{proof}

\subsection{The elements of $S(\mathfrak{h})^{\circledast} \otimes_{\mathbb{C}} S(\mathfrak{h}^*)_W^{\otimes}$ that 
are killed by $S(\mathfrak{h})^W_{<0}$ and fixed by $W$}  \label{subsection4.3}

The final step in our simplification of $\mathcal{S}$ is the following: 

\begin{theorem3.13}  \label{theorem3.13}
The elements of $S(\mathfrak{h})^{\circledast} \otimes_{\mathbb{C}} S(\mathfrak{h}^*)_W^{\otimes}$ that are killed by 
$S(\mathfrak{h})^W_{<0}$ and are fixed by all elements of $W$, are in one-to-one correspondence with the elements of the space $\textnormal{Hom}_{\mathbb{C}W} (S(\mathfrak{h}^*)_W, \textnormal{Har}^*)$. 

\end{theorem3.13}

\begin{proof}
We first identify the elements of $S(\mathfrak{h})^{\circledast} \otimes_{\mathbb{C}} S(\mathfrak{h}^*)_W^{\otimes}$ that are killed by $S(\mathfrak{h})^W_{<0}$. Since $S(\mathfrak{h})$ is only acting on the first factor of 
$S(\mathfrak{h})^{\circledast} \otimes_{\mathbb{C}} S(\mathfrak{h}^*)_W^{\otimes}$, it is sufficient to identify the elements of $S(\mathfrak{h})^{\circledast}$ that are
killed by $S(\mathfrak{h})^W_{<0}$. Denote these elements by $X$, and observe that $X = \{ f \in S(\mathfrak{h})^{\circledast} : f((S(\mathfrak{h}).S(\mathfrak{h})^W_{<0})) = 0 \}$ is a $W$-submodule of $S(\mathfrak{h})^{\circledast}$. We recall from Subsection \ref{subsection00.9} that we have a decomposition $S(\mathfrak{h}) = (S(\mathfrak{h}).S(\mathfrak{h})^W_{<0}) \oplus \textnormal{Har} $
of $S(\mathfrak{h})$ as a $W$-module, where Har denotes the set of harmonic polynomials in $S(\mathfrak{h})$.
The maps 
\begin{equation*}  \label{equation3.85} 
	X \to \textnormal{Har}^*, \quad f \mapsto f |_{\textnormal{Har}} \quad \textnormal{and} \quad \textnormal{Har}^* \to X, \quad g \mapsto 0 \oplus g
\end{equation*} 
are readily seen to be inverse to each other, and it is clear that these maps respect the action of $W$. 
This means that there is a one-to-one correspondence between elements of $S(\mathfrak{h})^{\circledast} \otimes_{\mathbb{C}} S(\mathfrak{h}^*)_W^{\otimes}$ that are killed by $S(\mathfrak{h})^W_{<0}$ and fixed by $W$, and elements of  $\textnormal{Har}^* \otimes_{\mathbb{C}} S(\mathfrak{h}^*)_W^{\otimes}$ that are fixed by $W$. 

In order to identify such elements, we 
recall the definition of the duality $\langle -,- \rangle: S(\mathfrak{h}^*) \times S(\mathfrak{h}) \to \mathbb{C}$ and the Jacobian $J \in S(\mathfrak{h})$  of $W$ given in Subsection \ref{subsection00.9}. In terms of these notions, it is proved in Theorem 4.25(2) of \cite{Bro} that there is an isomorphism of  $W$-modules $S(\mathfrak{h}^*)_W^{\otimes} \cong S(\mathfrak{h}^*)^*_W$, given by
\begin{equation}  \label{equation3.88} 
	 \phi :(a + \langle S(\mathfrak{h}^*)^W_{>0} \rangle )  \otimes_{\mathbb{C}} 1 \mapsto [(b + \langle S(\mathfrak{h}^*)^W_{>0} \rangle) \mapsto \langle ab, J \rangle ].
\end{equation} 
The elements of $\textnormal{Har}^* \otimes_{\mathbb{C}} S(\mathfrak{h}^*)_W^{\otimes}$ that are fixed by $W$ are therefore in one-to-one correspondence with the elements of $\textnormal{Har}^* \otimes_{\mathbb{C}} S(\mathfrak{h}^*)^*_W$ that are fixed by $W$, and we denote these by $(\textnormal{Har}^* \otimes_{\mathbb{C}} S(\mathfrak{h}^*)^*_W)^W$. 

Next, we observe that the space $\textnormal{Hom}_{\mathbb{C}} (S(\mathfrak{h}^*)_W, \textnormal{Har}^*)$ becomes a $W$-module by defining $(w.f)(u) = w.f(w^{-1}.u)$ for all $w \in W$, $u \in S(\mathfrak{h}^*)_W$ and $f \in \textnormal{Hom}_{\mathbb{C}} (S(\mathfrak{h}^*)_W, \textnormal{Har}^*)$. We define a map $\beta$ from $\textnormal{Har}^* \otimes_{\mathbb{C}} S(\mathfrak{h}^*)^*_W$ to $ \textnormal{Hom}_{\mathbb{C}} (S(\mathfrak{h}^*)_W, \textnormal{Har}^*) $ by 
\begin{equation*}  \label{equation3.90} 
	 \beta: f  \otimes_{\mathbb{C}} g \mapsto [u \mapsto g(u).f], 
\end{equation*} 
and it is easy to check that it is a $W$-module homomorphism. Furthermore, if we let $\mathscr{B}$ denote a basis for $S(\mathfrak{h}^*)_W$ and let $\mathscr{B}^{\vee} = 
\{ v^{\vee} : v \in \mathscr{B} \}$ be the corresponding dual basis of $S(\mathfrak{h}^*)_W^*$,  then the map $\alpha$ from $\textnormal{Hom}_{\mathbb{C}} (S(\mathfrak{h}^*)_W, \textnormal{Har}^*)$ to $\textnormal{Har}^* \otimes_{\mathbb{C}} S(\mathfrak{h}^*)^*_W$ defined by  
\begin{equation*}  \label{equation3.95} 
	 \alpha: f \mapsto \sum_{ v \in \mathscr{B} } f(v) \otimes v^{\vee}
\end{equation*} 
is readily seen to be inverse to the map $\beta$, so that $\textnormal{Hom}_{\mathbb{C}} (S(\mathfrak{h}^*)_W, \textnormal{Har}^*)$ and $\textnormal{Har}^* \otimes_{\mathbb{C}} S(\mathfrak{h}^*)^*_W$ are isomorphic as $W$-modules. 
This isomorphism restricts to an isomorphism of $ (\textnormal{Har}^* \otimes_{\mathbb{C}} S(\mathfrak{h}^*)^*_W)^W$ and $\textnormal{Hom}_{\mathbb{C}} (S(\mathfrak{h}^*)_W, \textnormal{Har}^*)^W$, and by definition of the action of $W$ on $\textnormal{Hom}_{\mathbb{C}} (S(\mathfrak{h}^*)_W, \textnormal{Har}^*)$, we see that $\textnormal{Hom}_{\mathbb{C}} (S(\mathfrak{h}^*)_W, \textnormal{Har}^*)^W$ is just $\textnormal{Hom}_{\mathbb{C}W} (S(\mathfrak{h}^*)_W, \textnormal{Har}^*)$, 
which completes the proof. 
\end{proof}

We now put the pieces together. 

\subsection{The correspondence between $\mathcal{S}$ and $\textnormal{Hom}_{\mathbb{C}W} (S(\mathfrak{h}^*)_W, \textnormal{Har}^*)$ } \label{subsection4.4}

Propositions \ref{theorem3.12}, \ref{theorem3.6} and \ref{theorem3.13} combine to give the following:

\begin{theorem000.1} \label{theorem000.1}
There is a one-to-one correspondence between the elements of $\mathcal{S}$ and the elements of $\textnormal{Hom}_{\mathbb{C}W} (S(\mathfrak{h}^*)_W, \textnormal{Har}^*)$. 
\end{theorem000.1}

In the next section, we analyse the correspondence of Theorem A in the setting of category $\mathcal{O}_c (\mathbb{Z}/n, \mathbb{C})$, to find a set of parameter values $c$ for which we have that $P_{\textnormal{KZ}} \cong \Delta_{c}(S(\mathfrak{h})_W)$.

\section{$P_{\textnormal{KZ}}$ in category $\mathcal{O}_{\boldsymbol{c}} (\mathbb{Z}/n, \mathbb{C})$} \label{chapter 5}

In this section, we show that $P_{\textnormal{KZ}} \cong \Delta_{\boldsymbol{c}}(S(\mathfrak{h})_W)$ in category $\mathcal{O}_{\boldsymbol{c}} (\mathbb{Z}/n, \mathbb{C})$
for a good set $\mathcal{F}$ of parameter values $\boldsymbol{c}$. We do this by using the results from Section \ref{chapter 3} to find an isomorphism in $\mathcal{S} = \text{Hom}_{H_{\boldsymbol{c}} (W, \mathfrak{h})}(\Delta_{\boldsymbol{c}} (S(\mathfrak{h})_W), \nabla_{\boldsymbol{c}}(S(\mathfrak{h}^*)_W^{\otimes}) )$ in the setting of category $\mathcal{O}_{\boldsymbol{c}} (\mathbb{Z}/n, \mathbb{C})$. 
Specifically, we describe the set $\text{Hom}_{\mathbb{C} W}(S(\mathfrak{h}^*)_W, \textnormal{Har}^*)$, and compute the lifting of its elements to the module $M := S(\mathfrak{h})^{\circledast} \otimes_{\mathbb{C}} S(\mathfrak{h}^*)_W^{\otimes}$. As described in Proposition  \ref{theorem3.14}, the lifting of such an element to $\mathcal{S}$ is an isomorphism,
 if and only if the lifting of the same element to $M$ generates $M$ as 
an $H_{\boldsymbol{c}}  (\mathbb{Z}/n, \mathbb{C})$-module. In order to determine for which parameter values $\boldsymbol{c}$ there exists an element of  $\text{Hom}_{\mathbb{C} W}(S(\mathfrak{h}^*)_W, \textnormal{Har}^*)$ that lifts to an isomorphism in  $\mathcal{S}$, we therefore need to understand the action of $H_{\boldsymbol{c}} (W, \mathfrak{h})$ on $M$. 
The action of the subalgebra $S(\mathfrak{h}) \rtimes W$ on $M$ was determined in Proposition  \ref{theorem3.6}, and the action of $S(\mathfrak{h}^*)$ that $M$ inherits from $\nabla_{\boldsymbol{c}}(S(\mathfrak{h}^*)_W^{\otimes})$ is computed in Proposition  \ref{theorem5.6}.
 It turns out that in the setting of category $\mathcal{O}_{\boldsymbol{c}} (\mathbb{Z}/n, \mathbb{C})$, it is enough to study the lifting of a single distinguished element of $\text{Hom}_{\mathbb{C} W}(S(\mathfrak{h}^*)_W, \textnormal{Har}^*)$, and that the question for which parameter values $\boldsymbol{c}$ this lifting is an isomorphism in  $\mathcal{S}$ reduces to determining for which values of the parameter $\boldsymbol{c}$ a 
certain family of matrices are all non-singular. 

We start by describing  $\text{Hom}_{\mathbb{C} W}(S(\mathfrak{h}^*)_W, \textnormal{Har}^*)$, and calculate the liftings of its elements to $M = S(\mathfrak{h})^{\circledast} \otimes_{\mathbb{C}} S(\mathfrak{h}^*)_W^{\otimes}$ in the setting of category $\mathcal{O}_{\boldsymbol{c}} (\mathbb{Z}/n, \mathbb{C})$.

\subsection{Lifting elements of $\text{Hom}_{\mathbb{C} W}(S(\mathfrak{h}^*)_W, \textnormal{Har}^*)$  to $S(\mathfrak{h})^{\circledast} \otimes_{\mathbb{C}} S(\mathfrak{h}^*)_W^{\otimes}$} \label{subsection5.1}

Recall from Subsection \ref{subsection00.9} that $S(\mathfrak{h}^*).S(\mathfrak{h}^*)^W_{>0} = \langle x^n \rangle $ and $\textnormal{Har} = \textnormal{Span}_{\mathbb{C}} \{\xi^i : 0 \leq i \leq n-1  \}$. Let $\{ (\xi^i)^* : 0 \leq i \leq n-1 \}$ be the dual basis of $\textnormal{Har}^*$ corresponding to the basis $\{ \xi^i : 0 \leq i \leq n-1 \}$ of Har, and for every $y \in S(\mathfrak{h}^*) $, let  $\bar{y} = y + \langle x^n \rangle \in S(\mathfrak{h}^*)_W$. 

\begin{theorem5.17}  \label{theorem5.17}
A basis for $\textnormal{Hom}_{\mathbb{C}W} (S(\mathfrak{h}^*)_W, \textnormal{Har}^*)$ is given by $\{ \psi_i : 0 \leq i \leq n-1  \}$, where $\psi_i$ is the unique linear map such that
\begin{equation*} \label{equation5.201}
	\psi_i (\bar{x}^j) = 
	\begin{cases}
	(\xi^i)^* & \textnormal{if} \,\, i=j  \\
	0 & \textnormal{otherwise}. 
	\end{cases}
\end{equation*} 
\end{theorem5.17}

\begin{proof}
As graded $\mathbb{Z}/n$-modules, both $S(\mathfrak{h}^*)_W$ and $\textnormal{Har}^*$ are concentrated in degrees between $0$ and $n-1$. As both $(S(\mathfrak{h}^*)_W)_i$ and $(\textnormal{Har}^*)_i$ are isomorphic to $E_i$ as $\mathbb{Z}/n$-modules for $0 \leq i \leq n-1$ (cf. Section \ref{subsection00.3}), the result follows. 
\end{proof}

Next, we lift the $\psi_i \in \text{Hom}_{\mathbb{C} W}(S(\mathfrak{h}^*)_W, \textnormal{Har}^*)$ to the module $M = S(\mathfrak{h})^{\circledast} \otimes_{\mathbb{C}} S(\mathfrak{h}^*)_W^{\otimes}$ by analysing the isomorphism of Propositions \ref{theorem3.13}. We denote by  $ \mathscr{B}^{\vee} = \{ (\xi^i)^{\vee}: i \in \mathbb{N} \}$
the set of dual vectors in $S(\mathfrak{h})^*$ corresponding to the basis $\mathscr{B}= \{ \xi^i: i \in \mathbb{N} \}$ of $S(\mathfrak{h})$. Observe that $\mathscr{B}^{\vee}$ is a basis for $S(\mathfrak{h})^{\circledast}$. 

\begin{theorem5.61} \label{theorem5.61}
The element $\psi_i \in \textnormal{Hom}_W(S(\mathfrak{h}^*)_W, \textnormal{Har}^*)$, where $0 \leq i \leq n-1$, lifts to the element 
\begin{equation*} \label{equation5.23}
	 \psi^M_i =\frac{1}{(n-1)!} (\xi^i)^{\vee} \otimes \bar{x}^{(n-1)-i}  \otimes 1
\end{equation*}
of $M = S(\mathfrak{h})^{\circledast} \otimes_{\mathbb{C}} S(\mathfrak{h}^*)_W^{\otimes}$. 
\end{theorem5.61}

\begin{proof}
Let $0 \leq i \leq n-1$. We begin by lifting $\psi_i$ to $\textnormal{Har}^* \otimes_{\mathbb{C}} S(\mathfrak{h}^*)^*_W$. It was seen in the proof of Proposition \ref{theorem3.13} that, if we let $\{(\bar{x}^j)^{\vee} : 0 \leq j \leq n-1 \}$ be the dual basis of $S(\mathfrak{h}^*)^*_W$ corresponding to the basis $\{\bar{x}^j : 0 \leq j \leq n-1 \}$ of $S(\mathfrak{h}^*)_W$, then $\psi_i$ lifts to the element 
\begin{equation} \label{equation5.25}
	 \sum^{n-1}_{j = 0} \psi_i (\bar{x}^j) \otimes (\bar{x}^j)^{\vee} = (\xi^i)^* \otimes (\bar{x}^i)^{\vee}
\end{equation}
of $\textnormal{Har}^* \otimes_{\mathbb{C}} S(\mathfrak{h}^*)^*_W$. 
Next, we lift this element to $\textnormal{Har}^* \otimes_{\mathbb{C}} S(\mathfrak{h}^*)_W^{\otimes}$. This amounts to finding the preimage of $(\bar{x}^i)^{\vee} $ under the isomorphism $\phi$ defined in \eqref{equation3.88}. Since it was seen in Subsection \ref{subsection00.9} that the Jacobian $J = \xi^{n-1}$, it follows that $\phi^{-1}((\bar{x}^i)^{\vee}) = \bar{a} \otimes 1$, where $\bar{a} \in  S(\mathfrak{h}^*)_W = \mathbb{C}[x] / \langle x^n \rangle$ is such that $ \langle ab, \xi^{n-1} \rangle = (\bar{x}^i)^{\vee} (\bar{b}) $
for every $\bar{b} \in \mathbb{C}[x]/  \langle x^n \rangle$. Recalling from Subsection \ref{subsection00.9} how $\langle -, - \rangle$ behaves with respect to a basis and its dual basis, we see that this forces 
\begin{equation*} \label{equation5.27}
	 \phi^{-1}((\bar{x}^i)^{\vee}) = \frac{1}{(n-1)!} (\bar{x}^{(n -1) -i} \otimes 1).  
\end{equation*}
Hence, the element in \eqref{equation5.25}, and thus $\psi_i$, lifts to the element 
\begin{equation} \label{equation5.28}
	\frac{1}{(n-1)!} (\xi^i)^*  \otimes \bar{x}^{(n -1) -i} \otimes 1
\end{equation}
of $\textnormal{Har}^* \otimes_{\mathbb{C}} S(\mathfrak{h}^*)_W^{\otimes}$. 
Finally, recall from Subsection \ref{subsection00.8} that $S(\mathfrak{h}).S(\mathfrak{h})^W_{>0} = \mathbb{C}[\xi].  \mathbb{C}[\xi]^{\mathbb{Z}/n}_{>0} = \langle \xi^n \rangle$. Therefore, the decomposition $S(\mathfrak{h}) = (S(\mathfrak{h}).S(\mathfrak{h})^W_{>0}) \oplus \textnormal{Har}$ together with the discussion at the beginning of the proof Proposition \ref{theorem3.13} implies that the element $(\xi^i)^*$ of $\textnormal{Har}^*$ lifts to the element $(\xi^i)^{\vee}$ in $S(\mathfrak{h})^{\circledast}$. It follows that the element of \eqref{equation5.28}, and hence $\psi_i$, lifts to the element 
\begin{equation*} \label{equation5.29}
	 \psi^M_i =  \frac{1}{(n-1)!} (\xi^i)^{\vee} \otimes \bar{x}^{(n-1)-i} \otimes 1
\end{equation*}
of $M = S(\mathfrak{h})^{\circledast} \otimes_{\mathbb{C}} S(\mathfrak{h}^*)_W^{\otimes}$. This completes the proof. 
\end{proof}

We make the following definition in order to simplify the notation.

\begin{definition5.4} \label{definition5.4}

\begin{itemize}

\item[(i)]
For $i, j \geq 0$ we define 
\begin{equation*}
	v_{i,j} = (\xi^i)^{\vee} \otimes \bar{x}^j \otimes 1.  
\end{equation*}
\item[(ii)]

For $\boldsymbol{t} = (t_0, t_1, \dots, t_{n-1}) \in \mathbb{C}^n$, we let 
\begin{equation*} \label{equation5.202}
	 \psi_{\boldsymbol{t}} =  \sum_{i = 0}^{n-1} t_i \psi_i \quad \textnormal{and}  \quad  \psi^M_{\boldsymbol{c}, \boldsymbol{t}} = \sum_{i = 0}^{n-1} t_i \psi^M_i = \sum_{i = 0}^{n-1} \frac{t_i}{(n-1)!} v_{i, (n-1)-i} . 
\end{equation*}
\end{itemize}
\end{definition5.4}

We want to understand for which parameter values $\boldsymbol{c}$ there exists an $n$-tuple $\boldsymbol{t}$ such that $\psi^M_{\boldsymbol{c}, \boldsymbol{t}}$ generates $M$ as an $H_{\boldsymbol{c}} (\mathbb{Z}/n, \mathbb{C})$-module. In order to do this we need to understand the action of $H_{\boldsymbol{c}} (\mathbb{Z}/n, \mathbb{C})$ on $M$. We begin by noting the following.

\begin{theorem5.9}  \label{theorem5.9}

A basis for $M = S(\mathfrak{h})^{\circledast} \otimes_{\mathbb{C}} S(\mathfrak{h}^*)_W^{\otimes}$ is given by $\{ v_{i,j} : i \geq 0, 0 \leq j \leq n-1 \}$. 
Furthermore, $M$ is concentrated in degrees greater than or equal to $1-n$, and for each $k \geq 1-n$, $\{ v_{i,j} : i \geq 0, 0 \leq j \leq n-1, i+j = k + (n -1)  \}$
 is a basis for the kth graded component $M_k$ of $M$. In particular 
 \begin{equation*}  \label{equation5.4}
	\textnormal{dim}_{\mathbb{C}} M_k = 
	\begin{cases}
	n + k & \textnormal{if }1-n \leq  k < 0, \\
	n & \textnormal{if } k \geq 0. 
	\end{cases}  
\end{equation*}
\end{theorem5.9}

\begin{proof}
As $\xi$ has degree $-1$, the space $S(\mathfrak{h})^{\circledast} = \mathbb{C}[\xi]^{\circledast}$ is concentrated in non-negative degrees, and as $x$ has degree $1$, $S(\mathfrak{h}^*)_W = \mathbb{C}[x] / \langle x^n \rangle $ is concentrated in degrees between 0 and $n-1$. Since $\text{det}^{-1}_\mathfrak{h} \langle n-1 \rangle $ is concentrated in degree $1-n$ and $S(\mathfrak{h}^*)_W^{\otimes} = S(\mathfrak{h}^*)_W \otimes \text{det}^{-1}_\mathfrak{h} \langle n-1 \rangle  $, the claim then follows easily. 
\end{proof}
We proceed to describe the action of $H_{\boldsymbol{c}} (\mathbb{Z}/n, \mathbb{C})$ on this basis. 
 
\subsection{The action of $H_{\boldsymbol{c}} (\mathbb{Z}/n, \mathbb{C})$ on $S(\mathfrak{h})^{\circledast} \otimes_{\mathbb{C}} S(\mathfrak{h}^*)_W^{\otimes}$ } \label{subsection5.2} 
 
The action of $H_{\boldsymbol{c}} (\mathbb{Z}/n, \mathbb{C})$ on $M$ is inherited from the isomorphism $M = S(\mathfrak{h})^{\circledast} \otimes_{\mathbb{C}} S(\mathfrak{h}^*)_W^{\otimes} \cong \nabla_{\boldsymbol{c}}(S(\mathfrak{h}^*)_W^{\otimes})$
of Proposition \ref{theorem3.6}. It was noted in that proposition that $S(\mathfrak{h})$ acts 
naturally on the first factor of $M$, leaving the second factor fixed, and that $W$ acts diagonally. In order to translate the action of $S(\mathfrak{h}^*)$ on $\nabla_{\boldsymbol{c}}(S(\mathfrak{h}^*)_W^{\otimes})$ to $M$, we will need the following result. Recall the definition of the idempotents $\varepsilon_i$ from Subsection \ref{subsection00.3}. 

\begin{theorem5.15} \label{theorem5.15}
The following holds in $H_{\boldsymbol{c}} (\mathbb{Z}/n, \mathbb{C})$:
\begin{itemize}
\item[(i)] For all $i \in \mathbb{Z}$ we have that  $\xi \varepsilon_i = \varepsilon_{i-1} \xi$. 

\item[(ii)]
For all $j \geq 0$ we have that 
\begin{equation*} \label{equation5.34}
	\xi^{j} x = x \xi^{j} + j \xi^{j-1} + \Big[\sum^{n-1}_{i=0} (c_{i+j} - c_i) \varepsilon_i \Big] \xi^{j-1}. 	
\end{equation*}	
\end{itemize}

\end{theorem5.15}

\begin{proof}
Suppose that $i \in \mathbb{Z}$. As $w^{-1} \xi w = w^{-1}.\xi$ so that $\xi w = \textnormal{det}^{-1}(w) w \xi$, it follows that
\begin{equation*} \label{equation5.047}
	\xi \varepsilon_i = \xi \Big(\frac{1}{n}\sum_{w \in \mathbb{Z}/n} \text{det}(w)^i w \Big) =  \frac{1}{n} \sum_{w \in \mathbb{Z}/n} \text{det}(w)^i \xi w = \Big(\frac{1}{n} \sum_{w \in \mathbb{Z}/n} \text{det}(w)^{i-1} w \Big) \xi = \varepsilon_{i-1} \xi. 
\end{equation*}
which proves (i). 
We prove part (ii) by induction on $j \geq 0$. It is clearly true for $j = 0$ (we let $\xi^{-1} = 0$), and it follows from \eqref{equation00.3}, \eqref{equation00.4} and the fact that $c_i = nk_i$ that it holds for $j = 1$.  Suppose that the claim holds true for some $j \geq 1$. Then, 
%{\allowdisplaybreaks
\begin{IEEEeqnarray}{rCl} 
 \xi^{j+1} x & = & \xi (\xi^{j} x)    \label{equation5.41}  \nonumber \\
& = & \xi \Big(x \xi^{j} + j \xi^{j-1} + \Big[\sum^{n-1}_{i=0} (c_{i+j} - c_i) \varepsilon_i \Big] \xi^{j-1}\Big)    \label{equation5.42} \\
& = & (\xi x) \xi^{j} + j \xi^{j} + \Big[\sum^{n-1}_{i=0} (c_{i+j} - c_i) \xi \varepsilon_i \Big] \xi^{j-1} \label{equation5.43} \nonumber \\ 
& = & \Big(1 + x\xi + \sum^{n-1}_{i=0} (c_{i+1} - c_i) \varepsilon_i \Big) \xi^{j} + j \xi^{j} + \Big[\sum^{n-1}_{i=0} (c_{i+j} - c_i) \varepsilon_{i-1} \Big] \xi^{j} \label{equation5.44}  \\
& = & x\xi^{j+1} + (j+1) \xi^{j} + \Big[ \sum^{n-1}_{i=0} ((c_{i+1} - c_i) + (c_{i + j+1} - c_{i+1})) \varepsilon_i \Big] \xi^{j} \label{equation5.45}  \\
& = & x \xi^{j+1} + (j+1) \xi^{j} + \Big[\sum^{n-1}_{i=0} (c_{i+j+1} - c_i) \varepsilon_i \Big] \xi^{j} \label{equation5.46} \nonumber,
\end{IEEEeqnarray}%} 
where equation \eqref{equation5.42} follows from our inductive hypothesis, equation \eqref{equation5.44} follows from our inductive hypothesis and part (i), and equation \eqref{equation5.45} follows from the fact that $\varepsilon_k =  \varepsilon_l$ and $c_k = c_l$ whenever $k \equiv l \pmod n$. The claim thus follows by induction, and this completes the proof. \qedhere
\end{proof}
With the help of this result, we are now able to compute the action of $S(\mathfrak{h}^*)$ on $M = S(\mathfrak{h})^{\circledast} \otimes_{\mathbb{C}} S(\mathfrak{h}^*)_W^{\otimes}$. We summarise the action of $H_{\boldsymbol{c}} (\mathbb{Z}/n, \mathbb{C})$ on $M$ in the following claim:

\begin{theorem5.6} \label{theorem5.6}
Let $i \geq 0$ and $0 \leq j \leq n-1$. The action of the generators $\xi$, $s$ and $x$ of $H_{\boldsymbol{c}} (\mathbb{Z}/n, \mathbb{C})$ on the basis elements of  $S(\mathfrak{h})^{\circledast} \otimes_{\mathbb{C}} S(\mathfrak{h}^*)_W^{\otimes}$ is given by
\begin{equation*} \label{equation5.048}
	\xi.v_{i,j} = v_{i-1, j}, \quad s.v_{i,j} = q^{-(i + j +1)} v_{i, j} \quad \textnormal{and} \quad x.v_{i,j} = v_{i, j+1} + (i+1 + c_{i + j +2} - c_{j+1}) v_{i+1, j}
\end{equation*}
where we let $v_{-1, j} = 0$.

\end{theorem5.6}

\begin{proof}
It was seen in Proposition \ref{theorem3.6} that $S(\mathfrak{h}) = \mathbb{C}[\xi]$ acts naturally on the first factor of $M = S(\mathfrak{h})^{\circledast} \otimes_{\mathbb{C}} S(\mathfrak{h}^*)_W^{\otimes}$, leaving the second factor fixed, and that $W = \mathbb{Z}/n$ acts diagonally, and this gives the action of $\xi$ and $s$. For the action of $x$, let  
\begin{equation*} \label{equation5.50}
	\phi_1: \textnormal{Hom}_{S(\mathfrak{h}^*) \rtimes W} (H_c (W, \mathfrak{h}), S(\mathfrak{h}^*)_W^{\otimes})^{ln} \to \textnormal{Hom}_{\mathbb{C}} (S(\mathfrak{h}), S(\mathfrak{h}^*)_W^{\otimes})^{ln}, 
\end{equation*} 
and 
\begin{equation*} \label{equation5.51}
	\phi_2: \textnormal{Hom}_{\mathbb{C}} (S(\mathfrak{h}), S(\mathfrak{h}^*)_W^{\otimes})^{ln} \to S(\mathfrak{h})^{\circledast} \otimes_{\mathbb{C}} S(\mathfrak{h}^*)_W^{\otimes}
\end{equation*} 
denote the isomorphisms described in the proof of Proposition \ref{theorem3.6}. The action of $S(\mathfrak{h}^*) = \mathbb{C}[x]$ on the module $M = S(\mathfrak{h})^{\circledast} \otimes_{\mathbb{C}} S(\mathfrak{h}^*)_W^{\otimes} $ is inherited from the action of $S(\mathfrak{h}^*)$ on the module  $\nabla_{\boldsymbol{c}}(S(\mathfrak{h}^*)_W^{\otimes}) = \textnormal{Hom}_{S(\mathfrak{h}^*) \rtimes W} (H_c (W, \mathfrak{h}), S(\mathfrak{h}^*)_W^{\otimes})^{ln}$, in accordance with the composition of these isomorphisms. Using the definitions of $\phi_1$ and $\phi_2$, we therefore have that for $i \geq 0$ and $0 \leq j \leq n-1$, 
%{\allowdisplaybreaks
\begin{IEEEeqnarray}{rCl} 
 x.v_{i,j} & = & (\phi_2 \circ \phi_1) (x.(\phi_1^{-1} \circ \phi_2^{-1})(v_{i,j}))  \label{equation5.52}  \nonumber \\
& = &  \sum_{k \geq 0} (\xi^k)^{\vee} \otimes (x.(\phi_1^{-1} \circ \phi_2^{-1})(v_{i,j}))(\xi^k)    \label{equation5.53} \nonumber \\
& = & \sum_{k \geq 0} (\xi^k)^{\vee} \otimes (\phi_1^{-1} \circ \phi_2^{-1})(v_{i,j})(\xi^k x)    \label{equation5.54} \nonumber \\ 
& = & \sum_{k \geq 0} (\xi^k)^{\vee} \otimes (\phi_1^{-1} \circ \phi_2^{-1})(v_{i,j})\Big(x \xi^{k} + k \xi^{k-1} + \Big[\sum^{n-1}_{l=0} (c_{l+k} - c_l) \varepsilon_l \Big] \xi^{k-1}\Big) \label{equation5.55}   \\
& = & \sum_{k \geq 0} (\xi^k)^{\vee} \otimes \Big[x.\phi_2^{-1}(v_{i,j})(\xi^{k}) + k \phi_2^{-1}(v_{i,j})(\xi^{k-1}) + \Big[\sum^{n-1}_{l=0} (c_{l+k} - c_l) \varepsilon_l \Big]. \phi_2^{-1}(v_{i,j})(\xi^{k-1})\Big] \label{equation5.56} \nonumber \\
& = & \sum_{k \geq 0} (\xi^k)^{\vee} \otimes \Big[(\xi^i)^{\vee}(\xi^{k}) (\bar{x}^{j+1} \otimes 1) + (\xi^i)^{\vee}(\xi^{k-1}) \Big[k + \sum^{n-1}_{l=0} (c_{l+k} - c_l) \varepsilon_l\Big].(\bar{x}^j \otimes 1)\Big]   \label{equation5.57} \nonumber \\
& = &  (\xi^i)^{\vee} \otimes \bar{x}^{j+1} \otimes 1 + (\xi^{i+1})^{\vee} \otimes \Big(\Big[i + 1 + \sum^{n-1}_{l=0} (c_{l+i+1} - c_l) \varepsilon_l\Big].(\bar{x}^j \otimes 1)\Big)  \label{equation5.58} \nonumber \\
& = &  v_{i, j+1} + (i + 1 + c_{i+j+2} - c_{j+1}) v_{i+1, j}   \label{equation5.59},  
\end{IEEEeqnarray}%} 
where equation \eqref{equation5.55} follows from Proposition  \ref{theorem5.15}(ii), and equation \eqref{equation5.59} follows from the fact that $\text{Span}_{\mathbb{C}} \{\bar{x}^j \otimes 1 \}$ is isomorphic to $E_{j+1}$ as a $\mathbb{C}(\mathbb{Z}/n)$-module, and hence is fixed by $\varepsilon_{j+1}$ and killed by the other $\varepsilon_l$ (see Subsection \ref{subsection00.3}). This completes the proof. 
\end{proof}

We see from Proposition  \ref{theorem5.6} that the action of $x$ on $M$ depends on the values of the parameters $\boldsymbol{c}$. We now go on to determine for which values of $\boldsymbol{c}$ one of the elements $\psi^M_{\boldsymbol{c}, \boldsymbol{t}}$ generates $M$ as an $H_{\boldsymbol{c}} (\mathbb{Z}/n, \mathbb{C})$-module. 

\subsection{A criterion for when $\psi^M_{\boldsymbol{c}, \boldsymbol{t}}$ generates $M$ } \label{subsection5.3} 

In this section we describe how whether $\psi^M_{\boldsymbol{c}, \boldsymbol{t}}$ generates $M$ as an $H_{\boldsymbol{c}} (\mathbb{Z}/n, \mathbb{C})$-module is equivalent to whether a family of matrices, whose entries are expressions in the coordinates of $\boldsymbol{c}$ and $\boldsymbol{t}$, are all non-singular. We begin by making the following definition:

\begin{definition5.3}
We let $M_{\boldsymbol{c}, \boldsymbol{t}}(\mathbb{Z}/n, \mathbb{C}) = H_{\boldsymbol{c}} (\mathbb{Z}/n, \mathbb{C}).\psi^M_{\boldsymbol{c}, \boldsymbol{t}}$ be the $H_{\boldsymbol{c}} (\mathbb{Z}/n, \mathbb{C})$-submodule of $M$ generated by $\psi^M_{\boldsymbol{c}, \boldsymbol{t}}$.

\end{definition5.3}

We will generally just write $M_{\boldsymbol{c}, \boldsymbol{t}}$ for $M_{\boldsymbol{c}, \boldsymbol{t}}(\mathbb{Z}/n, \mathbb{C})$. 
\begin{theorem5.8} \label{theorem5.8}
We have that $M_{\boldsymbol{c}, \boldsymbol{t}}= \textnormal{Span}_{\mathbb{C}} \{x^i \xi^j.\psi^M_{\boldsymbol{c}, \boldsymbol{t}}: i \geq 0, 0 \leq j \leq n-1 \}$ and $(M_{\boldsymbol{c}, \boldsymbol{t}})_k = \textnormal{Span}_{\mathbb{C}} \{x^i \xi^j. \psi^M_{\boldsymbol{c}, \boldsymbol{t}}: i \geq 0, 0 \leq j \leq n-1, i-j = k \}$ for all $k \geq 1-n$. Furthermore, 
 \begin{equation*}  \label{equation5.3}
	(M_{\boldsymbol{c}, \boldsymbol{t}})_{k+1} = 
	\begin{cases}
	x.(M_{\boldsymbol{c}, \boldsymbol{t}})_k + \mathbb{C}.\xi^{-(k+1)}.\psi^M_{\boldsymbol{c}, \boldsymbol{t}}& \textnormal{if }1-n \leq  k < 0, \\
	x.(M_{\boldsymbol{c}, \boldsymbol{t}})_k & \textnormal{if } k \geq 0. 
	\end{cases}  
\end{equation*}
\end{theorem5.8}

\begin{proof}
By definition, $M_{\boldsymbol{c}, \boldsymbol{t}}=  \{a.\psi^M_{\boldsymbol{c}, \boldsymbol{t}}: a \in H_{\boldsymbol{c}} (\mathbb{Z}/n, \mathbb{C}) \}$. Since $H_{\boldsymbol{c}} (\mathbb{Z}/n, \mathbb{C}) \cong \mathbb{C}[x] \otimes_{\mathbb{C}} \mathbb{C} (\mathbb{Z}/n) \otimes_{\mathbb{C}} \mathbb{C}[\xi] $ as a vector space by the PBW Theorem, it follows that $M_{\boldsymbol{c}, \boldsymbol{t}}$ is spanned by vectors of the form $(x^i w \xi^{j}). \psi^M_{\boldsymbol{c}, \boldsymbol{t}}$, where $w \in \mathbb{Z}/n$. Now, 
\begin{equation*}  \label{equation5.03}
	 (x^i w \xi^{j}).\psi^M_{\boldsymbol{c}, \boldsymbol{t}} = (x^i (w. \xi)^{j}w). \psi^M_{\boldsymbol{c}, \boldsymbol{t}}  = q^j (x^i  \xi^{j}). (w.\psi^M_{\boldsymbol{c}, \boldsymbol{t}}) = q^j (x^i  \xi^{j}).\psi^M_{\boldsymbol{c}, \boldsymbol{t}} 
\end{equation*}
where the last equality follows from the fact that $\psi^M_{\boldsymbol{c}, \boldsymbol{t}}\in M^W$, that is, it is fixed by all elements of $\mathbb{Z}/n$. Hence $M_{\boldsymbol{c}, \boldsymbol{t}}= \textnormal{Span}_{\mathbb{C}} \{x^i \xi^j.\psi^M_{\boldsymbol{c}, \boldsymbol{t}}: i \geq 0, 0 \leq j \leq n-1 \}$, and as $x$ has degree $1$, $\xi$ has degree $-1$, $\psi^M_{\boldsymbol{c}, \boldsymbol{t}}$ has degree $0$ and $M$ is concentrated in degrees $\geq 1-n$, the rest of the claim then follows. 
\end{proof}

This result lets us formulate a criterion for when $\psi^M_{\boldsymbol{c}, \boldsymbol{t}}$ generates $M$ as an $H_{\boldsymbol{c}} (\mathbb{Z}/n, \mathbb{C})$-module, in terms of the action of $x$ on the graded components of $M$. 

\begin{theorem5.10}  \label{theorem5.10}

\begin{itemize}

\item[(i)]

If $t_{n-1} \neq 0$, the element $\psi^M_{\boldsymbol{c}, \boldsymbol{t}}$ generates $M$ as an $H_{\boldsymbol{c}} (\mathbb{Z}/n, \mathbb{C})$-module if and only if 
\begin{equation}  \label{equation5.5}
 	M_{k+1} = 
	\begin{cases}
	x.M_k + \mathbb{C}.\xi^{-(k+1)}.\psi^M_{\boldsymbol{c}, \boldsymbol{t}} & \textnormal{if }1-n \leq  k < 0, \\
	x.M_k & \textnormal{if } k \geq 0. 
	\end{cases}  
\end{equation}
\item[(ii)]

If $t_{n-1} = 0$, then $\psi^M_{\boldsymbol{c}, \boldsymbol{t}}$ does not generates $M$ as an $H_{\boldsymbol{c}} (\mathbb{Z}/n, \mathbb{C})$-module. 

\end{itemize}

\end{theorem5.10} 

\begin{proof}
Assume first that $t_{n-1} \neq 0$ and that \eqref{equation5.5} holds. We will prove that $(M_{\boldsymbol{c}, \boldsymbol{t}})_k = M_k$ for all $k \geq 1-n$ by induction. For the base step $k = 1-n$, observe that $\xi^{n-1}.\psi^M_{\boldsymbol{c}, \boldsymbol{t}} = \frac{t_{n-1}}{(n-1)!}. v_{0,0}$. Hence it follows from  
Propositions \ref{theorem5.9} and \ref{theorem5.8} that $(M_{\boldsymbol{c}, \boldsymbol{t}})_{1-n} = \textnormal{Span}_{\mathbb{C}} \{v_{0,0}\} = M_{1-n}$.  Next assume that $1-n \leq k <0$ and that $(M_{\boldsymbol{c}, \boldsymbol{t}})_k = M_k$. It then follows from Proposition  \ref{theorem5.8}, and our assumption that \eqref{equation5.5} holds, that
\begin{equation*} \label{equation5.66}
	 (M_{\boldsymbol{c}, \boldsymbol{t}})_{k+1} =  x.(M_{\boldsymbol{c}, \boldsymbol{t}})_k + \mathbb{C}.\xi^{-(k+1)}. \psi^M_{\boldsymbol{c}, \boldsymbol{t}} = x.M_k + \mathbb{C}.\xi^{-(k+1)}. \psi^M_{\boldsymbol{c}, \boldsymbol{t}} = M_{k+1} 
\end{equation*}
Finally assume that $(M_{\boldsymbol{c}, \boldsymbol{t}})_k = M_k$ for some $k \geq 0$. Then it again follows from Proposition  \ref{theorem5.8}(ii),  and our assumption that \eqref{equation5.5} holds,  that  
\begin{equation*} \label{equation5.661}
	(M_{\boldsymbol{c}, \boldsymbol{t}})_{k+1} =  x.(M_{\boldsymbol{c}, \boldsymbol{t}})_k  = x.M_k  = M_{k+1} 
\end{equation*}
Hence it follows by induction that $(M_{\boldsymbol{c}, \boldsymbol{t}})_k = M_k$ for all $k \geq 1-n$, so that $M_{\boldsymbol{c}, \boldsymbol{t}} = M$, and $\psi^M_{\boldsymbol{c}, \boldsymbol{t}}$ generates $M$ as an $H_{\boldsymbol{c}} (\mathbb{Z}/n, \mathbb{C})$-module. 
For the converse, assume that \eqref{equation5.5} does not hold. Then either $x.M_k + \mathbb{C}.\xi^{-(k+1)}.\psi^M_{\boldsymbol{c}, \boldsymbol{t}} \subsetneq M_{k+1} $
for some $1-n \leq k < 0$, or $x.M_k \subsetneq M_{k+1}$
for some $k \geq 0$. In the first case we have that, according to Proposition  \ref{theorem5.8},
\begin{equation*} \label{equation5.662}
	(M_{\boldsymbol{c}, \boldsymbol{t}})_{k+1}  =  x.(M_{\boldsymbol{c}, \boldsymbol{t}})_k + \mathbb{C}.\xi^{-(k+1)}. \psi^M_{\boldsymbol{c}, \boldsymbol{t}}   \subseteq x.M_k + \mathbb{C}.\xi^{-(k+1)}. \psi^M_{\boldsymbol{c}, \boldsymbol{t}}  \subsetneq M_{k+1} 
\end{equation*}
and, similarly, in the second case we have that 
\begin{equation*} \label{equation5.663}
	(M_{\boldsymbol{c}, \boldsymbol{t}})_{k+1}  =  x.(M_{\boldsymbol{c}, \boldsymbol{t}})_k  \subseteq x.M_k  \subsetneq M_{k+1} 
\end{equation*}
In either case $M_{\boldsymbol{c}, \boldsymbol{t}} \neq M$, so that $\psi^M_{\boldsymbol{c}, \boldsymbol{t}}$ does not generate $M$ as an $H_{\boldsymbol{c}} (\mathbb{Z}/n, \mathbb{C})$-module. This completes the proof of (i). Next assume that $t_{n-1} = 0$. Then, according to Propositions \ref{theorem5.9} and \ref{theorem5.8}, we have that   
\begin{equation*} \label{equation5.664}
	 (M_{\boldsymbol{c}, \boldsymbol{t}})_{1 -n} =  \textnormal{Span}_{\mathbb{C}} \{ \xi^{n-1}. \psi^M_{\boldsymbol{c}, \boldsymbol{t}}  \} = \textnormal{Span}_{\mathbb{C}} \{ \frac{t_{n-1}}{(n-1)!}.v_{0,0} \}    = 0 \neq  M_{1-n} 
\end{equation*}
Therefore $M_{\boldsymbol{c}, \boldsymbol{t}} \neq M$, which means that $\psi^M_{\boldsymbol{c}, \boldsymbol{t}}$ does not generate $M$ as an $H_{\boldsymbol{c}} (\mathbb{Z}/n, \mathbb{C})$-module. This completes the proof. \qedhere
\end{proof}

With the help of Proposition  \ref{theorem5.10}, we can formulate a criterion for when $\psi^M_{\boldsymbol{c}, \boldsymbol{t}}$ generates $M$ as an $H_{\boldsymbol{c}} (\mathbb{Z}/n, \mathbb{C})$-module in terms of the non-singularity of a set of  matrices, whose 
entries are expressions in the coordinates of $\boldsymbol{c}$ and  $\boldsymbol{t}$. In order to simplify the notation of the following result, we make the following definition:

\begin{definition5.5} \label{definition5.5}
For $i,j \in \mathbb{Z}$ we define $\Delta_{i,j} = c_i - c_j.$
\end{definition5.5}

We then have the following:

\begin{theorem5.11} \label{theorem5.11}

Let $t_{n-1} \neq 0$. 

\begin{itemize} 
\item[(i)] For $1-n \leq k \leq  0$, $\{ \xi^{-k}.\psi^M_{\boldsymbol{c}, \boldsymbol{t}},  v_{k + n -2,1}, \ldots, v_{1, k + n -2}, v_{0, k + n -1} \}$ is a basis for $M_k$. 
%Observe that this is the same as the basis for $M_k$ given in equation \eqref{equation5.6} with the element $v_{0, k +n -1}$ replaced by the element $\xi^{-k} y$. 

\item[(ii)] For $1-n \leq k < 0$, the map from $M_k$ to $M_{k+1}/(\mathbb{C}.\xi^{-(k+1)}. \psi^M_{\boldsymbol{c}, \boldsymbol{t}})$ given by first multiplying by $x$ and then projecting onto   
$M_{k+1}/(\mathbb{C}.\xi^{-(k+1)}. \psi^M_{\boldsymbol{c}, \boldsymbol{t}})$, 
\begin{equation*} \label{equation5.82}
	M_k \stackrel{x.}{\longrightarrow} M_{k+1} \stackrel{\pi}{\longrightarrow} M_{k+1}/(\mathbb{C}.\xi^{-(k+1)}.\psi^M_{\boldsymbol{c}, \boldsymbol{t}}),
\end{equation*}
has $(n + k) \times (n +k)$-matrix $F_k^{\boldsymbol{c}, \boldsymbol{t}}$ given by 
\begin{footnotesize}
\begin{equation*} \label{equation5.94}
\begin{bmatrix}
		1 - ( n + k +  \Delta_{k+1, 1})t_{n-2}t^{-1}_{n-1}  & n + k -1 + \Delta_{k+1, 2} & 0 & \ldots  & 0 \\
		- ( n + k +  \Delta_{k+1, 1})t_{n-3}t^{-1}_{n-1} & 1  & n + k -2  + \Delta_{k+1, 3}  & \ldots  & 0 \\
		- ( n + k +  \Delta_{k+1, 1})t_{n-4}t^{-1}_{n-1} & 0 & 1 & \ldots & 0  \\
		\vdots & \vdots & \vdots & \ddots & \vdots \\
		- ( n + k +  \Delta_{k+1, 1})t_{-k}t^{-1}_{n-1} & 0 & 0 & \ddots & 1 + \Delta_{k+1, n+k}  \\
		- ( n + k +  \Delta_{k+1, 1})t_{-(k+1)}t^{-1}_{n-1} & 0 & 0 & \ldots & 1 
	\end{bmatrix}
\end{equation*} \end{footnotesize}		
\noindent with respect to the basis $\{ v_{k + n -1,0}, v_{k + n -2, 1}, \ldots, v_{0, k + n -1} \}$ of  $M_k$ 
and the basis 
$\{ \bar{v}_{k + n -1 ,1}, \bar{v}_{k + n -2, 2}, \ldots, \bar{v}_{0, k + n } \}$ of $M_{k+1}/(\mathbb{C}.\xi^{-(k+1)}. \psi^M_{\boldsymbol{c}, \boldsymbol{t}})$, 
where we let $\bar{v}_{i,j} = v_{i,j} + \mathbb{C}.\xi^{-(k+1)}. \psi^M_{\boldsymbol{c}, \boldsymbol{t}} $. 
\item[(iii)] For $k \geq 0$, the map from $M_k$ to $M_{k+1}$ given by multiplication by $x$, 
\begin{equation*} \label{equation5.84}
	M_k \stackrel{x.}{\longrightarrow} M_{k+1}, 
\end{equation*}
has $n \times n$-matrix  $D_k^{\boldsymbol{c}, \boldsymbol{t}}$ given by 
\begin{footnotesize}
\begin{equation*}  \label{equation5.93}
\begin{bmatrix}
		n+k+ \Delta_{k+1, 1} & 0 & \ldots & 0 & 0 \\
		1 & n+k-1 + \Delta_{k+1, 2}  & \ldots & 0 & 0  \\
		0 & 1 & \ldots & 0 & 0 \\
		\vdots & \vdots & \ddots & \vdots & \vdots \\
		0 & 0 & \ldots & k+2 + \Delta_{k+1, n-1} & 0 \\
		0 & 0 & \ldots & 1 & k+1 + \Delta_{k+1, n} 
	\end{bmatrix}
\end{equation*}
\end{footnotesize}
\noindent with respect to the basis $\{ v_{k + n -1,0}, v_{k + n -2, 1}, \ldots, v_{k, n -1} \}$ of $M_k$ and to the basis 
$\{ v_{k + n ,0}, v_{k + n -1, 1}, \ldots, v_{k+1, n-1} \}$ of $M_{k+1}$. 
\end{itemize}
\end{theorem5.11}

\begin{proof}
Assume that $1-n \leq k \leq 0$. We know from Proposition  \ref{theorem5.9} that $M_k$ has a basis  $\{ v_{n + k -1,0}, v_{n + k -2, 1}, \ldots, v_{1, n + k -2}, v_{0, n + k -1} \}$. Since $t_{n-1} \neq 0$,  the coefficient of $v_{ n + k -1, 0}$ in 
\begin{equation} \label{equation5.87}
	\xi^{-k}.\psi^M_{\boldsymbol{c}, \boldsymbol{t}}   =  \xi^{-k}.(\sum^{n-1}_{i = 0} \frac{t_i}{(n-1)!} v_{i, n -1 -i}) = \sum^{n-1}_{i = -k} \frac{t_i}{(n-1)!} v_{i + k, n -1 -i}
\end{equation}
is non-zero, and therefore $\{ \xi^{-k}.\psi^M_{\boldsymbol{c}, \boldsymbol{t}},  v_{n + k -2,1}, \ldots, v_{1, n + k -2}, v_{0, n +  k -1} \}$
is a basis for $M_k$ as well, which proves (i). For part (ii), assume that $1-n \leq k < 0$. We see from Proposition  \ref{theorem5.6} that for $0 \leq i \leq k + n -1$ we have that 
\begin{equation*} \label{equation5.206}
	x.v_{n + k - 1 - i, i} = v_{n + k - 1 - i,i + 1} + (n + k - i + c_{k+1} - c_{i+1})v_{n + k - i, i},  
\end{equation*} 
so that 
\begin{equation*} \label{equation5.207}
	\pi (x.v_{n + k - 1 - i, i}) = \bar{v}_{n + k - 1 - i,i + 1} + (n + k - i + c_{k+1} - c_{i+1}) \bar{v}_{n + k - i, i}.
\end{equation*}  
For $1 \leq i \leq k + n -1$, this gives the last $n + k -1$ columns of $F_k^{\boldsymbol{c}, \boldsymbol{t}}$. When $i = 0$, we observe that equation \eqref{equation5.87}, with $k+1$ in place of $k$,  implies that 
\begin{equation*} \label{equation5.208}
	v_{n+k, 0} = \frac{(n-1)!}{t_{n-1} } \xi^{-(k+1)}. \psi^M_{\boldsymbol{c}, \boldsymbol{t}} - \sum^{n-2}_{i = -(k+1)} t_{n-1}^{-1} t_i v_{i + k +1, n -1 -i}. 
\end{equation*}  
Therefore, 
\begin{equation*} \label{equation5.209}
	\bar{v}_{n+k, 0} =  - \sum^{n-2}_{i = -(k+1)} t_{n-1}^{-1} t_i \bar{v}_{i + k +1, n -1 -i}, 
\end{equation*}  
so that 
\begin{IEEEeqnarray}{rCl} 
	\pi (x.v_{n + k - 1 , 0})  & = & \bar{v}_{n + k - 1 ,1} + (n + k  + c_{k+1} - c_{1}) \bar{v}_{n + k , 0} \label{equation5.86} \nonumber \\
& = & \bar{v}_{n + k - 1 ,1} - (n + k  + c_{k+1} - c_{1})  \sum^{n-2}_{i = -(k+1)} t_{n-1}^{-1} t_i \bar{v}_{i + k +1, n -1 -i}  \label{equation5.1187} \nonumber 	
\end{IEEEeqnarray} 
which gives the first column of $F_k^{\boldsymbol{c}, \boldsymbol{t}}$ and completes the proof of (ii). 
Finally, part (iii) is a direct consequence of Proposition  \ref{theorem5.6}.  
\end{proof}

Combining Propositions \ref{theorem5.10} and \ref{theorem5.11} we obtain: 

\begin{theorem5.12} \label{theorem5.12}
The element $\psi^M_{\boldsymbol{c}, \boldsymbol{t}}$ generates $M$ as an $H_{\boldsymbol{c}} (\mathbb{Z}/n, \mathbb{C})$-module, if and only if the matrices $D_k^{\boldsymbol{c}, \boldsymbol{t}}$ are non-singular for every $k \geq 0$, and the matrices 
$F_k^{\boldsymbol{c}, \boldsymbol{t}}$ are non-singular for every $1-n \leq k <0$. 
\end{theorem5.12}

We go on to analyse for which values of the parameters  $\boldsymbol{c}$ and $\boldsymbol{t}$ the matrices $D_k^{\boldsymbol{c}, \boldsymbol{t}}$ and $F_k^{\boldsymbol{c}, \boldsymbol{t}}$ are all non-singular.

\subsection{A set of parameter values $\boldsymbol{c}$ for which $P_{\textnormal{KZ}} \cong \Delta_{\boldsymbol{c}}(S(\mathfrak{h})_W)$ } \label{subsection5.4} 

To determine for which parameter values $\boldsymbol{c}$ and $\boldsymbol{t}$ the matrices $D_k^{\boldsymbol{c}, \boldsymbol{t}}$ and $F_k^{\boldsymbol{c}, \boldsymbol{t}}$ are all non-singular, we begin by noting the following. 

\begin{theorem5.18}  \label{theorem5.18}
For $\boldsymbol{t} = (0, 0, \dots, 0, 1)$, that is $t_i = 0$ for $0 \leq i \leq n-2$ and $t_{n-1} = 1$, the matrices $F_k^{\boldsymbol{c}, \boldsymbol{t}}$ are non-singular for all parameter values $\boldsymbol{c}$ and all $1-n \leq k < 0$. 

\end{theorem5.18}

\begin{proof}
Let $\boldsymbol{c} \in \mathbb{C}^n$ with $c_n = 0$ (cf. end of Subsection \ref{subsection00.3}) and let $1-n \leq k < 0$.  We see from Proposition  \ref{theorem5.11} (ii) that when $\boldsymbol{t} = (0, 0, \dots, 0, 1)$, the matrix $F_k^{\boldsymbol{c}, \boldsymbol{t}}$ is upper-triangular with all its diagonal entries equal to 1, which completes the proof. 
\end{proof}

We go on to study the matrices $D_k^{\boldsymbol{c}, \boldsymbol{t}}$. As they are lower triangular, they are invertible exactly when all the diagonal entries are non-zero. This allows us to find a simple criterion for determining 
when all the matrices $D_k^{\boldsymbol{c}, \boldsymbol{t}}$ are non-singular. 

\begin{theorem5.1} \label{theorem5.1}
Let $\boldsymbol{c}, \boldsymbol{t} \in \mathbb{C}^n$ with $c_n = 0$.  The matrices $D_k^{\boldsymbol{c}, \boldsymbol{t}}$ are non-singular for every $k \geq 0$, if and only if for all $1 \leq i, j \leq n$ we have that whenever $c_i - c_j \equiv j - i \pmod n$, then $c_i - c_j = j - i$. 

\end{theorem5.1}

\begin{proof}
%We prove the contrapositive statement, namely that there exists at least one $k \geq 0$ such that $D_k$ is singular, if and only if there exist integers $1 \leq i,j \leq n$ such that $c_i - c_j \equiv j - i \pmod n$ but $c_i - c_j \neq j-i$. 
First, assume that $1 \leq i,j \leq n$ are integers such that $c_i - c_j \equiv j - i \pmod n$ but $c_i - c_j \neq j-i$. Possibly interchanging $i$ and $j$, we can assume that $c_i- c_j <0$, so that $c_i - c_j \leq (j - i) -n$ since $c_i - c_j \neq j-i$. If we let $k = j - (n+1) - (c_i- c_j)$,  then by the above remarks
\begin{equation*} \label{equation5.100}
	k \geq   j - (n+1) - [(j-i) - n] =  i-1 \geq 0. 
\end{equation*}  
We show that this implies that $D_k^{\boldsymbol{c}, \boldsymbol{t}}$ is singular. Since 
\begin{equation*} \label{equation5.101}
	k + 1 =  j - n -  (c_i- c_j) \equiv  j -  (j-i)  = i \pmod n,
\end{equation*}  
the $j$th diagonal entry of $D_k^{\boldsymbol{c}, \boldsymbol{t}}$ is 
\begin{equation*} \label{equation5.104}
	(n + k + 1) -j + \Delta_{k+1,j}  =  (j -  (c_i- c_j)) - j + \Delta_{i,j}  = 0.
\end{equation*}  
Therefore, as $D_k^{\boldsymbol{c}, \boldsymbol{t}}$ is lower triangular, it follows that $D_k^{\boldsymbol{c}, \boldsymbol{t}}$ is singular.
Conversely, assume there exists an integer $k \geq 0$ such that $D_k^{\boldsymbol{c}, \boldsymbol{t}}$ is singular. Then, one of its diagonal entries must be zero, so that there exists an integer $1 \leq j \leq n$ such that $(n+k+1) - j + \Delta_{k+1, j} = 0$. 
Let $1 \leq i \leq n$ be such that $i \equiv k+ 1 \pmod n $. Then 
\begin{equation*} \label{equation5.107}
	c_i - c_j  =  \Delta_{k+1, j}  =  j - (n +k+1) \equiv  j - i \pmod n
\end{equation*}  
but 
\begin{equation*} \label{equation5.108}
	c_i - c_j =  j - (n +k+1) \neq  j - i 
\end{equation*}  
as $n +k+1 \geq n+1 >i$. This completes the proof. 
\end{proof}

In light of this claim, we make the following definition. 

\begin{definition5.7}
We define $\mathcal{F}$ to be the subset of $\mathbb{C}^n$ such that $\boldsymbol{c} = (c_1, c_2, \dots, c_n) \in \mathcal{F}$ if $c_n = 0$ and for all $1 \leq i, j \leq n$ we have that whenever $c_i - c_j \equiv j - i \pmod n$, then $c_i - c_j = j - i$
\end{definition5.7}

Combining Corollaries  \ref{theorem00.5} and \ref{theorem00.8}, Proposition  \ref{theorem3.14}, Corollary \ref{theorem5.12} and Propositions \ref{theorem5.18} and \ref{theorem5.1}, we thus obtain the following: 

\begin{theorem5.19}  \label{theorem5.19}
We have that  $P_{\textnormal{KZ}} \cong \Delta_{\boldsymbol{c}}(S(\mathfrak{h})_W)$
in category $\mathcal{O}_{\boldsymbol{c}} (\mathbb{Z}/n, \mathbb{C})$ when $\boldsymbol{c} \in \mathcal{F}$.

\end{theorem5.19}

\begin{proof}
Suppose that $\boldsymbol{c} \in \mathcal{F}$ and let $\boldsymbol{t} = (0, 0, \dots, 0, 1)$. Then the matrices $D_k^{\boldsymbol{c}, \boldsymbol{t}}$ and $F_k^{\boldsymbol{c}, \boldsymbol{t}}$ are all non-singular by Propositions \ref{theorem5.18} and \ref{theorem5.1}. Therefore, the element $\psi^M_{\boldsymbol{c}, \boldsymbol{t}} \in M = S(\mathfrak{h})^{\circledast} \otimes_{\mathbb{C}} S(\mathfrak{h}^*)_W^{\otimes}$ generates $M$ as an $H_{\boldsymbol{c}} (\mathbb{Z}/n, \mathbb{C})$-module by Corollary \ref{theorem5.12}. Hence, the element $\psi_{\boldsymbol{t}} \in \text{Hom}_{\mathbb{C} W}(S(\mathfrak{h}^*)_W, \textnormal{Har}^*)$ corresponds to an isomorphism in the space $\mathcal{S} = \text{Hom}_{H_c(W, \mathfrak{h})} (\Delta_{\boldsymbol{c}} (S(\mathfrak{h})_W), \nabla_{\boldsymbol{c}}(S(\mathfrak{h}^*)_W^{\otimes}) )$ by Proposition  \ref{theorem3.14}. Therefore, $\Delta_{\boldsymbol{c}} (S(\mathfrak{h})_W)$ is a tilting module, and it follows from Corollaries \ref{theorem00.5} and \ref{theorem00.8} that $P_{\textnormal{KZ}} \cong \Delta_{\boldsymbol{c}}(S(\mathfrak{h})_W)$
in category $\mathcal{O}_{\boldsymbol{c}} (\mathbb{Z}/n, \mathbb{C})$, which completes the proof. 
\end{proof}

Next, we observe that the set $\mathcal{F}$ of parameter values is good in the sense of Definition \ref{definition00.1} (recall that $c_i = nk_i$). 

\begin{theorem5.2}  \label{theorem5.2}
For every $\boldsymbol{c} = (c_1, c_2, \dots, c_n) \in \mathbb{C}^n$ with $c_n = 0$, there exists $\boldsymbol{c}'  = (c'_1, c'_2, \dots, c'_n)\in \mathbb{C}^n$ with $c'_n = 0$ such that $c_i \equiv c'_i \,\pmod n$ for all $1 \leq i \leq n$ and such that $\boldsymbol{c}' \in \mathcal{F}$.

\end{theorem5.2}

\begin{proof}
Let $\boldsymbol{c} = (c_1, c_2, \dots, c_n) \in \mathbb{C}^n$ with $c_n = 0$. We introduce an equivalence relation on $\{1, \ldots, n \}$ by declaring that $i \sim j$ if 
$c_i - c_j \equiv j - i \pmod n$. Let $[i_1], \ldots, [i_l]$ be the equivalence classes of $\{1, \ldots, n \}$ with respect to $\sim$. Without loss of generality, we can assume that $i_l = n$. For $j \in [i_k]$, we define $c_j' := c_{i_k} + i_k - j$. Then
\begin{equation*} \label{equation5.999} 
	c_j - c_j'  =  c_j - (c_{i_k} + i_k - j)   =    (c_j - c_{i_k}) - (i_k - j) \equiv     0 \pmod n,
\end{equation*} 
where the last equality follows from the fact that $j \sim i_k$. Also $c_n' = 0$ as $i_l = n$. 
Next, observe that  if $i \in [i_u]$ and $j \in [i_v]$, then 
\begin{equation*} \label{equation5.993} 
	c'_i - c'_j = (c_{i_u} + i_u - i) - (c_{i_v} + i_v - j)   =    [(c_{i_u} - c_{i_v}) - (i_v - i_u)] +(j - i).
\end{equation*} 
Therefore, if $c'_i - c'_j \equiv j - i \pmod n$, then $c_{i_u} - c_{i_v} \equiv i_v - i_u \pmod n$ so that $i_u \sim i_v$ and hence $u = v$. Thus
\begin{equation*} \label{equation5.114} 
	c_i' - c_j'  =  (c_{i_u} + i_u - i) - (c_{i_u} + i_u - j)   = j-i,
\end{equation*} 
so we can pick $\boldsymbol{c}'  = (c'_1, c'_2, \dots, c'_n)$, and this completes the proof. 
\end{proof}

Taken together, Propositions \ref{theorem5.19} and \ref{theorem5.2} give the main result of this section: 

\begin{theorem000.2} \label{theorem000.2}
There is a good set of parameter values $\mathcal{F}$ such that $P_{\textnormal{KZ}} \cong \Delta_{\boldsymbol{c}}(S(\mathfrak{h})_W)$ 
in category $\mathcal{O}_{\boldsymbol{c}} (\mathbb{Z}/n, \mathbb{C})$ when $\boldsymbol{c} \in \mathcal{F}$, where $\mathcal{F}$ consists of those parameters $\boldsymbol{c}$ which satisfy that $c_i - c_j = j - i$ whenever $c_i - c_j \equiv j - i \,\pmod n$ for $1 \leq i,j \leq n$. 
\end{theorem000.2}

We will see in the next section that $P_{\textnormal{KZ}} \cong \Delta_{\boldsymbol{c}}(S(\mathfrak{h})_W)$ in category $\mathcal{O}_{\boldsymbol{c}} (\mathbb{Z}/n, \mathbb{C})$ if and only if $\boldsymbol{c} \in \mathcal{F}$.

\section{The endomorphism ring of $\Delta_{\boldsymbol{c}}(S(\mathfrak{h})_W)$} \label{chapter 6}

We saw in the last section that there exists a good set $\mathcal{F}$  of parameter values $\boldsymbol{c}$ for which $P_{\textnormal{KZ}} \cong \Delta_{\boldsymbol{c}}(S(\mathfrak{h})_W)$ in category $\mathcal{O}_{\boldsymbol{c}} (\mathbb{Z}/n, \mathbb{C})$. 
In this section, we will show that this is best possible, in the sense that  $P_{\textnormal{KZ}}$ is not isomorphic to $\Delta_{\boldsymbol{c}}(S(\mathfrak{h})_W)$ for parameter values $\boldsymbol{c}$ outside of $\mathcal{F}$. We will demonstrate this by examining the endomorphism rings of $P_{\textnormal{KZ}}$ and $\Delta_{\boldsymbol{c}}(S(\mathfrak{h})_W)$.   

\subsection{The endomorphism ring of $P_{\textnormal{KZ}}$} \label{subsection7.1}
 We begin by noting the following result:

\begin{theorem7.1} \label{theorem7.1}
The dimension of the endomorphism ring of $P_{\textnormal{KZ}}$ satisfies 
\begin{equation*} \label{equation7.1}
	 \textnormal{dim}( \textnormal{End}_{H_{\boldsymbol{c}} (\mathbb{Z}/n, \mathbb{C})} (P_{\textnormal{KZ}})) = n
\end{equation*}
for all parameter values $\boldsymbol{c}$. 

\end{theorem7.1}

\begin{proof}
This is a consequence of the decomposition $P_{\textnormal{KZ}} = \bigoplus_{E \in \textnormal{Irr} (W)} \textnormal{dim} (\textnormal{KZ}(L(E))) P(E)$ and is proved as part of Proposition 5.15 in \cite{GGOR}. 
\end{proof}

We proceed to show that the dimension of $\textnormal{End}_{H_{\boldsymbol{c}} (\mathbb{Z}/n, \mathbb{C})} (\Delta_{\boldsymbol{c}}(S(\mathfrak{h})_W))$ is $n$ exactly when $\boldsymbol{c} \in \mathcal{F}$. 

\subsection{The correspondence between  $\textnormal{End}_{H_{\boldsymbol{c}} (\mathbb{Z}/n, \mathbb{C})} (\Delta_{\boldsymbol{c}}(S(\mathfrak{h})_W))$ and $\Delta_{\boldsymbol{c}}(S(\mathfrak{h})_W)^{\circ}$ } \label{subsection7.2}

\begin{theorem7.2} \label{theorem7.2}
There is a linear one-to-one correspondence between the elements of $\textnormal{End}_{H_{\boldsymbol{c}} (\mathbb{Z}/n, \mathbb{C})} (\Delta_{\boldsymbol{c}}(S(\mathfrak{h})_W))$ and the elements of $\Delta_{\boldsymbol{c}}(S(\mathfrak{h})_W)$ that are killed by $S(\mathfrak{h})^W_{<0}$ and fixed by $W$. In particular, $ \textnormal{dim} (\textnormal{End}_{H_{\boldsymbol{c}} (\mathbb{Z}/n, \mathbb{C})} (\Delta_{\boldsymbol{c}}(S(\mathfrak{h})_W))) = \textnormal{dim} (\Delta_{\boldsymbol{c}}(S(\mathfrak{h})_W)^{\circ})$
where $\Delta_{\boldsymbol{c}}(S(\mathfrak{h})_W)^{\circ} =  \{y \in \Delta_{\boldsymbol{c}}(S(\mathfrak{h})_W) : S(\mathfrak{h})^W_{<0}.y = 0 \,\, \textnormal{and} \,\,  w.y = y \,\, \textnormal{for all} \,\, w \in W   \} $.

\end{theorem7.2}

\begin{proof}
As $\Delta_{\boldsymbol{c}}(S(\mathfrak{h})_W)  =  H_{\boldsymbol{c}} (\mathbb{Z}/n, \mathbb{C}) \otimes_{S(\mathfrak{h}) \rtimes W} S(\mathfrak{h})_W$, it follows from induction and restriction that $\textnormal{End}_{H_{\boldsymbol{c}} (\mathbb{Z}/n, \mathbb{C})} (\Delta_{\boldsymbol{c}}(S(\mathfrak{h})_W)) \cong \textnormal{Hom}_{S(\mathfrak{h}) \rtimes W} (S(\mathfrak{h})_W, \Delta_{\boldsymbol{c}}(S(\mathfrak{h})_W))$. 
Now, it was seen in the proof of Proposition \ref{theorem3.12} that $S(\mathfrak{h})_W \cong (S(\mathfrak{h}) \rtimes W )/ \langle S(\mathfrak{h})^W_{<0}, w-1 \rangle_{w \in W}$ 
as $S(\mathfrak{h}) \rtimes W$-modules. Therefore, the map $f \mapsto f(1 + \langle S(\mathfrak{h})^W_{<0} \rangle)$
is a linear isomorphism between $\textnormal{Hom}_{S(\mathfrak{h}) \rtimes W} (S(\mathfrak{h})_W, \Delta_{\boldsymbol{c}}(S(\mathfrak{h})_W))$ and $\Delta_{\boldsymbol{c}}(S(\mathfrak{h})_W)^{\circ}$ and this completes the proof. 
\end{proof}

\subsection{The structure of $\Delta_{\boldsymbol{c}}(S(\mathfrak{h})_W)^{\circ}$ } \label{subsection7.4}
By the PBW-Theorem,  $\Delta_{\boldsymbol{c}}(S(\mathfrak{h})_W) \cong  \mathbb{C} [x] \otimes_{\mathbb{C}} \mathbb{C}[\xi] / \langle \xi^n \rangle $ as a vector space, 
so that $\Delta_{\boldsymbol{c}}(S(\mathfrak{h})_W)$ has a basis given by $ \{v_{i,j} =  x^i \otimes \bar{\xi}^j : i \geq 0, 0 \leq j \leq n-1\}$, where $\bar{\xi} = \xi + \langle \xi^n \rangle$. In order to determine for which parameter values $\boldsymbol{c}$ the dimension of $\Delta_{\boldsymbol{c}}(S(\mathfrak{h})_W)^{\circ}$ is $n$, we calculate the action of the generators $x$, $s$ and $\xi$ of $H_{\boldsymbol{c}} (\mathbb{Z}/n, \mathbb{C})$ on this basis. 
\begin{theorem7.3}  \label{theorem7.3}
The action of the generators $x$, $s$ and $\xi$ of $H_{\boldsymbol{c}} (\mathbb{Z}/n, \mathbb{C})$ on the basis $ \{v_{i,j} =  x^i \otimes \bar{\xi}^j : i \geq 0, 0 \leq j \leq n-1\}$ of $\Delta_{\boldsymbol{c}}(S(\mathfrak{h})_W)$ is given by 
\begin{equation} \label{equation7.07}
	x.v_{i,j} = v_{i+1,j}, \quad s.v_{i,j} = q^{-(i-j)} v_{i,j} \quad \textnormal{and} \quad \xi.v_{i,j} = v_{i,j+1} + (i + c_{i-j} - c_{-j})v_{i-1, j}
\end{equation}
for $i \geq 0$ and $0 \leq j \leq n-1$, where we let  $v_{i,j} = 0$ if $i < 0$ or $j \geq n$
\end{theorem7.3}

\begin{proof}
The action of $x$ and $s$ on the basis follows immediately from the definitions.
We prove the action of $\xi$ by induction on $i$.  The statement is clearly true when $i = 0$. To prove it also holds for $i > 0$,  observe that \eqref{equation00.3}, \eqref{equation00.4} and the fact that $c_r = nk_r$ implies that $[\xi, x] = 1 + \sum^{n-1}_{k = 0} (c_{k+1} - c_k) \varepsilon_k$, and recall from Subsection \ref{subsection00.3} that $\varepsilon_k$ fixes the simple $\mathbb{Z}/n$-module $E_l$ if $k \equiv l \pmod n$ and kills it otherwise. As $\mathbb{C}.v_{i,j} \cong E_{i-j}$ as a $\mathbb{Z}/n$-module, it therefore follows that for $i > 0$, 
\begin{equation*} \label{equation7.007}
	\xi.v_{i,j}  = (x\xi +  [\xi, x] ).v_{i-1, j} = (x\xi +  1 + \sum^{n-1}_{k = 0} (c_{k+1} - c_k) \varepsilon_k).v_{i-1, j}   = v_{i, j+1}  + (i +  c_{i-j}- c_{-j}).v_{i-1, j}
\end{equation*}
from which the claim follows by induction. 
\end{proof}

We now show that whether the dimension of the endomorphism ring of $\Delta_{\boldsymbol{c}}(S(\mathfrak{h})_W)$ is $n$ is determined by whether $\Delta_{\boldsymbol{c}}(S(\mathfrak{h})_W)^{\circ}$ contains any non-zero homogeneous elements having degree a positive multiple of $n$. 

\begin{theorem7.4}  \label{theorem7.4}
The dimension of $\Delta_{\boldsymbol{c}}(S(\mathfrak{h})_W)^{\circ}$ satisfies 
\begin{equation*} \label{equation7.12}
	\textnormal{dim} (\Delta_{\boldsymbol{c}}(S(\mathfrak{h})_W)^{\circ}) = n + \textnormal{dim} (\bigoplus_{k > 0} \Delta_{\boldsymbol{c}}(S(\mathfrak{h})_W)^{\circ}_{kn}). 
\end{equation*}
\end{theorem7.4}

\begin{proof}
It is clear that
\begin{IEEEeqnarray}{rCl} 
	\Delta_{\boldsymbol{c}}(S(\mathfrak{h})_W)^{\circ} & = & \{y \in \Delta_{\boldsymbol{c}}(S(\mathfrak{h})_W) : S(\mathfrak{h})^W_{<0}.y = 0 \,\, \textnormal{and} \,\,  w.y = y \,\, \textnormal{for all} \,\, w \in W   \}    \label{equation7.8}  \nonumber \\
	& = & \{y \in \Delta_{\boldsymbol{c}}(S(\mathfrak{h})_W) : \xi^n.y = 0 \,\, \textnormal{and} \,\,  s.y = y \},    \label{equation7.9} \nonumber 
\end{IEEEeqnarray} 
and that it is a homogeneous subspace of $\Delta_{\boldsymbol{c}}(S(\mathfrak{h})_W)$. Furthermore it is clear from \eqref{equation7.07} that the generator $s$ of $\mathbb{Z}/n$ acts on $\Delta_{\boldsymbol{c}}(S(\mathfrak{h})_W)_i$ by multiplication by $q^{-i}$. This implies that $\Delta_{\boldsymbol{c}}(S(\mathfrak{h})_W)^{\circ}$ is concentrated in degrees divisible by $n$. Next, $\Delta_{\boldsymbol{c}}(S(\mathfrak{h})_W)$ is concentrated in degrees greater than or equal to $1-n$, and since $\xi$ sits in degree $-1$, it follows that $ \Delta_{\boldsymbol{c}}(S(\mathfrak{h})_W)_0 \subseteq \Delta_{\boldsymbol{c}}(S(\mathfrak{h})_W)^{\circ}$.
Since $ \Delta_{\boldsymbol{c}}(S(\mathfrak{h})_W)_0 = \textnormal{Span}_{\mathbb{C}} \{v_{i,i} : 0 \leq i \leq n-1\}$, 
this completes the proof.  
\end{proof}

We go on to show that $\Delta_{\boldsymbol{c}}(S(\mathfrak{h})_W)^{\circ}$ is concentrated in degree zero precisely for parameter values $\boldsymbol{c} \in \mathcal{F}$. Recall the notation $\Delta_{i, j} = c_i - c_j$ from Definition \ref{definition5.5}.

\begin{theorem7.5}  \label{theorem7.5}

\begin{itemize}

\item[(i)]

For each $k \in \mathbb{Z}_{>0}$ and for each $0 \leq i \leq n-1$  we have that 
\begin{equation*} \label{equation7.13}
	\xi^n.v_{kn + i, i}  = \prod^{n-1}_{j=0} (kn + [((n-j) - (n - i)) - \Delta_{n-i, n-j}]) v_{(k-1)n +i, i} + L_i   
\end{equation*}
where $L_i \in \textnormal{Span}_{\mathbb{C}} \{ v_{(k-1)n +j, j} : i < j \leq n-1 \}$. 
\item[(ii)]
For $k \in \mathbb{Z}_{>0}$, the map 
\begin{equation*} \label{equation7.14}
	(\xi^n)^*_k : \Delta_{\boldsymbol{c}}(S(\mathfrak{h})_W)_{kn} \to \Delta_{\boldsymbol{c}}(S(\mathfrak{h})_W)_{(k-1)n}
\end{equation*}
given by left multiplication by $\xi^n$ has determinant 
\begin{equation*} \label{equation7.15}
	\textnormal{det} (\xi^n)^*_k = \prod_{1 \leq i,j \leq n} (kn + [(i - j) - \Delta_{j, i}]).
\end{equation*}
\end{itemize}
\end{theorem7.5}

\begin{proof}
Let $k \in \mathbb{Z}_{>0}$ and $0 \leq i \leq n-1$. Using \eqref{equation7.07}, an easy induction argument shows that for $1 \leq r \leq n$,
\begin{equation*} \label{equation7.44}
	\xi^r.v_{kn + i, i}  = \prod^{r-1}_{j=0} (kn + [((n-j) - (n - i)) - \Delta_{n-i, n-j}]) v_{(kn-r) +i, i} + L^{(r)}_i,    
\end{equation*}
where $L^{(r)}_i \in \textnormal{Span}_{\mathbb{C}} \{ v_{(kn -r) +j, j} : i < j \leq n-1 \}$, and this proves (i). 
It follows that the matrix of $(\xi^n)^*_k$ with respect to the bases $\{v_{kn + i, i} : 0 \leq i \leq n-1   \}$ and $\{v_{(k-1)n + i, i} : 0 \leq i \leq n-1   \}$ of 
$\Delta_{\boldsymbol{c}}(S(\mathfrak{h})_W)_{kn}$ and $\Delta_{\boldsymbol{c}}(S(\mathfrak{h})_W)_{(k-1)n}$ respectively is lower triangular with its $i$th diagonal entry equal to  
\begin{equation*} \label{equation7.51}
	d_i = \prod^{n-1}_{j=0} (kn + [((n-j) - (n - (i-1))) - \Delta_{n-(i-1), n-j}])
\end{equation*}
which proves (ii).
\end{proof}
Combining Propositions \ref{theorem7.4} and \ref{theorem7.5} gives the following: 
\begin{theorem7.6}  \label{theorem7.6}
The dimension of the endomorphism ring of $\Delta_{\boldsymbol{c}}(S(\mathfrak{h})_W)$ satisfies 
\begin{equation*} \label{equation7.16}
	\textnormal{dim} (\textnormal{End}_{H_{ \boldsymbol{c}} (\mathbb{Z}/n, \mathbb{C})} (\Delta_{\boldsymbol{c}}(S(\mathfrak{h})_W))) = n
\end{equation*}
if and only if $\boldsymbol{c} \in \mathcal{F}$. 
\end{theorem7.6}

\begin{proof}
From equation \eqref{equation7.07} it follows that all elements of $\Delta_{\boldsymbol{c}}(S(\mathfrak{h})_W)_{kn}$ are fixed by $\mathbb{Z}/n$ for every integer $k$. It therefore follows from Propositions \ref{theorem7.4} and \ref{theorem7.5}(ii) that we have $\textnormal{dim} (\textnormal{End}_{H_{ \boldsymbol{c}} (\mathbb{Z}/n, \mathbb{C})} (\Delta_{\boldsymbol{c}}(S(\mathfrak{h})_W))) = n$ if and only if 

\begin{equation*} \label{equation7.52}
	\textnormal{det} (\xi^n)^*_k = \prod_{1 \leq i,j \leq n} (kn + [(i - j) - (c_j - c_i)]) \neq 0
\end{equation*}
for all positive integers $k$. This in turn happens if and only if for all $1 \leq i,j \leq n$ we have that either $c_i - c_j \not\equiv j - i \pmod n$ or $c_i - c_j = j - i$, or equivalently if and only if $\boldsymbol{c}  \in \mathcal{F}$, which completes the proof. 
\end{proof}

\subsection{The set of parameter values $\boldsymbol{c}$ for which $P_{\textnormal{KZ}} \cong \Delta_{\boldsymbol{c}}(S(\mathfrak{h})_W)$ } \label{subsection7.5}

Taken together with the results from Section \ref{chapter 5}, Corollary \ref{theorem7.6} gives us the main result of this section: 

\begin{theorem000.3} \label{theorem000.3}
In category $\mathcal{O}_{\boldsymbol{c}} (\mathbb{Z}/n, \mathbb{C})$, the modules $P_{\textnormal{KZ}}$ and $\Delta_{\boldsymbol{c}}(S(\mathfrak{h})_W)$ are isomorphic if and only if  $\boldsymbol{c} \in \mathcal{F}$. 
\end{theorem000.3}

\begin{proof}

It was seen in Proposition \ref{theorem5.19} that $P_{\textnormal{KZ}} \cong \Delta_{\boldsymbol{c}}(S(\mathfrak{h})_W)$ when $\boldsymbol{c} \in \mathcal{F}$. Furthermore, we see from Propositions \ref{theorem7.1} and \ref{theorem7.6} that if $\boldsymbol{c} \notin \mathcal{F}$, then the endomorphism rings of $P_{\textnormal{KZ}}$ and  $\Delta_{\boldsymbol{c}}(S(\mathfrak{h})_W)$ have different dimensions, so that necessarily  $P_{\textnormal{KZ}} \not\cong \Delta_{\boldsymbol{c}}(S(\mathfrak{h})_W)$. This completes the proof. 
\end{proof}

In the next section, we will complete the algebraic description of the KZ-functor when $\boldsymbol{c} \in \mathcal{F}$ by determining the action of the the Hecke algebra $\mathcal{H}_{\boldsymbol{c}}(W, \mathfrak{h}))$ on 
$\Delta_{\boldsymbol{c}}(S(\mathfrak{h})_W)$.

\section{The action of $\mathcal{H}_{\boldsymbol{c}} (\mathbb{Z}/n, \mathbb{C})$ on $\Delta_{\boldsymbol{c}}(S(\mathfrak{h})_W)$}  \label{chapter 7}

In the last two sections, we saw that $P_{\textnormal{KZ}} \cong \Delta_{\boldsymbol{c}}(S(\mathfrak{h})_W)$ as $H_{\boldsymbol{c}}(\mathbb{Z}/n, \mathbb{C})$-modules precisely when $\boldsymbol{c} \in \mathcal{F}$, and that $\mathcal{F}$ is a good set of parameter values in the sense of Definition \ref{definition00.1}. Recall from Subsection \ref{subsection00.7} that $P_{\textnormal{KZ}}$ represents the KZ-functor through its structure as an $(H_{\boldsymbol{c} }(W, \mathfrak{h}), \mathcal{H}_{\boldsymbol{c} }(W, \mathfrak{h}))$-bimodule. In order to have an algebraic description of the KZ-functor for $\boldsymbol{c} \in \mathcal{F}$, it thus remains to describe the action of $ \mathcal{H}_{\boldsymbol{c}}(\mathbb{Z}/n, \mathbb{C}) $ on $\Delta_{\boldsymbol{c}}(S(\mathfrak{h})_W)$, or in other words to determine the algebra homomorphism 
\begin{equation}  \label{equation6.101}
\phi: \mathcal{H}_{\boldsymbol{c}}(\mathbb{Z}/n, \mathbb{C}) \longrightarrow \textnormal{End}_{H_{\boldsymbol{c}}(\mathbb{Z}/n, \mathbb{C})} (\Delta_{\boldsymbol{c}}(S(\mathfrak{h})_W))^{\textnormal{opp}}
\end{equation}
induced by the KZ-functor (observe that $\phi$ is in fact an isomorphism by Theorem 5.15 of \cite{GGOR}). 

\subsection{Change of rings} \label{section6.00} In order to describe the isomorphism of \eqref{equation6.101}, we will consider the rational Cherednik algebra over rings other than $\mathbb{C}$, as described in Subsection \ref{subsection00.10}. For $\boldsymbol{c} = (c_1, \ldots, c_n) \in \mathcal{F}$, let $R'$ be the ring of germs of holomorphic functions in $n-1$ complex variables $z_1, \ldots, z_{n-1}$ about the point $(c_1, \ldots, c_{n-1}) \in \mathbb{C}^{n-1}$. Then $R'$ is a so called analytic algebra (cf. Section 5.1 of \cite{RSVV}). In particular, $R'$ is a Noetherian, regular, local integral domain, and we denote its unique maximal ideal by $\mathfrak{m}'$ so that $R'/\mathfrak{m}' \cong \mathbb{C}$. We also denote by $K' = \textnormal{Frac}(R')$ the field of fractions of $R'$, and by $H^{R'}_{\boldsymbol{c} }(\mathbb{Z}/n, \mathbb{C})$ and $H^{K'}_{\boldsymbol{c} }(\mathbb{Z}/n, \mathbb{C})$ the rational Cheredink algebra associated with $\mathbb{Z}/n$, $\mathbb{C}$ and $\boldsymbol{c}$ considered over $R'$ and $K'$ respectively, with parameters chosen such that $c_i^{R'} = c_i^{K'} = z_i$ for $1 \leq i \leq n-1$ (here $c_i^{R'}$ and $c_i^{K'}$ are the constants by which the element $z$ acts on $R' \otimes_{\mathbb{C}} E_i$ and $K' \otimes_{\mathbb{C}} E_i$ respectively, cf. Subsection \ref{subsection00.3}). For the purpose of simplifying notation below, we also set $z_n := 0$ so that $c_n^{R'} = c_n^{K'} = z_n = 0$. 

%[We make this choice of parameters to have compatibility with the ideal $\mathfrak{m}'$, in the sense that $c_i^{R'}$ becomes $c_i$ in $R'/\mathfrak{m}'$. This is important for the action of $\mathbf{eu}$ on $\Delta_{\boldsymbol{c}}(S(\mathfrak{h})_W)$ below. Observe from Remark 3.2 of \cite{GGOR} that the $c_i^{R'}$ and $c_i^{K'}$ can be chosen arbitrarily, as the $c^{R'}(s^i)$ and $c^{K'}(s^i)$ can be recovered from them.] 

\subsection{The cyclotomic Hecke algebra $\mathcal{H}_{\boldsymbol{c}} (\mathbb{Z}/n, \mathbb{C})$} \label{section6.01}

Equation \eqref{equation00.9} of Subsection \ref{subsection00.6}  and the fact that $c_i = nk_i$, as was seen in Subsection \ref{subsection00.3}, imply that 
\begin{equation*}  \label{equation6.83}
	\mathcal{H}_{\boldsymbol{c}} (\mathbb{Z}/n, \mathbb{C}) = \mathbb{C}[T]/ \langle \prod^n_{j = 1} (T - q^{-j}.q_j^{-1}) \rangle,  
\end{equation*}
where $q = \textnormal{exp} (\frac{2 \pi i}{n})$ and $q_j = \textnormal{exp} (\frac{2 \pi i}{n} c_j)$ for $1 \leq j \leq n$ and that for $k = R', K'$, 
\begin{equation}  \label{equation6.01}
	\mathcal{H}^k_{\boldsymbol{c}} (\mathbb{Z}/n, \mathbb{C}) = k[T]/ \langle \prod^n_{j = 1} (T - q^{-j}.{q'}_j^{-1}) \rangle,  
\end{equation}
where again $q = \textnormal{exp} (\frac{2 \pi i}{n})$ and $q'_j = \textnormal{exp} (\frac{2 \pi i}{n} z_j)$ for $1 \leq j \leq n$. (Observe that $\mathcal{H}_{\boldsymbol{c}} (\mathbb{Z}/n, \mathbb{C}) = (R'/\mathfrak{m}') \otimes_{R'} \mathcal{H}^{R'}_{\boldsymbol{c}} (\mathbb{Z}/n, \mathbb{C})$ and $\mathcal{H}^{K'}_{\boldsymbol{c}} (\mathbb{Z}/n, \mathbb{C}) = K' \otimes_{R'} \mathcal{H}^{R'}_{\boldsymbol{c}} (\mathbb{Z}/n, \mathbb{C})$). 

\subsection{The maps $\phi$, $\phi_{R'}$ and $\phi_{K'}$} \label{section6.05} 
Let $\mathfrak{m}$ be the maximal ideal of the polynomial ring $\mathbb{C}[\mathbf{k}_{H, i}]$ (cf. Subsection \ref{subsection00.3}) given by  the kernel of the canonical morphism from $\mathbb{C}[\mathbf{k}_{H, i}]$ to $\mathbb{C}$ defined by $\mathbf{k}_{H, i} \mapsto k_{H, i}$, and denote by $R$ the completion of $\mathbb{C}[\mathbf{k}_{H, i}]$ at $\mathfrak{m}$. In Section 5.3 of \cite{GGOR}, the KZ-functor is constructed over $R$, but the same construction works over $R'$ (cf. Section 6.1.2 of \cite{RSVV}). There are therefore functors $\textnormal{KZ}_{R'}: \mathcal{O}^{R', \Delta}_{\boldsymbol{c}} (\mathbb{Z}/n, \mathbb{C}) \longrightarrow \mathcal{H}^{R'}_{\boldsymbol{c}} (\mathbb{Z}/n, \mathbb{C})$-mod (here $\mathcal{O}^{R', \Delta}_{\boldsymbol{c}} (\mathbb{Z}/n, \mathbb{C})$ denotes the subcategory of category $\mathcal{O}^{R'}_{\boldsymbol{c}} (\mathbb{Z}/n, \mathbb{C})$ consisting of those objects that have a $\Delta$-filtration) and $\textnormal{KZ}_{K'}: \mathcal{O}^{K'}_{\boldsymbol{c}} (\mathbb{Z}/n, \mathbb{C}) \longrightarrow \mathcal{H}^{K'}_{\boldsymbol{c}} (\mathbb{Z}/n, \mathbb{C})$-mod and corresponding bimodules $P^{R'}_{\textnormal{KZ}}$ and $P^{K'}_{\textnormal{KZ}}$ as well as algebra homomorphisms  $\phi_{R'}: \mathcal{H}^{R'}_{\boldsymbol{c}} (\mathbb{Z}/n, \mathbb{C})  \longrightarrow \textnormal{End}_{H^{R'}_{\boldsymbol{c}}(\mathbb{Z}/n, \mathbb{C})} (P^{R'}_{\textnormal{KZ}})^{\textnormal{opp}}$ and $\phi_{K'}: \mathcal{H}^{K'}_{\boldsymbol{c}} (\mathbb{Z}/n, \mathbb{C})  \longrightarrow \textnormal{End}_{H^{K'}_{\boldsymbol{c}}(\mathbb{Z}/n, \mathbb{C})} (P^{K'}_{\textnormal{KZ}})^{\textnormal{opp}}$, which as $\mathcal{H}^{R'}_{\boldsymbol{c}} (\mathbb{Z}/n, \mathbb{C})$ and $\mathcal{H}^{K'}_{\boldsymbol{c}} (\mathbb{Z}/n, \mathbb{C}) $ are commutative by \eqref{equation6.01}, can be considered as homomorphisms into $\textnormal{End}_{H^{R'}_{\boldsymbol{c}}(\mathbb{Z}/n, \mathbb{C})} (P^{R'}_{\textnormal{KZ}})$ and $\textnormal{End}_{H^{K'}_{\boldsymbol{c}}(\mathbb{Z}/n, \mathbb{C})} (P^{K'}_{\textnormal{KZ}})$ respectively. In order to relate the maps $\phi$, $\phi_{R'}$ and $\phi_{K'}$ or equivalently the modules $P_{\text{KZ}}$, $P^{R'}_{\text{KZ}}$ and $P^{K'}_{\text{KZ}}$, we need the following lemma: 
\begin{theorem6.05} \label{theorem6.05}
For every projective module $P$ and every $\Delta$-filtered module $M$ in category $\mathcal{O}^{R'}_{\boldsymbol{c}} (\mathbb{Z}/n, \mathbb{C})$, the homomorphisms
\begin{equation*}
\theta_{P, M}: (R'/\mathfrak{m}') \otimes_{R'} \textnormal{Hom}_{\mathcal{O}^{R'}_{\boldsymbol{c}} (\mathbb{Z}/n, \mathbb{C})} (P, M) \to  \textnormal{Hom}_{\mathcal{O}_{\boldsymbol{c}} (\mathbb{Z}/n, \mathbb{C})} ((R'/\mathfrak{m}') \otimes_{R'} P, (R'/\mathfrak{m}') \otimes_{R'} M)
\end{equation*}
and 
\begin{equation*} 
\tau_{P, M}: K' \otimes_{R'} \textnormal{Hom}_{\mathcal{O}^{R'}_{\boldsymbol{c}} (\mathbb{Z}/n, \mathbb{C})} (P, M) \to  \textnormal{Hom}_{\mathcal{O}^{K'}_{\boldsymbol{c}} (\mathbb{Z}/n, \mathbb{C})} (K' \otimes_{R'} P, K' \otimes_{R'} M),
\end{equation*}
defined by $\theta_{P, M}(\overline{r} \otimes f)(\overline{s} \otimes p) = (\overline{rs}\otimes f(p))$ and $\tau_{P, M}(k \otimes f)(k' \otimes p) = (kk' \otimes f(p))$ respectively, are isomorphisms. 
\end{theorem6.05}
\begin{proof}
We start with $\theta_{P, M}$. Observe that there is a surjective map of $R'$-algebras $\pi: H^{R'}_{\boldsymbol{c}} (\mathbb{Z}/n, \mathbb{C}) \to (R'/\mathfrak{m}') \otimes_{R'} H^{R'}_{\boldsymbol{c}} (\mathbb{Z}/n, \mathbb{C}) \cong H_{\boldsymbol{c}} (\mathbb{Z}/n, \mathbb{C})$
given by $\pi: h \mapsto 1 \otimes h$. As a result, for every $H^{R'}_{\boldsymbol{c}} (\mathbb{Z}/n, \mathbb{C})$-module $N$, the module $(R'/\mathfrak{m}') \otimes_{R'} N$ is also an $H^{R'}_{\boldsymbol{c}} (\mathbb{Z}/n, \mathbb{C})$-module and there is a natural epimorphism $\pi_N: N \to (R'/\mathfrak{m}') \otimes_{R'} N$ given by $\pi_N: n \mapsto 1 \otimes n$. Furthermore, every $H_{\boldsymbol{c}} (\mathbb{Z}/n, \mathbb{C})$-module homomorphism between two modules $(R'/\mathfrak{m}') \otimes_{R'} N$ and $(R'/\mathfrak{m}') \otimes_{R'} N'$ is in addition an $H^{R'}_{\boldsymbol{c}} (\mathbb{Z}/n, \mathbb{C})$-module homomorphism and vice versa. Let $\psi \in \textnormal{Hom}_{\mathcal{O}_{\boldsymbol{c}} (\mathbb{Z}/n, \mathbb{C})} ((R'/\mathfrak{m}') \otimes_{R'} P, (R'/\mathfrak{m}') \otimes_{R'} M)$. We have a commutative diagram 
\begin{equation}  \label{equation6.200}
	\xymatrix@C2pc@R+0.5pc{
	P \ar^{\phi}[r]  \ar_{\pi_P}[d]    & M   \ar_{\pi_M}[d]  	\\
	(R'/\mathfrak{m}') \otimes_{R'} P \ar^{\psi}[r] & (R'/\mathfrak{m}') \otimes_{R'} M}
\end{equation}
where $\phi$ is induced by $\psi \circ \pi_P$ as $\pi_M$ is onto, $P$ is projective and all the modules live in category $\mathcal{O}^{R'}_{\boldsymbol{c}} (\mathbb{Z}/n, \mathbb{C})$. It follows that $\theta_{P, M}(1 \otimes \phi) = \psi$ and hence that $\theta_{P, M}$ is surjective. Next, by Theorem 4.24 and the remark before Corollary 4.26 of \cite{Ari}, we have that 
\begin{equation} \label{equation6.201}
(R'/\mathfrak{m}') \otimes_{R'} \textnormal{Hom}_{\mathcal{O}^{R'}_{\boldsymbol{c}} (\mathbb{Z}/n, \mathbb{C})} (P, M) \cong  \textnormal{Hom}_{\mathcal{O}^{R'}_{\boldsymbol{c}} (\mathbb{Z}/n, \mathbb{C})} (P, (R'/\mathfrak{m}') \otimes_{R'} M). 
\end{equation}
As $P$ is projective, it is $\Delta$-filtered and therefore free as an $R'$-module. It follows that every map of $\textnormal{Hom}_{\mathcal{O}^{R'}_{\boldsymbol{c}} (\mathbb{Z}/n, \mathbb{C})} (P, (R'/\mathfrak{m}') \otimes_{R'} M)$ factors through $(R'/\mathfrak{m}') \otimes_{R'} P$, and hence \eqref{equation6.201} gives 
\begin{equation*} \label{equation6.202}
(R'/\mathfrak{m}') \otimes_{R'} \textnormal{Hom}_{\mathcal{O}^{R'}_{\boldsymbol{c}} (\mathbb{Z}/n, \mathbb{C})} (P, M) \cong  \textnormal{Hom}_{\mathcal{O}_{\boldsymbol{c}} (\mathbb{Z}/n, \mathbb{C})} ((R'/\mathfrak{m}') \otimes_{R'} P, (R'/\mathfrak{m}') \otimes_{R'} M). 
\end{equation*}
As these spaces are finite-dimensional and $\theta_{P, M}$ is surjective, it follows that it is an isomorphism. As for $\tau_{P, M}$, since $P$ is finitely generated as an $H^{R'}_{\boldsymbol{c}} (\mathbb{Z}/n, \mathbb{C})$-module, if $\psi \in \textnormal{Hom}_{\mathcal{O}^{K'}_{\boldsymbol{c}} (\mathbb{Z}/n, \mathbb{C})} (K' \otimes_{R'} P, K' \otimes_{R'} M)$, then there exists $r \in R'$ such that $(r\psi)(P) \subseteq M$ (as $P$ and $M$ are free as $R'$-modules, they embed into $K' \otimes_{R'} P$ and $K' \otimes_{R'} M$ respectively), and it is straight-forward to check that the map that takes $\psi$ to $(1/r) \otimes (r\psi)$ is independent of $r$ and an inverse to $\tau_{P, M}$ which completes the proof. 
\end{proof}
Now, by Theorem 5.13 of \cite{GGOR} (with $R'$ in place of $R$), there is a commutative diagram 
\begin{equation*}  \label{equation6.02}
	\xymatrix@C+2pc@R+0.5pc{
	\mathcal{O}^{K'}_{\boldsymbol{c}} (\mathbb{Z}/n, \mathbb{C}) \ar^{\textnormal{KZ}_{K'}}[r]  & \mathcal{H}^{K'}_{\boldsymbol{c}} (\mathbb{Z}/n, \mathbb{C})  \textnormal{-mod}  	\\
	\mathcal{O}^{R', \Delta}_{\boldsymbol{c}} (\mathbb{Z}/n, \mathbb{C}) \ar^{K' \otimes_{R'} - }[u] \ar_{(R'/ \mathfrak{m}') \otimes_{R'} - }[d]  \ar^{\textnormal{KZ}_{R'}}[r] & \mathcal{H}^{R'}_{\boldsymbol{c}} (\mathbb{Z}/n, \mathbb{C})  \textnormal{-mod}  \ar_{K' \otimes_{R'} - }[u]  \ar^{(R'/ \mathfrak{m}') \otimes_{R'} - }[d] \\
	\mathcal{O}_{\boldsymbol{c}} (\mathbb{Z}/n, \mathbb{C}) \ar^{\textnormal{KZ}}[r] & \mathcal{H}_{\boldsymbol{c}} (\mathbb{Z}/n, \mathbb{C})  \textnormal{-mod}. } 
\end{equation*}
This diagram implies that, for $M$ an object of category $\mathcal{O}^{R', \Delta}_{\boldsymbol{c}} (\mathbb{Z}/n, \mathbb{C})$,  
\begin{IEEEeqnarray}{rCl} 
	 \textnormal{KZ}_{K'} (K' \otimes_{R'} M) & \cong & K' \otimes_{R'} \textnormal{KZ}_{R'} (M)  \cong  K' \otimes_{R'} \textnormal{Hom}_{\mathcal{O}^{R'}_{\boldsymbol{c}} (\mathbb{Z}/n, \mathbb{C})} (P^{R'}_{\textnormal{KZ}}, M)   \label{equation6.06}  \nonumber \\
	& \cong & \textnormal{Hom}_{\mathcal{O}^{K'}_{\boldsymbol{c}} (\mathbb{Z}/n, \mathbb{C})} (K' \otimes_{R'} P^{R'}_{\textnormal{KZ}}, K' \otimes_{R'} M)   \label{equation6.07} 
\end{IEEEeqnarray} 
as $\mathcal{H}^{K'}_{\boldsymbol{c}}(W, \mathfrak{h})$-modules, and that 
\begin{IEEEeqnarray}{rCl} 
	 \textnormal{KZ} ((R'/\mathfrak{m}') \otimes_{R'} M) & \cong & (R'/\mathfrak{m}') \otimes_{R'} \textnormal{KZ}_{R'} (M)  \cong  (R'/\mathfrak{m}') \otimes_{R'} \textnormal{Hom}_{\mathcal{O}^{R'}_{\boldsymbol{c}} (\mathbb{Z}/n, \mathbb{C})} (P^{R'}_{\textnormal{KZ}}, M)   \label{equation6.08}  \nonumber \\
	& \cong & \textnormal{Hom}_{\mathcal{O}_{\boldsymbol{c}} (\mathbb{Z}/n, \mathbb{C})} ((R'/\mathfrak{m}') \otimes_{R'} P^{R'}_{\textnormal{KZ}}, (R'/\mathfrak{m}') \otimes_{R'} M)   \label{equation6.09} 
\end{IEEEeqnarray} 
as $\mathcal{H}_{\boldsymbol{c}}(W, \mathfrak{h})$-modules, where \eqref{equation6.07} and \eqref{equation6.09} follow from Lemma \ref{theorem6.05} (observe that since $P^{R'}_{\text{KZ}}$ is projective in $\mathcal{O}^{R', \Delta}_{\boldsymbol{c}} (\mathbb{Z}/n, \mathbb{C})$, it is projective in $\mathcal{O}^{R'}_{\boldsymbol{c}} (\mathbb{Z}/n, \mathbb{C})$). As $P_{\textnormal{KZ}} \cong \Delta_{\boldsymbol{c}}(S(\mathfrak{h})_W) \cong (R'/ \mathfrak{m}') \otimes_{R'} \Delta_{\boldsymbol{c}}(R' \otimes_{\mathbb{C}} S(\mathfrak{h})_W)$ is an object of $(R'/ \mathfrak{m}') \otimes_{R'} \mathcal{O}^{R', \Delta}_{\boldsymbol{c}} (\mathbb{Z}/n, \mathbb{C})$, it follows by the Yoneda lemma that as an $(H^{K'}_{\boldsymbol{c}}(W, \mathfrak{h}), \mathcal{H}^{K'}_{\boldsymbol{c}}(W, \mathfrak{h}))$-bimodule, $P^{K'}_{\textnormal{KZ}} \cong K' \otimes_{R'} P^{R'}_{\textnormal{KZ}}$, and that as an $(H_{\boldsymbol{c}}(W, \mathfrak{h}), \mathcal{H}_{\boldsymbol{c}}(W, \mathfrak{h}))$-bimodule, $P_{\textnormal{KZ}} \cong (R'/ \mathfrak{m}') \otimes_{R'} P^{R'}_{\textnormal{KZ}}$. This is equivalent to that
\begin{equation*} \label{equation 6.03}
\phi_{K'} = K' \otimes_{R'} \phi_{R'} \quad \textnormal{and} \quad  \phi = (R'/ \mathfrak{m}') \otimes_{R'} \phi_{R'}. 
\end{equation*}
A description of the map $\phi_{R'}$ therefore immediately gives a description of the map $\phi$. In order to describe $\phi_{R'}$, we will construct a map $\psi_{R'}: \mathcal{H}^{R'}_{\boldsymbol{c}} (\mathbb{Z}/n, \mathbb{C})  \longrightarrow \textnormal{End}_{H^{R'}_{\boldsymbol{c}}(\mathbb{Z}/n, \mathbb{C})} (P^{R'}_{\textnormal{KZ}})$, and verify that $K' \otimes_{R'} \psi_{R'} = K' \otimes_{R'} \phi_{R'}$, thus implying that $\psi_{R'} = \phi_{R'}$. The reason for this approach is that category $\mathcal{O}^{K'}_{\boldsymbol{c}} (\mathbb{Z}/n, \mathbb{C})$ is semisimple, which will be key to showing that $K' \otimes_{R'} \psi_{R'} = K' \otimes_{R'} \phi_{R'}$.  In order to construct the map $\psi_{R'}$, we first need a description of $P^{R'}_{\textnormal{KZ}}$ and $P^{K'}_{\textnormal{KZ}}$. For $\boldsymbol{c} \in \mathcal{F}$, Theorem B implies an analogous description to that over $\mathbb{C}$, in terms of the coinvariant algebra: 

\begin{theorem6.01}  \label{theorem6.01}
For $\boldsymbol{c} \in \mathcal{F}$, we have that $P^{R'}_{\textnormal{KZ}} \cong \Delta_{\boldsymbol{c}}(R' \otimes_{\mathbb{C}} S(\mathfrak{h})_W)$ as $H^{R'}_{\boldsymbol{c}}(\mathbb{Z}/n, \mathbb{C})$-modules and $P^{K'}_{\textnormal{KZ}} \cong \Delta_{\boldsymbol{c}}(K' \otimes_{\mathbb{C}} S(\mathfrak{h})_W)$ as $H^{K'}_{\boldsymbol{c}}(\mathbb{Z}/n, \mathbb{C})$-modules. 
\end{theorem6.01}

\begin{proof}
By Theorem B, $(R'/\mathfrak{m}') \otimes_{R'} P_{\textnormal{KZ}}^{R'} \cong (R'/\mathfrak{m}') \otimes_{R'} \Delta(R' \otimes_{\mathbb{C}} S(\mathfrak{h})_W)$ as $H_{\boldsymbol{c}} (\mathbb{Z}/n, \mathbb{C})$-modules for $\boldsymbol{c} \in \mathcal{F}$. It therefore follows from a diagram such as in \eqref{equation6.200} that there is a map 
\begin{equation*}
\Omega: P_{\textnormal{KZ}}^{R'} \rightarrow \Delta(R' \otimes_{\mathbb{C}} S(\mathfrak{h})_W)
\end{equation*}
such that $(R'/\mathfrak{m}') \otimes_{R'} \Omega$ is an isomorphism. We want to show that $\Omega$ is an isomorphism as well. By Corollary 2.8 of \cite{GGOR}, category $\mathcal{O}^{R'}_{\boldsymbol{c}} (\mathbb{Z}/n, \mathbb{C})$ has a progenerator $Q^{R'}$ and with $\Gamma^{R'} = \textnormal{End}_{H^{R'}_{\boldsymbol{c}} (\mathbb{Z}/n, \mathbb{C})} (Q^{R'})^{\textnormal{opp}}$ there is the standard equivalence 
\begin{equation*}
F^{R'}(-) := \textnormal{Hom}_{H^{R'}_{\boldsymbol{c}} (\mathbb{Z}/n, \mathbb{C})} (Q^{R'}, -): \mathcal{O}^{R'}_{\boldsymbol{c}} (\mathbb{Z}/n, \mathbb{C}) \stackrel{\sim}{\longrightarrow} \Gamma^{R'}\textnormal{-mod}.
\end{equation*}
Then $Q := (R'/\mathfrak{m}') \otimes_{R'} Q^{R'}$ is a progenerator for category $\mathcal{O}_{\boldsymbol{c}} (\mathbb{Z}/n, \mathbb{C})$, and with $\Gamma := \textnormal{End}_{H_{\boldsymbol{c}} (\mathbb{Z}/n, \mathbb{C})} (Q)^{\textnormal{opp}}$ there is the analogous equivalence
\begin{equation*}
F(-) := \textnormal{Hom}_{H_{\boldsymbol{c}} (\mathbb{Z}/n, \mathbb{C})} (Q, -): \mathcal{O}_{\boldsymbol{c}} (\mathbb{Z}/n, \mathbb{C}) \stackrel{\sim}{\longrightarrow} \Gamma \textnormal{-mod}.
\end{equation*}
Since $F^{R'}$ is an equivalence, it is enough to prove that $F^{R'}(\Omega)$ is an isomorphism, and as $\Gamma^{R'}$ is a finite-rank $R'$-module, $F^{R'}(P_{\textnormal{KZ}}^{R'})$ and $F^{R'}(\Delta(R' \otimes_{\mathbb{C}} S(\mathfrak{h})_W))$ are finitely generated as $R'$-modules, and so it is sufficient by Nakayama's Lemma to prove that $(R'/\mathfrak{m}') \otimes_{R'} F^{R'}(\Omega)$ is an isomorphism. It follows from Lemma \ref{theorem6.05} that we have a commutative diagram
\begin{equation*}  \label{equation14}
	\xymatrix@C8pc@R+0.5pc{
	(R'/\mathfrak{m}') \otimes_{R'} F^{R'}(P_{\textnormal{KZ}}^{R'}) \ar^{(R'/\mathfrak{m}') \otimes_{R'} F^{R'}(\Omega)}[r]  \ar_{\cong}[d]    & (R'/\mathfrak{m}') \otimes_{R'} F^{R'}(\Delta(R' \otimes_{\mathbb{C}} S(\mathfrak{h})_W))   \ar_{\cong}[d]  	\\
	F((R'/\mathfrak{m}') \otimes_{R'} P_{\textnormal{KZ}}^{R'}) \ar^{\cong}[r] &  F((R'/\mathfrak{m}') \otimes_{R'} \Delta(R' \otimes_{\mathbb{C}} S(\mathfrak{h})_W))  }
\end{equation*}
and this completes the proof that $P^{R'}_{\textnormal{KZ}} \cong \Delta_{\boldsymbol{c}}(R' \otimes_{\mathbb{C}} S(\mathfrak{h})_W)$ when $\boldsymbol{c} \in \mathcal{F}$. As $P^{K'}_{\textnormal{KZ}} = K' \otimes_{R'} P^{R'}_{\textnormal{KZ}}$ and $\Delta_{\boldsymbol{c}}(K' \otimes_{\mathbb{C}} S(\mathfrak{h})_W) = K' \otimes_{R'} \Delta_{\boldsymbol{c}}(R' \otimes_{\mathbb{C}} S(\mathfrak{h})_W)$, this concludes the proof of the claim. 
\end{proof}
We proceed to construct the map $\psi_{R'}$. 
\subsection{The construction of the map $\psi_{R'}$} \label{subsection6.02} 
In this subsection, we construct the homomorphism $\psi_{R'}: \mathcal{H}^{R'}_{\boldsymbol{c}} (\mathbb{Z}/n, \mathbb{C})  \stackrel{\sim}{\longrightarrow} \textnormal{End}_{H^{R'}_{\boldsymbol{c}}(\mathbb{Z}/n, \mathbb{C})} (\Delta_{\boldsymbol{c}}(R' \otimes_{\mathbb{C}} S(\mathfrak{h})_W))$. This is done by identifying an element of $\textnormal{End}_{H^{R'}_{\boldsymbol{c}}(\mathbb{Z}/n, \mathbb{C})} (\Delta_{\boldsymbol{c}}(R' \otimes_{\mathbb{C}} S(\mathfrak{h})_W))$ that is annihilated by the polynomial $\prod^n_{j = 1} (T - q^{-j}.{q'}_j^{-1})$ giving the defining relation of $\mathcal{H}^{R'}_{\boldsymbol{c}} (\mathbb{Z}/n, \mathbb{C})$ in \eqref{equation6.01}. To this end, we recall from Subsection \ref{subsection00.3} the element $\mathbf{eu} =  \mathbf{eu}_{\mathbb{C}} - z$ of  $H^{R'}_{\boldsymbol{c}} (\mathbb{Z}/n, \mathbb{C})$ which gives it its inner grading, and that $\mathbf{eu} =  x \xi - \sum_{i =1}^{n-1} c^{R'}_i \varepsilon_i = x \xi - \sum_{i =1}^{n-1} z_i \varepsilon_i $. The element $\mathbf{eu}$ acts on $\Delta_{\boldsymbol{c}}(R' \otimes_{\mathbb{C}} S(\mathfrak{h})_W)$ by multiplication, and since it sits in the degree zero component of $H^{R'}_{\boldsymbol{c}} (\mathbb{Z}/n, \mathbb{C})$, it preserves the graded components of $\Delta_{\boldsymbol{c}}(R' \otimes_{\mathbb{C}} S(\mathfrak{h})_W)$. Proposition \ref{theorem7.3} easily generalises (with $ \{v_{i,j} =  x^i \otimes \bar{\xi}^j : i \geq 0, 0 \leq j \leq n-1\}$ an $R'$-basis rather than a $\mathbb{C}$-basis) from $\Delta_{\boldsymbol{c}}(S(\mathfrak{h})_W)$ to $\Delta_{\boldsymbol{c}}(R' \otimes_{\mathbb{C}} S(\mathfrak{h})_W)$ giving the following: 

\begin{theorem6.02}  \label{theorem6.02}

\begin{itemize}
\item[(i)]
For $k \geq 0$, the action of $\mathbf{eu}$ on $\Delta_{\boldsymbol{c}}(R' \otimes_{\mathbb{C}} S(\mathfrak{h})_W)_k$ with respect to the $R'$-basis $\{v_{k + (n-1), n-1}, v_{k + (n-2), n-2}, \dots, v_{k, 0}  \}$ is given by the $n \times n$-matrix 
\begin{equation*}  \label{equation6.04}
	\begin{bmatrix}
	k + (n-1) - z_1 & 1 & 0 & \ldots & 0 & 0 \\
	0 & k + (n-2) - z_2 & 1 & \ldots & 0 & 0 \\
	0 & 0 & k + (n-3) - z_3 & \ldots & 0 & 0 \\
	\vdots & \vdots & \vdots & \ddots  & \vdots & \vdots \\
	0 & 0 & 0 & \ldots & k + 1 - z_{n-1} & 1\\
	0 & 0 & 0 & \ldots & 0 & k - z_n 
	\end{bmatrix}.
\end{equation*}
\item[(ii)]
For $1-n \leq k < 0$, the action of $\mathbf{eu}$ on $\Delta_{\boldsymbol{c}}(R' \otimes_{\mathbb{C}} S(\mathfrak{h})_W)_k$ with respect to the $R'$-basis $\{v_{k + (n-1), n-1}, v_{k + (n-2), n-2}, \dots, v_{0, -k}  \}$ is given by the $(n+k) \times (n+k)$-matrix 
\begin{equation*}  \label{equation6.05}
	\begin{bmatrix}
	k + (n-1) - z_1 & 1 & 0 & \ldots & 0 & 0 \\
	0 & k + (n-2) - z_2 & 1 & \ldots & 0 & 0 \\
	0 & 0 & k + (n-3) - z_3 & \ldots & 0 & 0 \\
	\vdots & \vdots & \vdots & \ddots  & \vdots & \vdots \\
	0 & 0 & 0 & \ldots & 1 - z_{n + k-1} & 1\\
	0 & 0 & 0 & \ldots & 0 &  - z_{n+k} 
	\end{bmatrix}.
\end{equation*} 
\end{itemize}
\end{theorem6.02}
\begin{proof}
Proposition \ref{theorem7.3}, applied over $R'$, implies that, using the notation of Subsection \ref{subsection00.3}, $R'.v_{i,j} \cong R' \otimes_{\mathbb{C}} E_{i-j}$ as a $\mathbb{Z}/n$-module. Therefore, as a $\mathbb{Z}/n$-module, $\Delta_{\boldsymbol{c}}(R' \otimes_{\mathbb{C}} S(\mathfrak{h})_W)_k$ is a direct sum of $R' \otimes_{\mathbb{C}} E_k$'s. Since $E_k$ is fixed by $\varepsilon_k$ and killed by all the other $\varepsilon_l$, it follows that $\sum_{i =1}^{n-1} z_i \varepsilon_i $ acts on $\Delta_{\boldsymbol{c}}(R' \otimes_{\mathbb{C}} S(\mathfrak{h})_W)_k$ by multiplication by $z_k$. The claim is then a direct consequence of equation \eqref{equation7.07} (recalling that $c^{R'}_k = z_k$). 
\end{proof}

We now introduce the element of $\textnormal{End}_{H^{R'}_{\boldsymbol{c}} (\mathbb{Z}/n, \mathbb{C})} (\Delta_{\boldsymbol{c}}(R' \otimes_{\mathbb{C}} S(\mathfrak{h})_W))$ which we will use to construct the map $\psi_{R'}$. We let $\eta := s.\textnormal{exp}(\frac{2\pi i}{n} \mathbf{eu})$. This is not well-defined as an element of $H^{R'}_{\boldsymbol{c}} (\mathbb{Z}/n, \mathbb{C})$, but we will see that it is well-defined as an element of $\textnormal{End}_{H^{R'}_{\boldsymbol{c}} (\mathbb{Z}/n, \mathbb{C})} (\Delta_{\boldsymbol{c}}(R' \otimes_{\mathbb{C}} S(\mathfrak{h})_W))$ and that it is annihilated by $\prod^n_{j = 1} (T - q^{-j}.{q'}_j^{-1})$. 
\begin{theorem6.03}    \label{theorem6.03}
\begin{itemize}
\item[(i)]
When considered as an element of $\textnormal{End}_{R'} (\Delta_{\boldsymbol{c}}(R' \otimes_{\mathbb{C}} S(\mathfrak{h})_W))$, $\eta := s.\textnormal{exp}(\frac{2\pi i}{n} \mathbf{eu})$ is well-defined. 
\item[(ii)] 
The element $\eta$ is annihilated by $\prod^n_{j = 1} (T - q^{-j}.{q'}_j^{-1})$ where $q = \textnormal{exp}(\frac{2\pi i}{n})$ and $q'_j = \textnormal{exp}(\frac{2\pi i}{n} z_j)$. 
\item[(iii)]
$\eta$ is an element of $\textnormal{End}_{H^{R'}_{\boldsymbol{c}} (\mathbb{Z}/n, \mathbb{C})} (\Delta_{\boldsymbol{c}}(R' \otimes_{\mathbb{C}} S(\mathfrak{h})_W))$. 
\end{itemize}
\end{theorem6.03}
\begin{proof}
In light of Proposition \ref{theorem6.02}, we need to define the exponential of a square matrix with entries in $R'$. Suppose that $M = ((f_{ij}, U_{ij})) \in \textnormal{Mat}_{l \times l}(R')$ where $U_{ij}$ is some open set containing $(c_1, \ldots, c_{n-1})$ and $f_{ij} : U_{ij} \to \mathbb{C}$ is holomorphic. Let $W$ be an open set containing $(c_1, \ldots, c_{n-1})$ such that $\overline{W}$ is compact and $\overline{W} \subset \bigcap U_{i,j}$. Restricting the $f_{i,j}$ to $\overline{W}$, we can think of them as elements of $B = C(\overline{W}, \mathbb{C})$, the Banach algebra of continuous functions from $\overline{W}$ to $\mathbb{C}$. $M$ can then be considered as an element of the Banach algebra $\textnormal{Mat}_{l \times l} (B)$. As exponentials are well defined in the setting of Banach algebras, we can define $N = \exp (M) = (g_{ij})$ as an element of $\textnormal{Mat}_{l \times l} (B)$. Since $\overline{W}$ is compact, the convergence of every entry of the matrix is uniform. In particular, the convergence is uniform on $W$. As the uniform limit of holomorphic functions is again holomorphic, each $(g_{ij}, W) \in R'$, so we can let $\textnormal{exp}(M) = ((g_{ij}, W)) \in \textnormal{Mat}_{l \times l}(R')$ and it is clear that this definition is independent of the particular representatives $(f_{ij}, U_{ij})$ above. We can therefore take the exponentials of the matrices of Proposition \ref{theorem6.02} to see that $s.\textnormal{exp}(\frac{2\pi i}{n} \mathbf{eu}|_{\Delta_{\boldsymbol{c}}(R' \otimes_{\mathbb{C}} S(\mathfrak{h})_W)_k})$ is well-defined for all integers $k$, and therefore that $\eta := s.\textnormal{exp}(\frac{2\pi i}{n} \mathbf{eu})$ is well-defined as an element of $\textnormal{End}_{R'} (\Delta_{\boldsymbol{c}}(R' \otimes_{\mathbb{C}} S(\mathfrak{h})_W))$. 

With respect to the $R'$-basis of $\Delta_{\boldsymbol{c}}(R' \otimes_{\mathbb{C}} S(\mathfrak{h})_W)_k$ specified in Proposition \ref{theorem6.02}, the matrix of $\textnormal{exp}(\frac{2\pi i}{n} \mathbf{eu})|_{\Delta_{\boldsymbol{c}}(R' \otimes_{\mathbb{C}} S(\mathfrak{h})_W)_k}$ is upper triangular, with $j$th diagonal entry equal to $\textnormal{exp} (\frac{2\pi i}{n} (k + (n - j) - z_j)) = q^{k-j} {q'}_j^{-1}$. 
As $\Delta_{\boldsymbol{c}}(R' \otimes_{\mathbb{C}} S(\mathfrak{h})_W)_k$ is isomorphic to a direct sum of $R' \otimes_{\mathbb{C}}E_k$'s as a $\mathbb{Z}/n$-module, it follows that $s$ acts on $\Delta_{\boldsymbol{c}}(R' \otimes_{\mathbb{C}} S(\mathfrak{h})_W)_k$ by multiplication by $q^{-k}$, and therefore, the matrix of $\eta|_{\Delta_{\boldsymbol{c}}(R' \otimes_{\mathbb{C}} S(\mathfrak{h})_W)_k}$ is upper triangular, with $j$th diagonal entry equal to $q^{-k}.q^{k-j} {q'}_j^{-1} = q^{-j} {q'}_j^{-1}$. It therefore follows that $\eta|_{\Delta_{\boldsymbol{c}}(S(\mathfrak{h})_W)_k}$ is annihilated by $\prod^n_{j = 1} (T - q^{-j}.{q'}_j^{-1})$ for all $k$, so that $\eta$ itself is annihilated by $\prod^n_{j = 1} (T - q^{-j}.{q'}_j^{-1})$.  

Finally, since ${H^{R'}_{\boldsymbol{c}} (\mathbb{Z}/n, \mathbb{C})}$ is generated as an algebra by $x$, $\xi$ and $s$, it is enough to prove that $\eta$ commutes with the action of these three elements on $\Delta_{\boldsymbol{c}}(R' \otimes_{\mathbb{C}} S(\mathfrak{h})_W)$. Recall that $\mathbf{eu}$ controls the grading of ${H^{R'}_{\boldsymbol{c}} (\mathbb{Z}/n, \mathbb{C})}$ in the sense that $[\mathbf{eu}, a] = \textnormal{deg}(a)a$ for all homogenous elements $a$ of ${H^{R'}_{\boldsymbol{c}} (\mathbb{Z}/n, \mathbb{C})}$. It follows by induction that $\mathbf{eu}^ja = a(\mathbf{eu} + \textnormal{deg}(a))^j$ for all positive integers $j$, and this
also holds with $\mathbf{eu}$ and $a$ viewed as  elements of  $\textnormal{End}_{R'} (\Delta_{\boldsymbol{c}}(R' \otimes_{\mathbb{C}} S(\mathfrak{h})_W))$. A routine calculation then shows that $\textnormal{exp} (\frac{2\pi i}{n} \mathbf{eu}). a =  \textnormal{exp} (\frac{2\pi i}{n}\textnormal{deg}(a)). [a. \textnormal{exp} (\frac{2\pi i}{n} \mathbf{eu} )]$ for all such $a$. The fact that $s.x = q^{-1}x$ and $s.\xi = q\xi$ then implies, through another straight-forward calculation, that $\eta$ commutes with the action of $x$, $s$ and $\xi$. Therefore $\eta \in \textnormal{End}_{H^{R'}_{\boldsymbol{c}} (\mathbb{Z}/n, \mathbb{C})} (\Delta_{\boldsymbol{c}}(R' \otimes_{\mathbb{C}} S(\mathfrak{h})_W))$, which completes the proof.
\end{proof}

With the help of Proposition \ref{theorem6.03} we then define $\psi_{R'}$ as follows: 

\begin{definition6.01}  \label{definition6.01}
The map $\psi_{R'}: \mathcal{H}^{R'}_{\boldsymbol{c}} (\mathbb{Z}/n, \mathbb{C})  \longrightarrow \textnormal{End}_{H^{R'}_{\boldsymbol{c}}(\mathbb{Z}/n, \mathbb{C})} (\Delta_{\boldsymbol{c}}(R' \otimes_{\mathbb{C}} S(\mathfrak{h})_W))$ is defined by $\psi(\overline{T}) = \eta$, where we recall that $\mathcal{H}^{R'}_{\boldsymbol{c}} (\mathbb{Z}/n, \mathbb{C}) = R'[T]/ \langle \prod^n_{j = 1} (T - q^{-j}.{q'}_j^{-1}) \rangle$. 
\end{definition6.01}

Observe that $\eta$ can also be considered as an element of $\textnormal{End}_{H^{K'}_{\boldsymbol{c}}(\mathbb{Z}/n, \mathbb{C})} (\Delta_{\boldsymbol{c}}(K' \otimes_{\mathbb{C}} S(\mathfrak{h})_W))$ and $\textnormal{End}_{H_{\boldsymbol{c}}(\mathbb{Z}/n, \mathbb{C})} (\Delta_{\boldsymbol{c}}(S(\mathfrak{h})_W))$ so that we also have maps $\psi_{K'} = K' \otimes_{R'} \psi_{R'}$ and $\psi = (R'/ \mathfrak{m}') \otimes \psi_{R'}$ defined similarly to $\psi_{R'}$. 

\subsection{The equality of the maps $\phi$ and $\psi$} \label{subsection6.03} As explained above, we show that $\psi_{R'} = \phi_{R'}$ by showing that $K' \otimes_{R'} \psi_{R'} = K' \otimes_{R'} \phi_{R'}$. 

\begin{theorem6.04} \label{theorem6.04}
We have that $K' \otimes_{R'} \psi_{R'} = K' \otimes_{R'} \phi_{R'}$. As a consequence, $\psi_{R'} = \phi_{R'}$ and $\psi = \phi$. 
\end{theorem6.04}

\begin{proof}
By  \eqref{equation6.01} it is sufficient to prove that $\psi_{K'} (\overline{T}) = \phi_{K'} (\overline{T})$ where $\psi_{K'} (\overline{T}), \phi_{K'} (\overline{T}) \in \textnormal{End}_{H^{K'}_{\boldsymbol{c}}(\mathbb{Z}/n, \mathbb{C})} (\Delta_{\boldsymbol{c}}(K' \otimes_{\mathbb{C}} S(\mathfrak{h})_W))$. As $c^{K'}_i = z_i$, it follows from Corollary 2.20 of \cite{GGOR} that category $\mathcal{O}^{K'}_{\boldsymbol{c}} (\mathbb{Z}/n, \mathbb{C})$ is semisimple. It therefore follows from Corollary \ref{theorem00.8} (which easily generalises to $K'$) that $\Delta_{\boldsymbol{c}}(K' \otimes_{\mathbb{C}} S(\mathfrak{h})_W) = \bigoplus_{j=1}^n \Delta_{\boldsymbol{c}}(K' \otimes_{\mathbb{C}} E_j)$ with the $\Delta_{\boldsymbol{c}}(K' \otimes_{\mathbb{C}} E_j)$ simple. Therefore $\textnormal{End}_{H^{K'}_{\boldsymbol{c}}(\mathbb{Z}/n, \mathbb{C})} (\Delta_{\boldsymbol{c}}(K' \otimes_{\mathbb{C}} S(\mathfrak{h})_W)) = \bigoplus_{j=1}^n \textnormal{End}_{H^{K'}_{\boldsymbol{c}}(\mathbb{Z}/n, \mathbb{C})} (\Delta_{\boldsymbol{c}}(K' \otimes_{\mathbb{C}} E_j))$ and it is sufficient to show that $\psi_{K'} (\overline{T})|_{\Delta_{\boldsymbol{c}}(K' \otimes_{\mathbb{C}} E_j)} = \phi_{K'} (\overline{T})|_{\Delta_{\boldsymbol{c}}(K' \otimes_{\mathbb{C}} E_j)}$ for $1 \leq j \leq n$. As mentioned in Subsection \ref{subsection00.5}, $\Delta_{\boldsymbol{c}}(K' \otimes_{\mathbb{C}} E_j)_0 = \mathcal{W}_{-c^{K'}_j}(\Delta_{\boldsymbol{c}}(K' \otimes_{\mathbb{C}} E_j))$ (cf. \cite{Gua} and Lemma 2.3 of \cite{GGOR}) and $\Delta_{\boldsymbol{c}}(K' \otimes_{\mathbb{C}} E_j)_0 = K' \otimes_{\mathbb{C}} E_j$ as $\mathbb{Z}/n$-modules. Since $\textnormal{End}_{\mathbb{Z}/n} (K' \otimes_{\mathbb{C}} E_j) \cong K'$, it follows that every element of $\textnormal{End}_{H^{K'}_{\boldsymbol{c}}(\mathbb{Z}/n, \mathbb{C})} (\Delta_{\boldsymbol{c}}(K' \otimes_{\mathbb{C}} E_j))$ has an eigenvector so that $\textnormal{End}_{H^{K'}_{\boldsymbol{c}}(\mathbb{Z}/n, \mathbb{C})} (\Delta_{\boldsymbol{c}}(K' \otimes_{\mathbb{C}} E_j)) \cong K'$ as well, since $\Delta_{\boldsymbol{c}}(K' \otimes_{\mathbb{C}} E_j)$ is simple. The action of $\psi_{K'} (\overline{T})$ on $\Delta_{\boldsymbol{c}}(K' \otimes_{\mathbb{C}} E_j)$ is given by $\eta = s.\textnormal{exp}(\frac{2\pi i}{n} \mathbf{eu})$. As $\mathbf{eu}$ and $s$ act on $\Delta_{\boldsymbol{c}}(K' \otimes_{\mathbb{C}} E_j)_0 = \mathcal{W}_{-c^{K'}_j}(\Delta_{\boldsymbol{c}}(K' \otimes_{\mathbb{C}} E_j))$ by multiplication by $- z_j$ and $q^{-j} = \textnormal{exp} (-\frac{2\pi i}{n}j))$ respectively, it follows that $\eta$ acts on $\Delta_{\boldsymbol{c}}(K' \otimes_{\mathbb{C}} E_j)$ by multiplication by $q^{-j} {q'}_j^{-1} = q^{-j} \textnormal{exp} (-\frac{2\pi i}{n} z_j)$. Next, the action of $\phi_{K'} (\overline{T})$ on $\Delta_{\boldsymbol{c}}(K' \otimes_{\mathbb{C}} E_j)$ corresponds to the action of $\overline{T}$ on $\textnormal{KZ}_{K'} (\Delta_{\boldsymbol{c}}(K' \otimes_{\mathbb{C}} E_j)) \cong \textnormal{Hom}_{H^{K'}_{\boldsymbol{c}}(\mathbb{Z}/n, \mathbb{C})} (\Delta_{\boldsymbol{c}}(K' \otimes_{\mathbb{C}} S(\mathfrak{h})_W), \Delta_{\boldsymbol{c}}(K' \otimes_{\mathbb{C}} E_j)) \cong \textnormal{End}_{H^{K'}_{\boldsymbol{c}}(\mathbb{Z}/n, \mathbb{C})} (\Delta_{\boldsymbol{c}}(K' \otimes_{\mathbb{C}} E_j))$, and it follows from Theorem 4.12 of \cite{BMR} that this is again given by the constant $q^{-j} {q'}_j^{-1}$. Therefore $K' \otimes_{R'} \psi_{R'} = K' \otimes_{R'} \phi_{R'}$ which completes the proof. 
%
%The action of $\psi_{K'} (\overline{T})$ on $\Delta_{\boldsymbol{c}}(K' \otimes_{\mathbb{C}} E_j)$ is given by $\eta = s.\textnormal{exp}(\frac{2\pi i}{n} \mathbf{eu})$. As mentioned in Subsection \ref{subsection00.5}, $\Delta_{\boldsymbol{c}}(K' \otimes_{\mathbb{C}} E_j)_k = \mathcal{W}_{k - c^{K'}_j}(\Delta_{\boldsymbol{c}}(K' \otimes_{\mathbb{C}} E_j))$ (cf. \cite{Gua} and Lemma 2.3 of \cite{GGOR}), and therefore $\mathbf{eu}$ acts on the 1-dimensional F-vector space $\Delta_{\boldsymbol{c}}(K' \otimes_{\mathbb{C}} E_j)_k$ by multiplication by $k - z_j$. As $s$ acts on $\Delta_{\boldsymbol{c}}(K' \otimes_{\mathbb{C}} E_j)_k$ by multiplication by $q^{-(k+j)} = \textnormal{exp} (-\frac{2\pi i}{n} (k+j))$ it follows that $\eta$ acts on $\Delta_{\boldsymbol{c}}(K' \otimes_{\mathbb{C}} E_j)_k$ by multiplication by $q^{-j} {q'}_j^{-1} = q^{-j} \textnormal{exp} (-\frac{2\pi i}{n} z_j)$ for every $k \geq 0$ so that $\psi_{K'} (\overline{T})$ acts on all of $\Delta_{\boldsymbol{c}}(K' \otimes_{\mathbb{C}} E_j)$ by multiplication by $q^{-j} {q'}_j^{-1}$. 
\end{proof}

As a corollary, we obtain: 
\begin{theorem000.4}  \label{theorem000.4}
For $\boldsymbol{c} \in \mathcal{F}$, the isomorphism $\phi: \mathcal{H}_{\boldsymbol{c}} (\mathbb{Z}/n, \mathbb{C})  \stackrel{\sim}{\longrightarrow} \textnormal{End}_{H_{\boldsymbol{c}} (\mathbb{Z}/n, \mathbb{C})} (\Delta_{\boldsymbol{c}}(S(\mathfrak{h})_W))$ induced by the $\KZ$-functor is given by 
$\phi(\overline{T}) = \eta = s.\textnormal{exp}(\frac{2\pi i}{n} \mathbf{eu})$, where $\overline{T}$ is the natural generator of  $\mathcal{H}_{\boldsymbol{c}} (\mathbb{Z}/n, \mathbb{C}) = \mathbb{C}[T]/ \langle \prod^n_{j = 1} (T - q^{-j}.q_j^{-1}) \rangle$.  
\end{theorem000.4}

\end{document}